\documentclass[11pt]{amsart}
\usepackage{amssymb}
\usepackage{fullpage}

\usepackage[OT2,OT1]{fontenc}

\def\cyr{\fontencoding{OT2}\fontfamily{wncyr}\selectfont}

\def\Ch{\textrm{\cyr CH}}

\newcommand{\bmu}{{\bar\mu}}
\newcommand{\bnu}{{\bar\nu}}

\newcommand{\crit}{{\mathrm{cr}}}

\newcommand{\ntop}[2]{\genfrac{}{}{0pt}{1}{#1}{#2}}

\let\newpf\proof \let\proof\relax 
\newenvironment{pf}{\newpf[\proofname]}{\qed\endtrivlist}

\def\area{\operatorname {area}}

\def\MD{\operatorname{MD}}

\def\be{\begin{equation}}
\def\ee{\end{equation}}

\def\ba{{\begin{align}}}
\def\ea{{\end{align}}}

\def\u{{\mathbb U}}

\def\d{{\underline d}}

\def\op{\overline\partial}

\def\0{{\mathbf 0}}

\def\cal{\mathcal}

\def\SSS{{\cal {S}}}

\newtheorem{thm}{Theorem}[section]
\newtheorem*{thmA}{Theorem A}%[section]
\newtheorem*{thmB}{Theorem B}%[section]
\newtheorem*{thmC}{Theorem C}%[section]
\newtheorem*{thmD}{Theorem D}%[section]
\newtheorem*{thmE}{Theorem E}%[section]
\newtheorem{cor}[thm]{Corollary}
\newtheorem{conj}[thm]{Conjecture}
\newtheorem{lem}[thm]{Lemma}
\newtheorem{lemma}[thm]{Lemma}

\newtheorem{prop}[thm]{Proposition}

\theoremstyle{remark}
\newtheorem{rem}{Remark}[section]

\newtheorem{example}{Example}[section]
\newtheorem{problem}{Problem}

\numberwithin{equation}{section}

\def \bn {\hfill \\ \smallskip\noindent}

\theoremstyle{definition}

\def\proof{\bn {\bf Proof.} }

\newcommand{\di}{\partial}
\newcommand{\dibar}{\bar\partial}
\newcommand{\ra}{\rightarrow}

\def\ssk{\smallskip}
\def\msk{\medskip}

\def\nin{\noindent}

\def\sm{\smallsetminus}

\def\ssm{\smallsetminus}

\renewcommand{\setminus}{\ssm}

\newcommand{\diam}{\operatorname{diam}}
\newcommand{\depth}{\operatorname{depth}}
\newcommand{\dist}{\operatorname{dist}}

\newcommand{\cl}{\operatorname{cl}}
\newcommand{\inter}{\operatorname{int}}
\renewcommand{\mod}{\operatorname{mod}}
\newcommand{\tl}{\tilde}

\newcommand{\orb}{\operatorname{orb}}
\newcommand{\HD}{\operatorname{HD}}
\newcommand{\supp}{\operatorname{supp}}
\newcommand{\id}{\operatorname{id}}
\newcommand{\length}{\operatorname{length}}
\newcommand{\dens}{\operatorname{dens}}

\newcommand{\Spec}{\operatorname{Spec}}

\newcommand{\esssup}{\operatorname{ess-sup}}

\newcommand{\eps}{{\epsilon}}
\newcommand{\De}{{\Delta}}
\newcommand{\de}{{\delta}}
\newcommand{\la}{{\lambda}}
\newcommand{\La}{{\Lambda}}
\newcommand{\si}{{\sigma}}
\newcommand{\Om}{{\Omega}}
\newcommand{\om}{{\omega}}

\newcommand{\AAA}{{\cal A}}

\newcommand{\CC}{{\cal C}}

\newcommand{\FF}{{\cal F}}
\newcommand{\GG}{{\cal G}}
\newcommand{\JJ}{{\cal J}}
\newcommand{\HH}{{\cal H}}

\newcommand{\LL}{{\cal L}}
\newcommand{\MM}{{\cal M}}

\newcommand{\OO}{{\cal O}}

\newcommand{\VV}{{\cal V}}
\newcommand{\WW}{{\cal W}}
\newcommand{\XX}{{\cal X}}
\newcommand{\YY}{{\cal Y}}
\newcommand{\ZZ}{{\cal Z}}

\newcommand{\A}{{\mathbb A}}
\newcommand{\C}{{\mathbb C}}
\newcommand{\D}{{\mathbb D}}

\newcommand{\N}{{\mathbb N}}

\newcommand{\R}{{\mathbb R}}

\newcommand{\V}{{\mathbb V}}
\renewcommand{\U}{{\Upsilon}}

\newcommand{\Z}{{\mathbb Z}}

\newcommand{\hyp}{{\mathrm{hyp}}}

\def\B0{{\bold{0}}}

\renewcommand{\lq}{``}

\catcode`\@=12

\def\Empty{}
\newcommand\oplabel[1]{
  \def\OpArg{#1} \ifx \OpArg\Empty {} \else
  	\label{#1}
  \fi}
		
\newcommand{\comm}[1]{}
\newcommand{\comment}[1]{}

\begin{document}

\bigskip\bigskip

\title[Hausdorff dimension of Feigenbaum Julia sets]{Hausdorff
dimension and conformal measures \\of Feigenbaum Julia sets}
\author {Artur Avila and Mikhail Lyubich}

\address{
Laboratoire de Probabilit\'es et Mod\`eles al\'eatoires\\
Universit\'e Pierre et Marie Curie--Boite courrier 188\\
75252--Paris Cedex 05, France
}
\email{artur@ccr.jussieu.fr}

\address{
Mathematics Department and IMS\\
SUNY Stony Brook\\
Stony Brook, NY 11794, USA
}
\email{misha@math.toronto.edu}

\address{
Department of Mathematics\\
University of Toronto\\
Ontario, Canada M5S 3G3
}
\email{mlyubich@math.sunysb.edu}

%\thanks{
%   This work was supported in part by Sloan Research Fellowship
% and NSF grants DMS-8920768 and DMS-9022140.}
\date{\today}

\begin{abstract}
We show that contrary to anticipation suggested by the dictionary
between rational maps and Kleinian groups and by the ``hairiness
phenomenon'', there exist many Feigenbaum
Julia sets $J(f)$ whose Hausdorff dimension is strictly smaller than two.
% In fact, it  can be arbitrarily close to 1.
We also prove that for any Feigenbaum Julia set,
the Poincar\'e critical exponent $\de_\crit$
is equal to the hyperbolic dimension $\HD_\hyp(J(f))$.
Moreover, if $\area J(f)=0$ then $\HD_\hyp (J(f))=\HD(J(f))$.
In the stationary case, the last statement can be reversed:
if $\area J(f)> 0$ then $\HD_\hyp (J(f))< 2$.
%Moreover, in the case of periodic combinatorics,
%$\HD_\hyp (J(f))=\HD(J(f))$ if and only if $\area J(f)=0$.
We also give a new construction of conformal measures on $J(f)$ that implies
that they exist for any $\de\in [\de_\crit, \infty)$,
and analyze their scaling and dissipativity/conservativity properties.
\end{abstract}

\setcounter{tocdepth}{1}

\maketitle
\thispagestyle{empty} \def\IMSmarkvadjust{0 pt}
\def\IMSmarkhadjust{0 pt}
\def\IMSmarkhpadding{0 pt}
\def\IMSpubltext{Published in modified form:}
\def\SBIMSMark#1#2#3{
 \font\SBF=cmss10 at 10 true pt
 \font\SBI=cmssi10 at 10 true pt
 \setbox0=\hbox{\SBF \hbox to \IMSmarkhpadding{\relax}
                Stony Brook IMS Preprint \##1}
 \setbox2=\hbox to \wd0{\hfil \SBI #2}
 \setbox4=\hbox to \wd0{\hfil \SBI #3}
 \setbox6=\hbox to \wd0{\hss
             \vbox{\hsize=\wd0 \parskip=0pt \baselineskip=10 true pt
                   \copy0 \break%
                   \copy2 \break% 
                   \copy4 \break}}
 \dimen0=\ht6   \advance\dimen0 by \vsize \advance\dimen0 by 8 true pt
                \advance\dimen0 by -\pagetotal
	        \advance\dimen0 by \IMSmarkvadjust
 \dimen2=\hsize \advance\dimen2 by .25 true in
	        \advance\dimen2 by \IMSmarkhadjust

%
%   Check for publication info
%
%  \newread\jref
  \openin2=publishd.tex
  \ifeof2\setbox0=\hbox to 0pt{}
  \else 
     \setbox0=\hbox to 3.1 true in{
                \vbox to \ht6{\hsize=3 true in \parskip=0pt  \noindent  
                {\SBI \IMSpubltext}\hfil\break
                \input publishd.tex 
                \vfill}}
  \fi
  \closein2
  \ht0=0pt \dp0=0pt
 \ht6=0pt \dp6=0pt
 \setbox8=\vbox to \dimen0{\vfill \hbox to \dimen2{\copy0 \hss \copy6}}
 \ht8=0pt \dp8=0pt \wd8=0pt
 \copy8
 \message{*** Stony Brook IMS Preprint #1, #2. #3 ***}
}

\SBIMSMark{2004/05}{August 2004}{}

\tableofcontents

 \section{Introduction}

\subsection{Statement of the results}

One of the first questions usually asked about  a fractal subset of $\R^n$ 
is whether it has  the maximal possible Hausdorff dimension, $n$. 
It certainly happens if the set has positive Lebesgue measure. On the other hand,
it is easy to construct fractal sets of zero measure
but of dimension $n$.  Moreover,  this phenomenon is
often observable for fractal sets produced by conformal dynamical systems,
iterated rational functions or Kleinian groups.
In particular, the analogy with Kleinian groups suggested that the Julia sets of  
Feigenbaum maps should have Hausdorff dimension two. 
In this paper we will show that this is not always the case.  
% describe a phenomenon that goes against this common intuition.

A quadratic polynomial (or more generally, a quadratic-like map)
 is called {\it Feigenbaum} if it is infinitely renormalizable of
bounded combinatorial type with {\it a priori} bounds (see \S 2 for a precise definition). 
 Feigenbaum polynomials are remarkable
dynamical systems whose geometry has  been a focus
of  research for the past 25 years. By analogy with Kleinian groups,
it was anticipated that their Julia sets  have  dimension 2,
and this problem has been around for a while.
However, in this paper we will show that this is not always the case:

\begin{thmA}[$\HD<2$] %\label{HD<2}
 There exist Feigenbaum quadratic polynomials whose Julia sets $J(f)$ have
Hausdorff dimension strictly less than 2. 
% In fact, it can be arbitrary close to 1.
\end{thmA}

Maps to which this theorem applies  are infinitely renormalizable maps with ``high'' combinatorics.
To the best of our knowledge, previous examples of Julia sets of Hausdorff
dimension smaller than $2$ were related to one of two mechanisms:
\begin{enumerate}
\item Definite expansion along the critical orbit,
\item Porosity (presence of definite holes in small scales) of the
Julia set.
\end{enumerate}
Note that both mechanisms fail for Feigenbaum maps.

\msk
{\it Remark. } 
We show in \cite{AL} that 
in fact $HD(J(f))$ can be arbitrary close to 1 for some Feigenbaum Julia sets.
Namely, when the combinatorics of the renormalization gets close to the Chebyshev one 
(of the map $z\mapsto z^2-2$), the dimension gets close to 1. 
\msk

It is still unknown whether there exist Feigenbaum Julia sets of positive area,
or the ones with  zero area but Hausdorff dimension two.
Curiously, the affirmative answer to the former question would imply the affirmative answer to 
the latter  one (subject of certain {\it a priori} bounds assumptions):

\begin{thmB} %\label{dim 2 area 0} 
Assume that a renormalization horseshoe $\AAA$ (see \S \ref{horseshoe})
contains   maps $g_+$ and $g_0$ such that $\area (J(g_+))>0$ while  $\area(J(g_0))=0$. 
Then it contains a map $f$ such that  $\area(J(f))=0$  but $\HD(J(f))= 2$.
\end{thmB}  

In the course  of this paper, we give a criterion for $\area J(f) > 0$
which can provide an efficient numerical test on this property. 

\msk
Along with the Hausdorff dimension, there are other natural geometric objects and quantities
associated with  Julia sets: 

\ssk \nin $\bullet$
The {\it Poincar\'e series}
\begin{equation}\label{Poincare series} 
   \Xi(\de) \equiv  \Xi_\de (f,z) = \sum_{n=0}^\infty\sum_{f^n\zeta=z} \frac{1}{|Df^n (\zeta)|^\de }, \quad z\in \C\sm \OO,
\end{equation}
where $\OO$ is the postcritical set of $f$. 

\ssk \nin $\bullet$
 The {\it critical exponent} $\de_\crit=\de_\crit(f)$ which separates convergent exponents of the Poincar\'e series
from the divergent ones.   

\ssk \nin $\bullet$
The {\it hyperbolic dimension} $\HD_\hyp (J)$ of the Julia set $J=J(f)$, 
that is, the supremum of the dimensions of all invariant hyperbolic subsets of $J$.

\ssk \nin $\bullet$  {\it Conformal measures} on the Julia sets, that is the measures that transform according to the
rule:
\begin{equation}\label{trans rule}
    \mu(fX) = \int_X |Df|^\de d\mu
\end{equation}
for any measurable set $X$ such that $f|X$ is injective. 
By  Sullivan's Theorem \cite{S-conformal},
the Julia set $J$ supports a conformal measure with exponent $\de_\crit$.

\ssk \nin $\bullet$ 
 The {\it minimum exponent} $\de_*= \de_*(f)$ is the infimum of the exponents
of all conformal measures on the Julia set. 

\msk
The following general relation between the above quantities has been known:
\begin{equation}\label{general relation}
    \de_* = \HD_\hyp (J) \leq \HD(J),
\end{equation}
where the equality is due to  Denker and Urbanski \cite{DU} while the inequality is trivial.
If moreover the Julia set has zero area then the inequality
\be
\HD(J) \leq \delta_\crit
\ee
holds for a large class of maps (including Feigenbaum maps),
a result due to Bishop\footnote{Actually Bishop has proved the {\it 
a priori} stronger result that (under the assumption that the Julia set has 
zero Lebesgue measure) $\delta_\crit \geq \overline \MD$ , the {\it upper
Minkowski dimension}.  Together with this result of Bishop, our Theorem~C
below shows that the Minkowski dimension of a Feigenbaum Julia set 
is equal to its Hausdorff dimension.} \cite{B} (dynamically interpreted in
\cite{GS}).  However, it is not known whether these quantities in
general coincide.  In the case of Feigenbaum maps,  
(and in fact, for quite a general class of rational maps
(see \S \ref{remarks})), we can resolve this issue:

\begin{thmC}[Critical Exponents and Dimensions] %\label{exp vs dim}
  Let $f$ be a  Feigenbaum map.  Then 
$$
   \de_\crit(f) =\de_*(f)  = \HD_\hyp(J).
$$
Moreover, if $\area J(f) =0$, then $$\HD_\hyp(J) = \HD(J).$$ 
\end{thmC}

It is quite surprising that the last assertion can be reversed in the case
of periodic combinatorics:

\begin{thmD}[Positive Area and Hyperbolic Dimension] %\label{positive area}
Let $f$ be a Feigenbaum map with periodic combinatorics.  
If $\area(J)>0$ then $$ \HD_\hyp (J) < \HD(J) = 2.$$
\end{thmD}

At this moment we cannot tell whether this result is positive or negative  
(that is, whether it gives an interesting property of certain Julia sets
or indicates that those Julia sets should not exist).

Putting Theorems~C and~D together we obtain the following amusing statement: 
{\it A Feigenbaum map with periodic combinatorics has a Julia set of zero Lebesgue
measure if and only if its hyperbolic and Hausdorff dimensions coincide.}

\msk 
It is known (Prado \cite{Pra}) that any conformal measure on a Feigenbaum
Julia set is ergodic (and hence there exists at most one conformal
measure for any given exponent $\de$).  Another traditional question
asked about quasi-invariant measures is whether they are
dissipative or conservative\footnote{For Kleinian groups this question is related to existence of 
the Green function on the associated 3-manifold (see \cite{Ah}).}.
Recall that an ergodic
measure is called {\it conservative} if it satisfies the conclusion of the 
Poincar\'e Recurrence Theorem and is called {\it dissipative} otherwise.  

\begin{thmE}[Conformal Measures] %\label{conformal measures}
Let $f$ be a  Feigenbaum  map.
Then for any $\de\geq \de_* (f)$,  the Julia set $J(f)$
supports a unique conformal measure $\mu_\de$ with  exponent $\de$.
If $\mu_\de$ is conservative then $\de = \de_* < 2$. 
Moreover, both dissipative and conservative alternatives with $\de=\de_*<2$
are realizable for some Feigenbaum quadratic-like maps.
\end{thmE}

This gives  first examples of  dissipative $\de_*$-conformal measures.
% with $\de_* < 2$.\footnote{We note that when $\delta_*=2$ (as is the case in the
%Shishikura examples \cite {Sh-dim}), dissipativity of the conformal measure
%is automatic.}

%\bignote{On second thinking, I am not sure that Shishikura's examples are obviously
%dissipative, since our arguments don't apply to the case when the critical point is 
%is dense in $J$.}

\subsection{Methods}
In this paper we develop an efficient method to estimate the Poincar\'e series 
by comparing  its values on different renormalization levels. 
For periodic points of renormalization this method
leads to  a Recursive Quadratic Estimate for the Poincar\'e series.
Let $A^n$ be the  fundamental annuli on consecutive renormalization levels. 
Then we show that the average Poincar\'e series,
$$
         \om\equiv  \om_m(\de)=  \frac{1}{\area A^m} \int_{A^m} \Xi_\de(x) dx
$$
satisfies an estimate $\om \leq P_\de (\om)$, where $P_\de$ is a quadratic
polynomial (with positive coefficients) whose coefficients continuously
depend on $\de$ near $2$.  It gives a bound of the average Poincar\'e series
$\om$ by a real fixed point of $P_\de$ if it exists. 
Moreover, coefficients of $P_2$ involve areas of points that escape
annuli $A^m$, and one can see that if substantial  masses of points escape
the $A^m$ then $P_2(x)$ has an attracting fixed point. Persistence of
this point under a perturbation of $\de$ yields  convergence of the
Poincar\'e series for $\de$ slightly smaller than 2,
implying that $\de_\crit <2$.

More generally, this leads us to the following {\it Trichotomy}. 
For any Feigenbaum map with periodic combinatorics
one of the following three possibilities occurs: 

\ssk\nin
1) {\it Lean} case: substantial  masses of points escape the annuli $A^n$
to the {\it outer} regions.
      In this case
                        $$\HD(J)=\HD_\hyp(J)<2.$$

\ssk\nin
2) {\it Balanced} case: there is a balance between masses of points escaping
{\it in} and {\it out} 
     the fundamental annuli. In this case 
              $$  \area J = 0\quad \text{but} \quad   \HD_\hyp(J)=2. $$  

\ssk\nin
3) {\it Black hole} case: the critical point acts like a black hole that
attracts huge masses of matter, so that
    substantial  masses  escape the annuli $A^n$ to the {\it inner} regions. 
    In this case $$ \area J>0\quad \text{ but}\quad  \HD_\hyp(J)<2.$$

\ssk 
In the case of Feigenbaum maps with high combinatorics, 
 the critical point is lean which yields Theorem~A. %\ref{HD<2}. 
%On the other hand, when $\area(J)>0$, 
%the critical point acts like a black hole that  attracts huge masses of matter.
%Then  substantial  masses  escape the annuli $A^n$ to the {\it inner} regions,
%which yields 
Theorem~D %\ref{positive area}
is incorporated into  the black hole case. 
{\it Heuristically}, phase transition from the lean case to the black hole case
can occur only through the balanced case
% when   there is a balance between masses of points escaping
%{\it in} and {\it out} the fundamental annuli
indicating that Theorem~B should hold. %\ref{dim 2 area 0}.
%\note {I do not think this sentence is valid except at a heuristic level,
%since we do not have the Trichotomy for aperiodic combinatorics}

%\msk
%In the case when combinatorics is close to Chebyshev, a more
%direct analysis gives a stronger recursive quadratic estimate: 
%it shows that the corresponding critical exponent is close to 1
%(which is the dimension of $J(z^2-2)= [-2,2]$).  
 
\msk
Theorem~C is based upon  a new construction of a
$\de$-conformal measure by taking a properly truncated Poincar\'e series
$\Xi_\de^{tr} (f,z)$, putting the corresponding measure on the preimages
of $z$, and letting $z\to 0$. This allows us to demonstrate that conformal
measures exist for any $\de\geq \de_*$ which is the key step of the proof
(see also Theorem~E). %\ref{conformal measures}).

Dissipativity  of the $\de$-conformal measure is equivalent to convergence
of the Poincar\'e series $\Xi_\de$.
A simple  argument shows that $\Xi_{\de_*}(Rf,z)$ 
 is convergent, provided   the dimension does not drop under the
renormalization: 
$$
   \HD_\hyp (J(Rf)) = \HD_\hyp(J(f))\equiv \de_*,
$$
which holds, e.g., for renormalization fixed  points $f= Rf$. 
On the other hand, if the dimension drops,
\begin{equation}\label{drop}
   \HD_\hyp (J(Rf)) < \HD_\hyp(J(f)),
\end{equation}
then we find a linear recursive estimate
for the Poincar\'e series (this time averaging over the conformal measure
$\mu_{\de_*}$ on the annuli $A^n$), which would yield convergence of
$\Xi_\de$ for some $\de< \de_*$ if  the measure  $\mu_{\de_*}$ 
were dissipative (contradicting that $\de_*$ is the minimum exponent). 

Finally, we give two proofs that the drop condition (\ref{drop}) is satisfied for some
Feigenbaum maps,  which completes the proof of Theorem~E.  %\ref{conformal measures},
One of them,  based on methods of quasiconformal deformations,   shows that  
$$ 
  \sup_{f\in \HH} \HD_\hyp(f)=2,
$$  
where $\HH$  is any hybrid class of Feigenbaum quadratic-like maps.
Another one is based on a dynamical interpolation in an appropriate renormalization horseshoe. 

\msk
Our methods also work in the generalized renormalization contexts, 
with Fibonacci maps as main examples.  They can be also applied to the real
dynamics, where it is known that ``wild attractors'' exist \cite{BKNS}
(real analogue of the phenomenon  of ``positive area Julia sets''),
which  yields a real life version of Theorems~B %\ref{dim 2 area 0}
and~D. %\ref{positive area}.

%Some heuristics (partially supported by the  results of \cite {BKNS}) indicate,
Also, in the search for Feigenbaum Julia sets of  positive area, 
it might be a good idea to consider maps of   higher degree (that is, $\geq 2$)  as well. 
% with a critical point of higher order.
Let us  remark with  this respect that all the results and proofs of this
paper are still  valid for (unicritical) Feigenbaum maps of arbitrary degree
% with a critical point of higher order 
(with only terminological and notational adjustments).

\subsection{Some history}
Geometric  measures and Hausdorff dimension of limit and Julia sets are popular themes in
holomorphic dynamics. 

First examples of a conformal dynamical system with the limit set of
measure zero but Hausdorff dimension 2
appeared in the works of Thurston \cite{Th}
and Sullivan \cite{S-groups} on Kleinian groups.
On the other hand, the limit set of a geometrically finite Kleinian group%
\footnote{We follow the classical convention that
the limit set of a Kleinian group is not the whole sphere.}
has Hausdorff dimension less than~2
 \cite {Bo}, \cite {S4}, \cite {T}.
In this context, the Hausdorff dimension  problem   was eventually
resolved in  full by Bishop and Jones \cite{BJ}: the limit set
$\La(\Gamma)$ of a (finitely generated) Kleinian group has
Hausdorff dimension 2 if and only if the group is geometrically infinite. 
Also, it has been recently announced that the limit set of {\it any}
finitely generated Kleinian group has zero area,
which settles a long-standing  Ahlfors Conjecture,  see \cite{Ag,CG}.

Another basic class of conformal dynamical systems is given by iterates
of a single endomorphism, rational or transcendental.
The exponential map $z\mapsto \la e^z$, $\la<e^{-1}$,  was
 the first example of  such an endomorphism whose Julia set has
zero area \cite{EL} but Hausdorff dimension 2
\cite{McM-exp}.  Examples of entire functions with Julia sets
of positive area were given in \cite{EL2,McM-exp}.
For other interesting properties of the measure and dimension
of transcendental Julia sets, see \cite{exp,R,K,SZ,UZ}. 

The first examples of  rational Julia sets with Hausdorff dimension two
were constructed by Shishikura \cite{Sh-dim} (see also \cite{McM-triptich}).
Indeed, Shishikura  demonstrated that this
phenomenon is generic for quadratic polynomials
$z\mapsto z^2+c$ with $c$ on the boundary of the
Mandelbrot set.  All these Julia sets have zero area \cite{measure,Sh-ICM}.

On the other hand, the works of McMullen \cite{McM-Siegel} (see also \cite{Pe}) 
 and Przytycki \cite{Prz}  (see also \cite{PR,GS}) showed that the situation in
the iteration theory is more complicated than in the theory of Kleinian
groups: There exist geometrically infinite Julia sets
with Hausdorff dimension strictly less than two.  Here a map is
called geometrically finite if all the critical points on the
Julia set are non-recurrent (see \cite{LMin}). The Hausdorff
dimension of any geometrically finite Julia set 
is strictly less than 2, provided it is not the whole sphere  (see \cite{CJY,U,McM-triptich}). 

As we have already mentioned,  it is still unknown whether there exists
a polynomial map whose Julia set has  positive area.
However, it was known that many Feigenbaum Julia sets
have zero area \cite{Y}. 
(It is essentially  the same class of maps for which {\it a priori}
bounds had been established in \cite{puzzle}
and to which our Theorem~A %\ref{HD<2}
applies.) 

 McMullen  \cite{McM-towers}  has demonstrated that in the dictionary
between rational maps and Kleinian groups,
Feigenbaum maps correspond to hyperbolic 3-manifolds that fiber over
the circle. Using this analogy, he proved
that the critical point lies ``deep" (in some precise geometric sense)
inside any Feigenbaum Julia set. This strongly suggested that these sets
should have the maximal possible Hausdorff dimension. However,
as we have already mentioned, these points can be lean in the
measure-theoretic sense which leads to the opposite conclusion. 

For Kleinian groups, conformal measures were introduced by Patterson \cite{Pa} 
and Sullivan \cite{S-density}.
% See \cite{Ah} for the discussion of their conservativity/dissipativity
% properties in relation with the potential theory on 3-manifolds.
See \cite{N} for a discussion of their ergodic properties.
Ergodic and conservativity/dissipativity properties of conformal  measures on Julia sets 
were  studied, among other papers, 
in \cite{old,McM-renorm,Pra,Ba,BM,GS} 
(part of this work was motivated by \cite{BL} concerning real dynamics).   

For the relations between the critical exponent and Hausdorff dimension for Kleinian groups
 see \cite{BJ} (compare Theorem C).

\subsection{Organization of the paper}

We begin the paper (\S 2) with a basic material concerning the
Hausdorff dimension, Poincar\'e series, conformal measures,
and renormalization of quadratic-like maps. In \S \ref{hyp dim}
we study scaling and ergodic  properties of conformal measures,
and prove Theorem~C  and part of Theorem E.

Section \ref{pert est} is the central in the paper:  here we derive the Recursive Quadratic Estimates. 
In \S \ref{trichotomy} we deduce from them the Trichotomy.
It immediately yields  Theorem~D and provides us with the ``leanness''
condition  for Theorem A. 

In \S \ref{conservativity} and \S \ref {varying dim}
we complete a proof of Theorem E.
In \S \ref{varying dim1} we study how the Hausdorff dimension
$h(g)\equiv \HD (J(g))$
varies within  hybrid classes of  Feigenbaum maps.
% We show that $\sup h(g)$ over any such hybrid class is always equal to 2. 
One of the consequences of our analysis is the following:
If $f$ has stationary combinatorics then the function $h(g)$ is either
identically equal to 2, or else $\operatorname{Im} h = [h(f_*), 2)$
where $f_*$ is the associated renormalization fixed point.

In \S \ref{interpolation}, we prove Theorem~B.

In \S \ref{Yarr} we describe a class of well controlled
Feigenbaum maps \cite{puzzle} whose Julia sets are lean
at the critical point \cite{Y},
which completes the proof of the first assertion of  Theorem~A.
The description is given in terms of the principal nest of the Yoccoz puzzle.

%In \S  \ref{Chebyshev} we analyze the case of almost Chebyshev combinatorics
%which gives us the second assertion of Theorem~A. %\ref{HD<2}. 
%The reader who is only interested to see a Feigenbaum Julia set with
%Hausdorff dimension  close to 1  can go directly to this section.

\S \ref{remarks} contains various remarks and some open problems.

The paper is concluded with two appendices.
Appendix A is of technical nature: 
here  we construct fundamental domains for Feigenbaum maps satisfying nice Markov properties.
In Appendix \ref{towers-sec} we analyze conformal measures
on towers which yields sharper scaling properties of conformal measures
on the Julia sets. 

\medskip{\bf Acknowledgement.}
This work was partly supported by the the Guggenheim Foundation, Clay Foundation, NSF and NSERC. 
It was done during authors' meetings at  SUNY at Stony Brook,  University of Toronto,
and Institut Henri Poincar\'e. The authors thank all these Institutions and Foundations
for hospitality and support.

 \section{Background}\label{background}

\subsection{Notations and terminology}
  $\C$ is the complex plane; \\
  $\N=\{0,1,\dots\}$ is the set of natural numbers;\\
  $\Z$ is the set of integers;\\ 
  $\D_r(z)=\{\zeta: |z-\zeta|<r\}$, $\D_r\equiv \D_r(0)$, $\D\equiv \D_1$;\\
$U\Subset V$ means that $U$ is {\it compactly contained } in $V$, i.e., 
$\cl U$ is compact and is contained in $V$. 

A domain is a connected open set.
By  a (topological) disk we will mean a simply connected domain $U\subset \C$. 
A Jordan disk is a topological disk whose boundary is a Jordan curve. 
{\it Filling} of a domain $D$ is the union of all Jordan disks with boundary
contained in $D$: it is the smallest simply connected set containing $D$.
Given a family of topological disks, we say that they have a
{\it bounded shape} if they are quasi-disks with uniformly bounded dilatation.  
A (topological) annulus $A\subset \C$ is a doubly connected domain. Any topological annulus
can be conformally uniformized by a round annulus $\{z: r<|z|<R\}$ (where possibly $r=0$
or $R=\infty$). The modulus of the topological  annulus, $\mod(A)$, is defined
as $\log(R/r)$.

The Lebesgue measure ($\equiv$ area)  of a measurable set $X \subset \C$ will be denoted by $\area X$ or $|X|$. 
 If $X,Y \subset \C$ are measurable and $|Y|>0$, we let $p(X|Y)=|X \cap Y|/|Y|$.

{\it Distortion} of a (holomorphic) univalent  map $f: U\ra \C$ is defined
as 
%D(f)=
    $$\sup_{z, \zeta\in U} \frac {|Df(z)|} {|Df(\zeta)|}.$$
%By the {\it nonlinearity} of $f$ we will understand the logarithm of the  distortion.

%\bignote{Are all these notations adequate?}

\lq{Quasi-conformal}" will be abbreviated as \lq{qc}''.

\ssk
For a natural $s\geq 2$, let $\Sigma_s$ 
 denote the space of two-sided sequences $(\eps_k)_{k=-\infty}^\infty$
in $s$ symbols:  $\eps_k \in \{1,\dots ,s\}$. 
 The shift transformation $\sigma$ on this space is called a
{\it Bernoulli shift}.

\comm{
Let $A=(a_{ij})_{i=1}^s$ be a {\it transition matrix} of 0's and 1's. It determines an invariant
subset $\Sigma^+_A\subset \Sigma^+_s$ consisting of sequences $\{i(k)\}$ with
$A_{i(k), i(k+1)}=1$. The restriction of the shift to $\Sigma^+_A$ is called
a {\it Markov shift} (or a {\it subshift of finite type}).

A matrix $A$ is called {\it primitive} if $A^n$ has strictly positive entries for some
$n\in \N$. It corresponds to the mixing property of the Markov shift.
}

\ssk The forward orbit of a point $z$ under $f$  is denoted as $\orb(z)\equiv \orb_f z= \{f^n z\}_{n=0}^\infty$;\\
Its limit set is denoted as $\om(z)$.

%\ssk
% The space of plane measures is supplied with topology of week convergence on compact subsets.
%In other words, $\mu_n\to \mu$ if for any continuous function with compact support,
%   $$\int \phi d\mu_n\to \int \phi d\mu.  $$

\ssk
$a\asymp b$ means that $C^{-1}<a/b< C$ with a constant $C>0$ independent
of particular $a$ and $b$ under consideration.

Assertion that ``some property is satisfied for $\de \approx a$'' 
 naturally  means that it is satisfied for all $\de$ sufficiently close to $a$. 

%As usual,  a constant, say $C$,   can change in the course of the argument without changing notations. 

\subsection{Hyperbolic sets}

Let $X\subset \C$ be a compact set invariant under a certain  holomorphic map $f$ defined in a
neighborhood of $X$.
% Let us call such a set \lq{dynamical}".
It is called {\it hyperbolic} if there exist $C>0$ and $\la>1$ such that for all $z\in X$
$$
         |Df^n(z)|\geq C\la^n,\; n=0,1,\dots
$$

\subsection{Quadratic-like maps} \label{q-l maps}
%  See \cite{DH,McM1} for a background in the theory of quadratic-like maps.
% Below we will briefly recall the definitions.
Let $\u\Subset \V$ be two topological disks. A (holomorphic) branched
covering of degree $d$ $f: \u\ra \V$ is called
a {\it polynomial-like map} of degree $d$.
If $d=2$, the map is called {\it quadratic-like}.  
%A quadratic-like map is called {\it real} if it preserves the real line.
% These maps are considered up to  affine conjugacy.
% Unless otherwise is assumed, we
% normalize such a  map to put its critical point at the origin.
%\marginpar{Should we introduce unicritical polynomial-like maps?}

The {\it filled  Julia} set $K(f)$ of a polynomial-like map is defined as
the set of non-escaping points (where \lq{escaping}" means landing at the
fundamental annulus
$\V\sm \u$ under some iterate of $f$), and the {\it Julia set} $J(f)$ is defined as the boundary of
$K(f)$.  In the quadratic-like case, 
these sets are either Cantor or connected depending on whether the critical point
itself is escaping or not.

 Any quadratic-like map has two fixed points counted with multiplicity. In the case of
connected Julia set these two points can be dynamically distinguished.  One of them, usually
denoted by
$\alpha$, is either non-repelling or {\it dividing}, i.e., removing of it makes the Julia set
disconnected. Another one, denoted by $\beta$, is  always {\it non-dividing}.

 If a quadratic-like map is considered up to choice of domains $\u$ and $\V$,
then it should be more carefully called  a quadratic-like {\it germ}.
More precisely, one says that two quadratic-like maps with connected Julia sets
 represent the same germ if they have a common Julia set and coincide in a neighborhood
of it. For such a  germ,
let $\mod(f)=\sup \mod(\V\sm \u)$
where the supremum is taken over all possible choices of domains $\u$ and $\V$.
We will not make notational difference between maps and germs.
% The bigger
% $\mod(f)$, the better geometric control over $f$ we have (the closer $f$ to being purely
% quadratic).

We will assume below that 0 is the critical point of $f$.

Two quadratic-like maps/germs $f$ and $g$ are called {\it hybrid equivalent} if there is  a qc
conjugacy $h$ between them with $\dibar h=0$ a.e. on the filled Julia set $K(f)$.
By the Straightening Theorem \cite{DH-pl}, any quadratic-like map $f$
with connected Julia set  is hybrid equivalent
to a unique quadratic polynomial $z\mapsto z^2+\chi(f)$, with $c=\chi(f)$ on the Mandelbrot set $M$.
Let $\HH_f$ stand for the hybrid class of $f$.

\subsection{Feigenbaum maps}
 A quadratic-like germ $f: \u\ra \V$ is called {\it renormalizable} if
there exist $p=p(f)>1$ and
topological disks $\u'\Subset \V'$ containing the critical point such that:
\begin{itemize}
  \item $g\equiv f^p: \u'\ra \V'$ is a quadratic-like map with connected Julia set $J(g)$ called a {\it pre-renormalization} of $f$;
  \item The \lq{little Julia sets}"
 $f^k J(g)$, $k=0,1, \dots, p-1$, are pairwise disjoint except perhaps
         touching at their $\beta$-fixed points.
\end{itemize}
   If the little Julia sets indeed touch then one says that $f$ is {\it immediately
renormalizable} (or  that the renormalization is of {\it satellite type}).
Otherwise the renormalization is called  {\it primitive}.  
% In the real case, only the doubling renormalization (with period $p=2$) is of satellite type.
% Assume that $f$ is renormalizable and $p>1$ is the minimal possible renormalization
% period.

A pre-renormalization
considered up to affine conjugacy is called a  {\it renormalization} of $f$.
The renormalization of $f$  with minimal possible period is denoted  $Rf$.

%If $f$ is a real renormalizable map with minimal period $p$, then
%the combinatorial type of its renormalization  is defined
%by the order of the orbit $\{f^k 0\}_{k=0}^{p-1}$ on the real line.

The Mandelbrot set contains a plenty of canonical homeomorphic copies
of itself.  A quadratic-like map $f$ is renormalizable if and only if its
straightening $\chi(f)$ belongs to one of these copies.  The choice
of this copy specifies the {\it combinatorics} of the renormalization.  

Now one can inductively define an infinitely renormalizable map  $f$;
its successive renormalizations are denoted by $R^m f$, $m=0,1,2,\dots.$
If all these renormalizations have the same combinatorial type,
one says that $f$ has a {\it stationary type}.
If the sequence of the combinatorics is periodic, then $f$ has a
{\it periodic type}.  If the periods of all renormalizations are
bounded then  $f$ has a {\it bounded type}.  

One says that an infinitely renormalizable map has {\it a priori} bounds
if $\mod(R^n f) \geq \eps >0$, $n=0,1,2,\dots$. 
Conjecturally {\it a priori} bounds  are valid for all infinitely
renormalizable maps of bounded type.  So far, this conjecture has been
confirmed for real maps (see \cite{S-renorm,MS}) 
and for a class of complex maps of high type (see \cite{puzzle} and
\S \ref{Yarr}).

An infinitely renormalizable quadratic-like map of
bounded type with {\it a priori} bounds
will  be also called a {\it Feigenbaum map}.

\subsection{Renormalization horseshoe}\label{horseshoe}

A renormalizable quadratic-like map is called the {\it renormalization fixed}
(respectively, {\it periodic}) {\it point} if $Rf=f$
(respectively $R^p f = f$ for some $p$).  More generally, let us pick
a finite family $\MM$ of little Mandelbrot copies $M_k$, $k=1,\dots, s$.
Assume that to any two-sided sequence
$\bar M = (M_{ k(n) })_{k=-\infty}^\infty\in \Sigma_s$ we can associate a
quadratic-like germ $g_ {\bar M}$ with uniform {\it a priori} bounds which
is renormalizable with combinatorics $M_{k(0)}$ and such that
$R f_{\bar M} = f_{\si (\bar M)}$ where $\si: \Sigma_s\ra \Sigma_s$ is the
Bernoulli shift. Then the family $\AAA_\MM$ of these quadratic-like germs
is called the {\it renormalization horseshoe} associated to the given
family of copies. Notice that by definition, $\AAA$ is invariant under
the renormalization, and $R: \AAA_\MM\ra \AAA_\MM$ is 
topologically conjugate to the Bernoulli shift $\si$. 
 
We say that a finite family $\MM$ of little Mandelbrot copies is {\it fine} if for any one-sided sequence
$\bar M_+= (M_{k(n)})_{k=0}^\infty $ of these copies there is a quadratic-like germ which is renormalizable
with combinatorics  $\bar M_+$ with {\it a priori} bounds. 
For instance, any family of real Mandelbrot copies  is nice.  

\begin{thm}[\cite{S-renorm,McM-towers}] 
 Any fine family $\MM$ of Mandelbrot copies has an associated renormalization horseshoe $\AAA_\MM$.
\end{thm}   

Notice that periodic points of renormalization  are dense in the horseshoe $\AAA_\MM$. 

%\begin{thm}[\cite{S-renorm,McM-towers}]\label{horseshoe}
%For any real combinatorics, there exists a unique real Feigenbaum map $f$ with this combinatorics
%such that $Rf=f$. Moreover, any real Feigenbaum map with this combinatorics is hybrid equivalent
%to  $f$. 
%\end{thm}

% This above map $f$ (\lq{the renormalization fixed point}") has the 
%scaling invariant Julia set. More precisely, there exists a constant $\rho>1$ such that 
%the \lq{Feigenbaum dilation}"
% $T_\rho: z\mapsto \rho z$ maps $J^{m+1}$ onto $J^m$ conjugating $R^{m+1} f$ to $R^m f$.  

\begin{thm}[\cite{FCT}]\label{hyperbolicity}
  The horseshoe $\AAA_\MM$ is hyperbolic with codimension-one stable foliation.
Moreover, local stable manifolds of points $f\in \AAA_\MM$ coincide with their hybrid classes.
\end{thm}  

\subsection{Nest of little Julia sets}
Let $f$ be a Feigenbaum quadratic-like map.  
Let us consider a sequence $f_m: \u^m\ra \V^m$ of $m$-fold pre-renormalizations of $f$.
Let 
$J\equiv J^0\supset J^1\supset J^2\supset\dots$ stand for the corresponding nest of little Julia sets,
$J^m = J(f_m)$. 
Then
$$
        \OO \equiv \OO(f)=  \om(0) =\bigcap_{m \geq 0} \bigcup_{i \geq 0} f^i(J^m)
$$
is  the {\it postcritical set} of $f$. 

\begin{lem}[see e.g.,  \cite{McM-towers}]\label{shapes}
There exist quadratic-like pre-renormalizations $f_m: \u^m \ra \V^m$
with the following properties:
\begin{itemize}
  \item[(P1)] $\cl(\V^m \sm \u^m)\cap \OO= \emptyset$; 
  \item[(P2)] {\it Unbranched  a priori bounds}: $\mod(\V^m\sm \overline \u^m) \geq \mu>0$ and $\mod(\u^m\sm J^m)\geq\mu>0$;
  \item[(P3)] The disks  $\u^m$ and $\V^m$ are quasidisks with  bounded dilatation;
% \item  The inverse branches $J^m\ra J^m_i$ admit univalent extensions  
%        $\u^m\ra \u^m_i$ with bounded distortion;
  \item[(P4)] $\diam J^m \asymp \diam \V^m$. %\asymp \rho^m $ for some $\rho\in (0,1)$. 
\end{itemize}
  All the constants depend only on the combinatorial
and {\it a priori} bounds for $f$.
\end{lem}

\ssk
If $f$ is a renormalization fixed point then
there exists a dilation $\la > 1$ such that $f_{m+1} = \la f_m (\la^{-1} z)$,
and the above picture becomes  scaling invariant:
$ \V^m= \la^m \V^0$, $\u^m = \la^m \u^0$, and $J^m = \la^m J^0$.  
In the case of renormalization periodic point, the scaling invariance holds on the subsequences of levels
congruent modulo the period. 

{\it In what follows, the nests of domains $\V^m$ and $\u^m$ will always be assumed to satisfy 
properties (P1)-(P4). In the case of periodic combinatorics, they will be
additionally assumed to be scaling invariant on the appropriates subsequences of levels.}   

\subsection{Nice fundamental domains} \label {nfd}

Let us consider a Feigenbaum map $f:\u \to \V$ with pre-renormalizations $f_n:\u^n \to \V^n$. 
% Without loss of generality, see Lemma~\ref {shapes}, we
%assume unbranched {\it a priori bounds}:
%\begin{enumerate}
%\item [(C0)] $\overline \V^n \setminus \u^n$ is disjoint from the
% postcritical set of $f$,
% \item [(G0)] $\mod(\V^n \setminus \overline \u^n)$ is bounded from below
% (such a lower bound is called an unbranched {\it a priori} bound for $f$).
% \end{enumerate}
Along with the domains $\V^n$ and $\u^n$, we will also need simply connected domains $V^n$ and $U^n$ 
satisfying the following ``nice'' topological and geometric properties:
\begin{enumerate}
\item [(C1)] $\V^n \cap \OO(f) \subset V^n \subset \u^n$,
\item [(C2)] $V^{n+1} \subset U^n \equiv f_n^{-1}(V^n)$,
\item [(C3)] $f^k(\partial V^n) \cap V^n=\emptyset$, $k \geq 0$ 
             (compare with ``nice intervals''  of Martens \cite{Ma}).
%\item [(C2)] $U^n \equiv f_n^{-1}(V^n) \subset V^n$,
%\item [(C3)] Any univalent pullback of $V^n$ (under iterates of $f$) which
%intersects $V^n$ is contained in $A^n \equiv V^n \setminus U^n$,
%\item [(C4)] $V^{n+1} \subset U^n$.
\end{enumerate}
%We shall also assume the following geometric properties:

\msk
\begin{enumerate}
\item [(G1)] $A^n \equiv V^n \setminus U^n$
is far from the postcritical set
$\OO(f)$: it has bounded hyperbolic diameter in $\V^n \setminus \OO(f)$.
\item [(G2)] $\area(A^n) \asymp \area(U^n) \asymp (\diam(U^n))^2 \asymp
(\diam(V^n))^2$.
\end{enumerate}
The bounds (G1)-(G2) together with the unbranched {\it a priori} bounds
%, together with an upper bound on the degree of $f$,
%and a lower bound on the modulus of $V^n \setminus \u^n$,
will be called {\it geometric bounds}.  By {\it combinatorial bounds} we will mean a bound on
the renormalization period of all renormalizations of $f$.

In Appendix \ref{mark} we will construct such a nest of nice fundamental domains for any Feigenbaum map.
% an infinitely renormalizable map with unbranched complex bounds
%\footnote {This includes all real maps, and conjecturally all infinitely
%renormalizable maps of bounded type.} and
%such that none of its renormalizations is immediately renormalizable.
%Immediately renormalizable maps can still be dealt with by our method
%if one relaxes the hypothesis that $J^n \subset V^n$ to allow the boundary
%$\di V^n$ to pass through the $\beta$-fixed point of $f_n|J^n$.
%In order not to obscure the basic constructions, we will often restrict
%ourselves to the primitive case.

\begin{rem} \label {bounded distortion}

Condition (C3) implies that the first landing map to $V^n$ is particularly
simple: its restriction to each component of its domain is a univalent map
onto $V^n$.  Due to (P1) each univalent pullback of $V^n$ extends to a
univalent pullback of $\V^n$, so by the Koebe Distortion Lemma
the first landing map to $V^n$ has bounded distortion.

Another consequence of (C3) is the structure of the first return map to
$V^n$.  We have that $U^n \subset V^n$ and that $f_n:U^n \to V^n$ is a
double covering.  Thus $U^n$ is one of the components of the domain of the
first return map to $V^n$.  All other components are univalent pullbacks of
$V^n$.

\end{rem}

\subsection{Hausdorff dimension}
% For a detailed account of the notion  of Hausdorff dimension see e.g., \cite{Mat}.
% To be definite, let everything below happen inside the complex plane $\C$.
%
Given a $\delta\geq 0$, the Hausdorff $\delta$-measure $h_\delta$ is defined as follows
$$
       h_\delta(X)=\lim_{\eps\to 0} \inf \sum (\diam U_i)^\delta,
$$ 
where the infimum is taken over all coverings of $X$ by sets $U_i$ of diameter at most
$\eps>0$.   For any $X$ there is a unique critical exponent separating infinite and vanishing
values of the measure $h_\delta$. This exponent is called  the Hausdorff dimension of $X$ and is
denoted by $\HD(X)$.

Given a Borel measure $\mu$, the Hausdorff dimension $\HD(\mu)$ is defined as the infimum of the
$\HD(X)$ as $X$ runs over all measurable sets of full measure.
{\it Local dimension} $\HD_\mu(z)$ of a measure $\mu$ at a point $z$ is defined as
\begin{equation}\label{local dim}
\lim_{r\to 0} {\frac {\log \mu(\D_r(z))} {\log r}},
\end{equation}
 provided this limit exists.
If the limit does not exist, one can still consider {\it upper} and {\it lower} local dimensions,
$\overline{\HD}_\mu(z)$ and $\underline\HD_\mu(z)$,
by taking upper and lower limits in (\ref{local dim}).

The {\it hyperbolic dimension} $\HD_\hyp(X)$ of an invariant set $X$ is  defined as  the supremum of
the dimensions of all invariant hyperbolic subsets of $X$.

\comm{
\begin{thm}[\cite{Z}]\label{Zdunik}
  For any polynomial-like map $f: U\ra V$ with connected Julia set, there is an alternative:
\begin{itemize}
  \item Either $f$ is conformally conjugate to $z\mapsto z^d$ or the Chebyshev polynomial $\Ch_d$;
  \item or $\HD(J(f)) >1$.
\end{itemize}
\end{thm}
}

\begin{thm}[\cite{Z}]\label{Zdunik}
  For any polynomial-like map $f: U\ra V$ with connected Julia set, there is an alternative:
\begin{itemize}
  \item Either $J(f)$ is a real analytic curve;
  \item or $\HD(J(f)) >1$.
\end{itemize}
\end{thm}

\subsection{Conformal measures and critical exponents}
The {\it Poincar\'e series} $\Xi_\de(z)\equiv \Xi_\de(f,z)$, $z\in \C\sm \OO$,
was defined in the Introduction (\ref{Poincare series}).
Note that by the Koebe Distortion Lemma it has a bounded oscillation on any compact set $K
\subset \V \setminus \OO$:
\begin{equation}\label{bounded oscillation}
        \Xi_\de(f,z) \leq C \Xi_\de(f,\zeta), \quad z, \zeta\in K,
\end{equation}
where the constant $C$ depends only on the hyperbolic diameter of $K$ in $\V \sm \OO$, 

By definition, the {\it critical exponent}  $\de_{cr}=\de_{cr}(f)$ 
 separates convergent  values of $\de$ in $\Xi_\de$  from the divergent ones: 
the Poincar\'e series is convergent for $\de>\de_\crit$ and divergent for $\de<\de_\crit$.
Notice that by (\ref{bounded oscillation}), $\de_\crit$ is independent of the particular choice of 
$z\in\V\sm \OO$.  
For $\de=\de_\crit$ the Poincar\'e series can behave in both ways, and $f$ is called
of {\it convergent} or {\it divergent} type depending on it.  

Taking a wandering disk $D\Subset \V \sm J$ (i.e., $f^{-n}D\cap f^{-m}D =\emptyset$ for $m\not=n$), we see that
$$
   \Xi_2(z) \asymp \bigcup_{n=0}^\infty \area( f^{-n} (D)) <\infty,
$$ 
so that $\de_\crit\leq 2$. One can also show that $\de_\crit >0$.
In fact, we have:

\begin{lem}\label{de ge 1}
   If the Julia set $J(f)$ is connected then $\Xi_1=\infty$, so that $\de_\crit\geq 1$. 
\end{lem}

\begin{pf}
Let $\Gamma_n =\di( f^{-n} \u$). Then 
$\Xi_1 \asymp \sum \length (\Gamma_n)= \infty$.
\end{pf}

\msk
A quasi-invariant measure $\mu$ 
is called $\delta$-{\it conformal} if its Radon-Nicodim Jacobian
$ d (f^* \mu)/ d\mu$ is equal to $|Df|^\delta$,
which is equivalent to the transformation rule (\ref{trans rule}). 
This transformation rule and density of preimages on the Julia set imply that 
$\supp(\mu)=J(f)$ for any conformal measure $\mu$  on the Julia set.

\begin{thm}[Sullivan \cite{S-conformal}]\label{Sullivan's meas}
  Any polynomial-like map $f$ has at least one $\de_\crit$-conformal measure $\mu$ on
its  Julia set. 
\end{thm}

We have also defined in the Introduction the {\it minimal exponent} $\de_*=\de_*(f)$,  the infimum of the
exponents $\de$ for which $f$ admits a $\de$-conformal measure on the Julia set.  

\begin{thm}[Denker-Urbansky \cite{DU}]\label{DU formula}
   For any polynomial-like map, $$\de_*(f)=\HD_\hyp(J(f)).$$ Moreover, $f$ admits a $\de_*$-conformal measure
on its Julia set. 
\end{thm}  

Let us also mention the following simple property:  

\begin{lem}\label{semi-cont}
  The exponent $\de_*(f)=\HD_\hyp(f)$ depends lower semicontinuously on the quadratic-like map $f$.
\end{lem}

Recall that an $f$-quasi-invariant measure  $\mu$ on $J$  is called {\it ergodic} if there
is no decomposition  $J=X_1\cup X_2$ into two invariant measurable subsets of positive
measure.

\begin{thm}[Prado \cite{Pra}]\label{ergodicity} Let $f$ be  a Feigenbaum map.  
Then any conformal measure $\mu$ on $J(f)$ is ergodic.
 Hence for any  given exponent $\de$,  there exists at most one  $\de$-conformal measure. 
\end{thm}

\subsection{Poincar\'e series in the hyperbolic case}
Let us collect here some basic estimates on Poincar\'e series.
They show that in the hyperbolic case the Poincar\'e series cannot blow up:
it depends continuously on both $\de$ and $f$. 

Let $\Delta$ be the set of all continuous maps $f$ defined over a
compact set $X \subset \C$ that admit a holomorphic extension to a
neighborhood of $X$.  We say that $f_n:X_n \to \C$ converge to
$f:X \to \C$ is $X_n \to X$ in the Hausdorff topology and there exists a
neighborhood $U$ of $X$ such that for every $n$ sufficiently big, $f_n$
admits a holomorphic extension to $U$ and the $f_n:U \to \C$ converge to 
$f:U \to \C$ uniformly.  We will assume that $X$ has no isolated points (so
that derivatives of a function $f:X \to \C$ of class
$\Delta$ are well defined).

\begin{lemma} \label {continuity}

Let $f:X \to \C$ be a map of class $\Delta$
without critical points.  If $\delta>0$, $K>0$ are such that
\be
\Xi_\de(f,x)<K,
\ee
for every $x \in X$.
Then for every $\tilde K>K$, there exists $\tilde \delta<\delta$ such that
if $\tilde f:\tilde X \to \C$ is close
to $f:X \to \C$ then for every $x \in \tilde X$, we have
$\Xi_{\tilde \delta} (\tilde f,x)<\tilde K$.

\end{lemma}

\begin{pf}

For every $x \in X$, let $N_x \geq 0$ be minimal such that
\be
\sum_{f^{N_x}(y)=x} |Df^{N_x}(y)|^{-\delta}<\frac {1} {2 K}.
\ee
Notice that $N_x$ is upper semicontinuous, so $N=\sup N_x<\infty$.  It
follows that
\be
\sum_{f^N(y)=x} |Df^N(y)|^{-\delta}=
\sum_{f^{N_x} z =x} |Df^{N_x}(z)|^{-\delta} \sum_{f^{N-N_x}(y)=z}
|Df^{N-N_x}(y)|^{-\delta}<\frac {1} {2}.
\ee
%Let $M>0$ satisfy $2^M>2 K/(\tilde K-K)$.  We have
%\be
%\sum_{m=M N}^{M N+N-1} |Df^m(y)|<\frac {K} {2^M}<\frac {\tilde K-K} {2}.
%\ee

Choose an integer $M>0$ such that $2^{M-1}>  K/(\tilde K-K)$.
If $\tilde f:\tilde X \to \C$ is close to $f:X \to \C$ and $\tilde \delta$
is close to $\delta$, we still have, for every $x \in \tilde X$,
\be
\sum_{\tilde f^N(y)=x} |D\tilde f^N(y)|^{-\tilde \delta}<\frac {1} {2},
\ee
%\be
%\sum_{m=M N}^{MN+N-1}
%\sum_{\tilde f^m(y)=x} |D\tilde f^m(y)|^{-\tilde \delta}<
%\frac {\tilde K-K} {2},
%\ee
\be
\sum_{m=0}^{M N-1} \sum_{\tilde f^m(y)=x} |D\tilde f^m(y)|^{-\tilde
\delta}<K.
\ee
Hence
\begin{align}
\Xi_{\tilde\de} (\tilde f,x)&=
\sum_{m=0}^{M N-1} \sum_{\tilde f^m(y)=x} |D\tilde f^m(y)|^{-\tilde\delta}+
\sum_{m=0}^{N-1}\sum_{\tilde f^m(z)=x} \left( |D\tilde f^m(z)|^{-\tilde \delta}
\sum_{k=M}^\infty \sum_{\tilde f^{kN}(y)=z}
|D\tilde f^{kN}(y)|^{-\tilde \delta} \right) \\
\nonumber
&<K+\sum_{k=M}^{\infty} K 2^{-k}<\tilde K.
\end{align}
\end{pf}

\begin{lemma}

Let $f:X \to \C$ be a map of class $\Delta$
without critical points.
If the maximal invariant set $Q=\cap_{n \geq 0} f^{-n}(X)$ is hyperbolic
then there exists $K>0$ such that for every $x \in X$, we have
$\Xi_2(f,x)<K$.

\end{lemma}

\begin{pf}

%Let $\epsilon_0,\epsilon_1>0$ be such that $f$ admits a
%holomorphic extension to a
%$2 \epsilon_0$ neighborhood of $X$ and that for every $x \in X$,
%$f|\D_{2 \epsilon_0}(x)$ is univalent and $f(\D_{\epsilon_0}(x)) \supset
%\D_{\epsilon_1}(f(x))$.

Let $\epsilon>0$ be such that $f$ admits a holomorphic extension to an
$\epsilon$-neighborhood $U$ of $X$ and that $f|U$ has no critical points.

Let $Q_m=\cap_{0 \leq j \leq m} f^{-j}(X)$ so that $Q=\cap Q_m$.
Hyperbolicity of $\cap_{m \geq 0} f^{-m}(X)$ implies that there exists
$\rho>0$ such that if $y \in Q_m$ then there exists
a Jordan domain $D^m_y \subset U$ such that
$f^m:D^m_y \to \D_\rho(f^m(y))$ is univalent
and the diameter of $f^k(D^m_y)$, $0 \leq k \leq m$ is exponentially small
in $m-k$.

%By a density point argument, one sees that either
%$|Q|=0$ or it contains a disk $\D_\rho(x)$.
%The latter possibility gives a contradiction (since $\{f^k|\D_\rho(x)\}$
%would form a normal family and $Df^k|\D_\rho(x)$ would be bounded), so
%$Q$ has zero Lebesgue measure.

It is easy to see that $Q$ has empty interior (otherwise $\{f^m|\inter Q\}$
would form a normal family, contradicting the hyperbolicity of $f|Q$).
For $x \in X$, let us fix a compact set
$Z_x \subset \D_\rho(x) \setminus Q$ of positive Lebesgue measure.  We may
assume that $\{Z_x\}_{x \in X}$ form a finite family.  If $y \in Q_m$, let
$Z^m_y=(f^m|D^m_y)^{-1}(Z_{f^m(x)})$.  Since the diameters of the
$\{y\} \cup Z^m_y$ go to $0$ as $m$ grows, there exists $M>0$
such that $Z^m_y \cap Z_x=\emptyset$ whenever
$y \in Q_m$, $x \in X$, and $m \geq M$.
It follows that $Z^m_y \cap Z^{m'}_{y'}=\emptyset$ if $|m'-m| \geq M$.  We
also have $Z^m_y \cap Z^m_{y'}=\emptyset$ whenever $f^m(y)=f^m(y')$ and $y
\neq y'$.
%Let $Z^n=\cup_{m \geq n}\cup_{y \in Q^m} Z^m_y$.  Since $\cap
%Z^n=\emptyset$, we have $|Z^n \to 0|$.

Notice that, by the Koebe Distortion Lemma,
$  |Df^m(y)|^{-2} \leq C |Z^m_y|/|Z_{f^m(y)}| $ for some $C>0$.
Thus, for every $x \in X$,
$$
\sum_{m \geq 0} \sum_{\ntop{f^m(y)=x,} {y \in Q_m}}
|Df^m(y)|^{-2} \leq C \sum_{m \geq 0} \sum_{\ntop{f^m(y)=x,}{y \in Q_m}}
\frac {|Z^m_y|} {|Z_x|} \leq C M \frac {|U|} {|Z_x|}.
$$
%\be
%\sum_{m \geq n} \sum_{f^m(y)=x} |Df^m(y)|^{-2} \leq C \sum_{m \geq 0} \frac
%{|Z^m_y|} {|Z_x|} \leq C M \frac {|Z^n|} {|Z_x|} \leq \kappa_n,
%\ee
%where $\kappa_n \to 0$.  Let $N$ be such that $\kappa_N<1/2$, and let
%$\delta<2$ be such that
%\be
%\sum_{f^N(y)=x} |Df^N(y)|^{-\delta}<1/2.
%\ee
%In particular, there exists $m \geq 0$ such that for every $x \in X$.
\end{pf}

\begin{cor} \label {Delta}

Let $f:X \to \C$ be a map of class $\Delta$ without critical points.
If the maximal invariant set $Q=\cap_{n \geq 0} f^{-n}(X)$ is hyperbolic
then there exists $\delta<2$, $K>0$ such that for every $x \in X$ we have
$\Xi_\de (f,x )<K$.
Moreover, $\delta$ and $K$ can be chosen uniform over a compact family of
maps as above.

\end{cor}

\section{Critical exponent,  hyperbolic dimension, and dissipativity}\label{hyp dim}

{\it In what follows, $f$ is assumed to be a Feigenbaum map, unless otherwise is explicitly stated.}

\subsection{Dissipativity}
Recall that a measure $\mu$ is called {\it dissipative} if there exists 
 a wandering set $X$ of positive measure, i.e., $\mu(X)>0$ and 
$f^{-n} X\cap f^{-m} X=\emptyset$ for all $m > n\geq 0$. 

\begin{lem}\label{diss criterion}  
Let $\mu$ be a $\de$-conformal measure on $J$.
Then the  following properties are equivalent:

\begin{itemize}
 \item[(i)] $\mu$ is dissipative;
   \item[(ii)] The Poincar\'e series $\Xi_\de(z)$, $z\in \V\sm \OO $,  is convergent;
 \item[(iii)] $f^n z\to \OO$ for almost all $z$;
 \item[(iv)] For any little Julia set $J^m$, almost any $\orb(z)$ is eventually absorbed by
   the cycle of little Julia sets, $\cup_{k=0}^{p-1} J^m$ (here $p$ is the  period of $J^m$);   
 \item[(v)] $\mu(J^m)>0$ for all $m\in \N$;
   \item[(vi)] $\mu(J^1)>0$.
\end{itemize}

If $\mu$ is conservative then almost all orbits are dense in $J(f)$. 

\end{lem}

\begin{pf}
(i) $\implies$ (ii).
Let $X$ be a wandering set of positive measure. Taking an appropriate local branch of $f^{-n}$ 
on $X$, we can easily construct a wandering set 
$X'\subset \D(z, r)$ such that $\D(z, 2r)\cap \OO =\emptyset$. 
Then by the Koebe Distortion Lemma,
$$
    \Xi_\de(z) \asymp \sum  \mu(f^{-n} (X')) <\infty.
$$ 

\msk\noindent
  (ii) $\implies$ (iii). Take any disk $D=\D(z, r)$  such that $\D(z, 2r)\cap \OO =\emptyset$. 
Then 
$$
   \sum  \mu(f^{-n} (D)) \asymp \Xi_\de(z) <\infty.
$$
By the Borel-Cantelli Lemma, for a.e. $\zeta\in J$, the $\orb (\zeta)$ visits $D$ only finitely many times,
which implies (iii).

\msk\noindent
    (iii) $\implies$ (iv). 
For a pre-renormalization $f^p: \u^m\ra \V^m$ with period $p$, let $\u_k^m = f^k (\u^m)$, $k=0,1,\dots, p$
(so that $\u_0^m= \u^m$ and $\u_p^m = \V^m$).  

Since the little Julia sets $f^k J^m$ may touch each other 
only at their $\beta$-fixed points, which are not contained in 
the postcritical set $\OO$, there is a choice of pre-renormalization
such that the sets $\u_1^m\cap \OO,\dots \u_p^m\cap \OO$ are pairwise disjoint. 
Then there exist neighborhoods $\Om_k=f^k(\Om_0) \subset \u_k^m $ 
of $\OO \cap \u_k^m$ such that each  
 $\Om_k$ is disjoint from all  $\Om_i$ with $i\not\equiv k\ \mod p$. 

 If $f^n z\to \OO$, then $$ f^n z\in \bigcup_{k=0}^{p-1} \Om_k , \quad n\geq N .$$
Let us take some moment $n\geq N$ for which $f^n z\in \Om_0$. 
Then $f^{n+p} z\in \Om_p$ which is disjoint from all $\Om_k$, $k=1,\dots, p-1$.
Hence $f^{n+p}\in \Om_0$.
It follows  that $f^{n+pl}\in \Om_0\subset \u^m$ for all $l=0,1,\dots$.
Thus $f^n z\in J^m$. 

\msk\noindent
    (iv) $\implies$ (v).
Take a little Julia set $J^m$. 
By (iv), $\bigcup_{n=0}^\infty f^{-n} (J^m)$ has  full measure. This implies (v) by the $\de$-covariance of $\mu$. 

\msk\noindent Of course, (v) $\implies$ (vi). 

\msk\noindent (vi) $\implies$ (i).
Full preimage $f^{-1}(J^1)$ consists of two symmetric components, $f^{p-1}(J^1)$ and $X$.
If $\mu(J^1)>0$ then   $X$ is a wandering set of positive measure.  

\msk
Let us prove the last assertion. 
If it is not valid then  there is a disk $D$ (intersecting $J$)  
and a forward invariant set $X\subset J$ of positive measure such that $D\cap X=\emptyset$. 
Since $f^N D= J$ for some $N$, $D$ contains a set $Y$ of positive measure such that $f^N Y \subset X$.
Hence $f^n Y \cap Y=\emptyset$ for $n\geq N$. This easily implies that
$$ \mu \{ y\in Y: \ f^n y \not\in Y, \ n=1,2, \dots \} >0,  $$  
which gives us a wandering set of positive measure. 
\end{pf}
   
Since $\Xi_2(z)<\infty$, we conclude: 

\begin{cor}[compare\cite{old}]\label{diss1}
  The 2-conformal measure $\de_2$ is always dissipative.
\end{cor}

\subsection{Minimal conformal measure}
 
Let $\Spec(f)$ denote the {\it conformal spectrum} of $f$, that is, the set of exponents $\de$
for which there exists a $\de$-conformal measure on $J(f)$. Since weak$^*$ limits of conformal measures
are conformal, we have:

\begin{lem}\label{Spec}
 For any quadratic-like map $f$, $\Spec (f)$ is a closed non-empty subset of $\R_+$.
\end{lem}

In particular, the minimal exponent 
$$ 
    \de_*=\de_*(f) =  \inf \Spec(f)
$$  
belongs to the conformal spectrum.  The corresponding $\de_*$ conformal measure
$\mu_*$ will be called {\it minimal}.
  
All existing conformal measures with the possible exception of the minimal one
(in case $\de_*<2$),  are bound to be dissipative: 

\begin{lem}\label{diss2}
  Any conformal measure $\mu_\de$ on $J$ is dissipative for $\de>\de_*$. 
\end{lem}

\begin{pf}
  Assume $\mu_\de$ is conservative. Then by the last assertion of Lemma \ref{diss criterion},
for $\mu_\de$-a.e. $z$, there exists a sequence of moments $n_k\to \infty$ such that  
$\dist(f^{n_k} z, \OO) >2r$, where $r>0$ is independent of $z$. Let 
$D_k(z)$ be the component of $f^{-n_k} (\D_r ( f^{n_k} z))$ containing $z$. Then by the Koebe Distortion
Lemma,
$D_k(z)$ is an oval of bounded shape such and the map $f^{n_k}: D_k \ra \D_r( f^{n_k} z)$ has a bounded distortion
(with an absolute bound). 
By covariance of conformal measures, 
\begin{equation}\label{mu-de}
    \mu_\de(D_k)\asymp (\diam D_k)^\de \leq \eps_k\, (\diam D_k)^{\de_*} \asymp \eps_k\,  \mu_{\de_*} (D_k),    
\end{equation}
where $\eps_k\to 0$ as $k\to \infty$.
 
By the Besikovich Covering Lemma (see \cite[\S 2.7]{Mat}), 
there is a covering of a set $X$ of full $\mu_\de$-measure
with a family of balls $D_k$ with intersection multiplicity bounded by some absolute constant $N$.  
Hence 
$$
 \mu_\de(J)=\mu_\de(X) \leq \sum \mu_\de(D_k) \leq \eps_k\sum\mu_{\de^*}(D_k) \leq  N\eps_k\, \mu_{\de_*}(J).
$$ 
Letting $k\to \infty$ we conclude that $\mu_\de(J) =0 $  -- contradiction.
\end{pf}

As the minimal measure $\mu_*$ is concerned, 
both dissipative and conservative options are realizable, at least for quadratic-like maps. 
In the following statement we give a condition for dissipativity;
later on (\S \ref{conservativity}) we will deal with the conservative case. 

\begin{prop}\label{diss3}
Assume $\de_*(f) = \de_*(Rf) \equiv \de_*$.
Then the $\de_*$-conformal measure of $Rf$ is dissipative 
(and thus all conformal measures of $Rf$ are dissipative).  
\end{prop}

\begin{pf} Let $g: \u'\ra \V'$ be a pre-renormalization of $f$. We let $\Xi_*\equiv \Xi_{\de_*}$,
and let $\mu_*$ and $\nu_*$ be $\de_*$-conformal measures for $f$ and $g$ respectively.
We will use the measure $\mu_*$ to prove dissipativity of $\nu_*$.

Take some disk $D=\D_r(z)\subset \V'\sm \u'$ that intersects $J(f)$. Since $D$ is wandering under the dynamics of $g$,
$$
    \Xi_*(g,z) \asymp \sum\mu_*(g^{-n}(D)) <\infty.
$$
By Lemma \ref{diss criterion}, $\nu_*$ is dissipative. 
\end{pf}

 \begin{lemma}\label{diss4}

  Let $f$ be a Feigenbaum map which is recurrent under the dynamics of the
  renormalization operator (for instance, a renormalization fixed point). 
  Then the $\delta_*$-conformal measure of $f$ is
  dissipative.

 \end{lemma}

  \begin{pf}
  Note that Theorem \ref{DU formula} implies that $\de_*$ decays under renormalization:
$\de_*(Rf)\leq\de_*(f)$ (since it is obviously true for the hyperbolic dimension). 

  Let $R^{n_k} f \to f$, and let $g$ be a limit of $R^{n_k-1} f$. 
  Then $R g=f$, so that $ \de_*(g) \geq \de_* (f)$.
  On the other hand,  by the lower semicontinuity of $\delta_*$ (Lemma \ref{semi-cont}) we have
  $\delta_*(g) \leq \lim \delta_* (R^i(f)) \leq \delta_* (f)$.
  Thus  $\de_* (f) = \de_* (g)$, and  Proposition \ref{diss3} completes the proof.  
  \end{pf}

\subsection{Conformal spectrum}

  \begin{prop}\label{spec}
 For any $\de\geq \HD_\hyp(J)$, there exists a unique $\de$-conformal measure $\mu_\de$.
 In particular, there exists always a unique $2$-conformal measure $\mu_2$.
\end{prop}

\begin{pf}
The uniqueness part follows from Theorem \ref{ergodicity}, so let us deal with the existence. 

\msk
Let 
$ S_n^r   = \{ \zeta:\  f^n \zeta= 0 \text{ and }
|f^k \zeta| \geq r,\ k=0, \dots n-1 \}$ and let $S^r=\cup S^r_n$.  
% If $\zeta \in S^r$ we will let $n \equiv n(\zeta)$ be such that $\zeta \in S^r_n$.
  Let us consider the cut-off Poincar\'e series
$$
   \Xi^r_\de (0) =  \sum_{n=0}^\infty \sum_{\zeta\in S^r_n } {\frac {1}{|Df^n (\zeta)|^\de}}.
$$

   Since the set $J^r= \{z: |f^n z|\geq r,\ n=0,1,\dots\} $ is hyperbolic, 
$$ 
  \de_\crit(f|J^r) = \HD (J^r) \leq \HD_\hyp(J), 
$$  
%where the last equality is provided by the Denker-Urbanski Theorem (\ref{DU formula}).
It follows that the cut-off Poincar\'e series is convergent for $\de> \HD_\hyp(J)$.
Hence  we can weigh the involved preimages of $0$ by the  terms of $\Xi^r_\de(0)$
and consider the corresponding probability measure: 
$$
 \mu_\de^r = {\frac {1} {\Xi^r_\de (0)}} \sum_{n=0}^\infty  \sum_{\zeta\in S^r_n}
{\frac {\de_\zeta} {|Df^n (\zeta)|^\de}}.
$$

Let now $r\to 0$ ($\de>\de_*$ being fixed). 
Then $ \Xi^r_\de (0) \to \infty$. To see this,  consider the pre-renormalizations
$f_m : \u^m \ra \V^m$ and  pick  some $f_m$-preimage $\zeta_m$ of $0$. 
Since the derivatives $|Df_m (\zeta_m)|$ are bounded, 
$$
   \lim_{r\to 0} \Xi^r_\de (0)\geq  \sum_{m=0}^\infty {\frac {1}
{|Df_m (\zeta_m)|^\de}} = \infty, 
$$  
as was asserted.  

Let us now consider a weak$^*$ limit $\mu_\de$ of measures $\mu_\de^r$ as $r\to 0$.
We claim that this measure is $\de$-conformal. Let us first verify this outside the critical point.
Let $U$ be a small disk whose closure does not contain $0$ and such that $f|U$ is injective.
Then for sufficiently small $r>0$, 
$$
    f( S^r_n \cap U )= S^r_{n-1} \cap f(U),
$$
and hence 
$$
   \mu_\de^r |f(U) = |D f|^\de (\mu_\de^r | U).
$$
Passing to the limit,  we conclude that 
$$
   \mu_\de  |f(U) = |D f|^\de (\mu_\de | U),
$$
provided $\mu_\de (\di U)=0$. Since there is a base of disks with the latter property, the measure $\mu$
is $\de$-conformal everywhere outside $0$. 

But the critical value $f(0)$
cannot carry any mass, for otherwise each of the $f_m(0)$ would carry a definite
mass, and the total mass of $\mu_\de$ would be infinite.  (We also note that
the critical point $0$ cannot carry any mass either,
for otherwise every preimage $\zeta_m$ would carry a definite mass 
$|Df_m|^{-\de} \mu(0)$,
and the total mass of $\mu_\de$ would be infinite.)
This completes the proof.  
\end{pf} 

We will make use of the easy part of  Denker-Urbanski's Theorem:   
 
\begin{lem}\label{easy DU}
For any quadratic-like map $f$,  $\HD_\hyp(J)\leq \de_*$.
\end{lem}

\begin{pf}
  Let us consider a $\de$-conformal measure $\mu$ and a hyperbolic set $X\subset J$.
We need to show that $\HD(X)\leq \de$.

By a standard distortion argument using $\de$-covariance of $\mu$, we obtain:
 $\mu_\de(\D_r(z))\asymp r^\de$ for all disks $\D_r(z)$ centered at $z\in J$ of radius $r\leq 1$.

By the Besikovich Covering Lemma, for any $\eps>0$, 
there exist a  covering of $X$ with disks $\D_{r_i}(z_i)$ of radii $r_i< \eps$
with intersection multiplicity bounded by $N$. 
Then $$ 
  \sum r_i^\de \asymp \sum  \mu_\de(\D_{r_i}(z_i))^\de \leq N \mu_\de(X) . 
$$
Hence $h_\de(X)<\infty$, so that $\de\geq \HD(X)$.   
\end{pf}

\begin{thm}\label{two exponents}
   For any Feigenbaum map, $\HD_\hyp(J) = \de_* =  \de_\crit$.
\end{thm}

\begin{pf}
We will show that  $\HD_\hyp(J) \leq \de_* \leq   \de_\crit\leq \HD_\hyp(J)$.

The first inequality is the content of  Lemma \ref{easy DU}.

The middle inequality, $\de_*\leq \de_\crit$, holds  
since by Theorem \ref{Sullivan's meas} 
there exists a $\de_\crit$-conformal measure, and by definition,  $\de_*$ is the
minimum exponent of conformal measures on $J$.

To prove the last inequality, $ \de_\crit\leq \HD_\hyp(J)$, take any $\de>\HD_\hyp(J)$.
Then   by Lemma \ref{spec}, there exists a  $\de$-conformal measure on $J$.
By Lemma \ref{diss2}, it is dissipative.
By Lemma \ref{diss criterion},  the Poincar\'e series $\Xi_\de(z)$, $z\in \C\sm \OO$,  is convergent.
By definition,  $\de_\crit\leq \de$, and the conclusion follows.   
\end{pf}

\begin{rem}
Our exposition is organized in such a way that it does not use explicitly the constructive part of
Theorem \ref{DU formula} of Denker and Urbanski (which provides a minimal conformal measure $\mu_*$), 
so that along the lines we obtain a new proof of it (in the Feigenbaum setting). 
\end{rem}

\subsection{Scaling at the critical point}

Let $\mu$ be a $\de$-conformal measure on the Julia set $J(f)$, and let
$$
     \si = \lim_{r\to 0}\sup {\frac {\log \mu(\D_r)} {\log r}}.
$$

\begin{lem}\label{si leq de}
If  $\area(J(f))=0$ then $\si\leq\de $.
\end{lem}

\begin{pf}
 Let us consider a fundamental annulus $A^1=V^1\sm U^1$ of the pre-renormalization $g$ of $f$,
where $V^1$ and $U^1$ are domains constructed in \S \ref{nfd}.
Then $\mu(A^1)>0$ since $A^1 \cap J(f)\not=\emptyset$ and $\supp\mu = J(f)$.
We also have $\mu(\di A^1)=0$.\footnote {Indeed, we have $f^n(\partial A^1) \cap
U^1=\emptyset$, $n \geq 0$, which implies that $\mu(\partial A^1)=0$ in both
dissipative and conservative cases.}
%\footnote{Indeed,  by slightly perturbing the fundamental 
%annulus $\V^1\sm \u^1$, we can ensure that $\mu(\di \V^1\cup \di \u^1)=0$. But the construction
%of the domains $V^1$ and $U^1$ given in  Appendix A shows that their boundaries are contained
%in the iterated preimages of $\di \V^1\cup \di \u^1$.}
Given a point  $z\in U^1 \sm J(g) $, there is a unique moment $n(z)$
such that $g^{n (z) } \in A^1$. Since $\area(J(g))=0$, $n(z)$ is well defined for a.e. $z\in U^1$.

Let us tile $\overline A^1$ by finitely many compact sets $\Delta_i \subset
\overline A^1$ of positive $\mu$-measure that have boundary of zero
$\mu$-measure, such that each $\Delta_i$ is contained in a Jordan domain
$\tilde \Delta_i \subset \V^1  \setminus \OO(f)$.

%Let us tile $\bar A$ by several topological rectangles $\De_i\subset A$ of positive $\mu$-measure 
%that have %piecewise smooth 
%boundary of zero $\mu$-measure.
Pulling the $\Delta_i$ back by the appropriate branches of $g^{-n(z)}|\tilde
\Delta_i$, we obtain
compact sets $\Pi_j$ covering almost all of $U^1$. Since $U^1$ is a nice domain,
the $\Pi_j$
have pairwise disjoint interiors.  Moreover, there exists $\tilde \Pi_j
\Supset \Pi_j$ such that if $g^{n_j}(\Pi_j)=\Delta_i$ then $g^{n_j}|\tilde
\Pi_j$ is a univalent map onto $\tilde \Delta_i$
%Moreover, for each tile $\De_i$, there is 
%a bigger disk $\tl \De_i\Supset \De_i$ such that $\tl \De_i\cap \OO=\emptyset$.
%These disks can be
%univalently pulled back together with the $\De_i$  
%providing a definite Koebe space $\tl \Pi_j$ around the $\Pi_j$.  

Since the $\tilde \Pi_j$ provide a definite Koebe space around $\Pi_j$,
we can apply the Koebe Distortion Lemma to conclude that
each $\Pi_j$ is mapped with
bounded distortion by some $g^{n_j}$ onto the appropriate $\De_i$.
Together with the $\de$-covariance of $\mu$, this yields:
\begin{equation}\label{comparison}
   \mu(\Pi_j)\asymp {\frac {\mu(\De_i)} {|D g^{n_j} (z)|^\de}}\asymp
         (\diam \Pi_j)^\de.    %\asymp |\Pi_j|^{\de/2} .
\end{equation}
Take some $\eps>0$. Then there exist  arbitrary small disks $\D_r$ such that $\mu(\D_r) \leq r^{\si - \eps}$.
Since the hyperbolic diameter of $\Pi_j$ in $\tilde \Pi_j$ is bounded
and $0\not\in \tl \Pi_j$, we have: $\diam \Pi_j \leq C \dist(\Pi_j,0)$. Hence there is
a scaling factor $k>1$ independent of $r$ such that
if $\Pi_j\cap \D_{r/k} \not=\emptyset$ then $\Pi_j\subset \D_r$. 
Since the $\Pi_j$ have disjoint interiors and have boundaries of
zero $\mu$-measure, we conclude:

\begin{equation}\label{si}
    \sum \mu (\Pi_j) \leq \mu(\D_r)\leq  r^{\si-\eps}% = |\D_r|^{(\si-\eps)/2},
\end{equation}
where the summation is taken over those $\Pi_j$ that intersect $\D_{r/k}$.
Since these $\Pi_j$ cover almost all of $\D_{r/k}$,
$$
  r^2\asymp  \area (\D_{r/k}) \leq \sum \area \Pi_j \asymp \sum (\diam \Pi_j)^2.
$$
Hence  $r^\de \leq C \sum (\diam \Pi_j)^\de \asymp \sum \mu (\Pi_j) $.
Comparing it with (\ref{si}), we conclude that $\si-\eps \leq \de$.
Since it holds for any $\eps>0$, we are done.
\end{pf}

\subsection{Hyperbolic dimension}

\begin{prop}\label{local dimension}
Let $\mu$ be a $\de$-conformal measure on $J$.
  If $\area(J)=0$ then for any $z\in J$, $\underline {\HD}_\mu(z) \leq \de $.
% \note {I believe it can be changed to the upper dimension}
\end{prop}

\begin{pf}
By Lemma \ref{si leq de},  for any $\eps>0$ there exists a $c>0$ such that 
\begin{equation}\label{scaling at origin}
    \mu(\D_r) \geq c r^{\de+\eps}, \quad  r\leq \diam J.
\end{equation}

Take a point $z\in J$. If the $\orb(z)$ does not accumulate on 0,
then obviously $\HD_\mu(z) = \de$. If it does, then consider the first entries of the $\orb(z)$
 to the domains $V^m$ from \S \ref{nfd} and pull them back to $z$.
We obtain a shrinking nest of neighborhoods $D^m$ of $z$ with bounded shape.

Let $l_m$ be such that $f^{l_m}(D^m)=V^m$.  Then
it follows from the properties
(C2-3) of the domains $V^m$ that either
$f^{l_m}(D^{m+1}) \subset U^m$, or $f^{l_m}(D^{m+1}) \subset (V^m
\setminus U^m)$.  In either case, by property (G2) of the domains $V^m$,
$\area(f^{l_m}(D_{m+1}))/\area(f^{l_m}(D_m))$ is bounded away from $1$.
Since $f^{l_m}:D_m \to V^m$ has bounded distortion (it admits a univalent
extension onto $\V^m$), we conclude that $\area(D_{m+1})/\area(D_m)$ is
bounded away from $1$.  So $\diam(D_m)$ decays exponentially fast.

Since $\diam(V^m)>C \lambda^k$ for some $C>0$, $\lambda>0$, it follows
that $\diam(D^m) \leq c \diam(V^m)^\kappa$ for some $c>0$, $\kappa>0$.
It easily follows from (\ref{scaling at origin}), bounded distortion,
and $\de$-covariance of $\mu$ that 
$\mu(D^m) \geq c (\diam D^m)^{\de+\eps}$, and the conclusion follows.
\end{pf}

 \begin{cor}\label{HD = hyp}
If $\area(J)=0$ then $\HD(J)= \HD_\hyp (J)$. 
  \end{cor}

\begin{pf} 
By Theorem \ref{DU formula}, it is enough to prove that 
for any $\de$-conformal measure $\mu$ on $J$, $\HD(J) \leq \de $. Take some $\eps>0$.
By Proposition \ref{local dimension}, any point $z\in J$ is a center of arbitrary small disks
$\D_r(z)$ such that $\mu(\D_r(z)) \geq r^{\de+\eps}$. 
By the Besicovich Covering Lemma, $J$ admits an arbitrary fine covering by such disks
$D_k=\D(z_k, r_k)$ with intersection multiplicity at most $N$. 
Hence $\sum r_k^{\de+\eps} \leq N \mu(J) $, $h_\de(J)<\infty$, and we are done. 
\end{pf} 

\begin{rem}

This result also follows from 
Corollary \ref {two exponents} and the work of Bishop \cite {B}.

\end{rem}

\subsection{All dimensions and exponents are equal}

Combining  Theorem \ref{two exponents} with Corollary \ref{HD = hyp}, we obtain: 

\begin{thm}
Assume that a Feigenbaum Julia set $J$ has zero area. Then
$$
    \de_\crit  = \de_* = \HD_\hyp(J)  =  \HD( J ). 
$$
\end{thm}

\subsection{Continuity properties of the critical exponent}
  Let us formulate a couple of useful technical lemmas. 
Theorem \ref{two exponents} and  Lemma \ref{semi-cont} yield:  

\begin{cor}\label{lower}
  The critical exponent $\de_\crit(f)$ depends lower semicontinuously on the
Feigenbaum map $f$.
\end{cor}

The situation within one hybrid class is even better: 

\begin{lem}\label{continuity1}
  The critical exponent $\de_\crit(g)$ is continuous on any hybrid class $\HH_f$.
Moreover, it is uniformly continuous over compact subsets $K$ of quadratic-like germs:
For any $\eps>0$ there exists $\si >0$ such that if $f,g \in K$, $g\in \HH_f$ and  $\dist(f,g)< \si$
then $|\de_\crit(f) - \de_\crit(g)| <\eps$. 
 
\end{lem}

\begin{pf}
If $f,g \in K$, $g\in \HH_f$ and  $\dist(f,g)<\si $ then $f$ and $g$ are 
conjugate by a  $K$-qc homeomorphism $h$, where
dilatation  $K$ is close to 1. This map is   $\kappa$-H\"older continuous,
with exponent $\kappa=1/K$.

%  It follows from the fact the the critical exponent does not change much under
% $(1+\eps)$-H\"older continuous conjugacies.

Let $f:\u\ra \V$ and  $\OO=\OO(f)$.
Consider a round disk $D$ compactly contained in $\V \setminus \OO$ centered around some point $x$.
If $f^m(y)=x$, let $D_y$ be the component of $f^{-m}(x)$ containing $y$.  Then by the Koebe Distortion Lemma,
the maps $f^m : D_y \to D$ and $g^m: h(D_y) \to h(D)$ have bounded distortion.  Thus 
$|Df^m(y)| \asymp (\diam D_y) ^{-1} $ and $|Dg^m (h(y))| \asymp (\diam h (D_y))^{-1} $.  This implies that
$$
   M^{-1} |Df^m (y)|^\kappa \leq |Dg^m_p(y)| \leq M |Df^m_p(y)|^K.
$$
Hence
$K^{-1} \leq \delta_\crit (f_p) / \delta_\crit (g_p)  \leq K $, and we are done.

\end{pf}

\section{Recursive Quadratic Estimates}\label{pert est}

Let $X_{m,n}$, $m<n$ be the set of points in $U^m$ that land in $V^n$
under iterates of $f_m$, and let $X_n=X_{0,n}$.  Let $Y_{m,n}$, $m<n$ be
the set of points in $A^n$ that never return to $V^n$ under iterates of
$f_m$, and let $Y_n=Y_{0,n}$.
%Let $Z_{n,m}$ be the set of points in $A^m$ which, under iteration by $f$,
%either do not return to $A^m$ or only return after passing by $A^{n+m}$.

We define the following quantities:
\be
\eta_{m,n} \equiv \eta_{m,n}(f)=\frac {|X_{m,n}|} {|U^m|},
\quad \eta_n=\eta_{0,n},
\ee
\be
\xi_{m,n} \equiv \xi_{m,n}(f)=\frac {|Y_{m,n}|} {|A^n|},
\quad \xi_n=\xi_{0,n},
\ee
%\be
%\kappa_{n,m} \equiv \kappa_{n,m}(f)=\frac {|Z_{n,m}|} {|A^m|},
%\quad \kappa_n=\kappa_{n,0},
%\ee
\be
\rho_{m,n} \equiv \rho_{m,n}(f)=\frac {|U^n|} {|U^m|}.
\ee

Let us start with a couple of relations between these quantities:

\begin{lemma} \label {exp}

Let $f$ be a Feigenbaum map.  There exists $C>0$ only depending on
the geometric bounds and $C_0>0$, depending also on the combinatorial bounds,
such that
%\be
%C^{-1}(\xi_m+\eta_{n,m}) \leq \kappa_{n,m} \leq C(\xi_m+\eta_{n,m}),
%\ee
\be
\max \{C_0^{-1},1-C_0\xi_{m,n}\}
\leq \frac {\eta_{m,n+1}} {\eta_{m,n}} \leq 1-C^{-1}\xi_{m,n}.
\ee
%\be
%\rho_{m,n} \leq \theta^{n-m}.
%\ee

\end{lemma}

\begin{pf}

We may assume that $m=0$.
The first landing map $\phi:X_n \to V^n$ sends each component of
$X_n$ univalently onto $V^n$ with bounded distortion.  Clearly,
$\phi^{-1}(Y_n) \cap X_{n+1}=\emptyset$, so that
\be
\frac {\eta_{n+1}} {\eta_n}=p(X_{n+1}|X_n) \leq 1-C^{-1} \xi_n.
\ee
Notice that $|V^{n+1}|/|V^n|>q$, where $q$ depends also on the
combinatorial bounds, so we get
\be
\frac {\eta_{n+1}} {\eta_n}>C_0^{-1}.
\ee

To conclude the remaining inequality, it is enough to show that
$p(X_{n+1}|V^n)>1-C_0 \xi_n$.  Let $S^n$ be the set of
points in $V^n$ which land in $A^n$ and later return to $V^n$, and let $W^n$
be the set of connected components of $S^n$ which are not contained in
$V^{n+1}$.  Then $|W^n|/|V^n|<1-q$, where $q$ depends also on the
combinatorial bounds. We set $\psi:W^n \to V^n$ to be the composition of the
first landing to $A^n$ followed by the first return to $V^n$.  The
restriction of $\psi$ to each connected component of $W^n$ is a univalent
map onto $V^n$, moreover any iterate of $\psi$ has bounded distortion by
the Koebe Distortion Lemma (the restriction of $\psi$ to each component of
$W^n$ extends to a univalent map onto $\V^n$).
Since $|W^n|/|V^n|<1-q$, we have $|\psi^{-k}(W^n)| \leq
e^{-q k} |V^n|$.

Since almost every $x \in V^n$ either escapes through $Y_n$ or eventually
passes through $V^{n+1}$, we have
\be
\frac {|V^n \setminus X_{n+1}|} {|V^n|} \leq
\sum_{k=0}^\infty \frac {|\psi^{-k}(Y_n)|} {|V^n|} \leq
C \sum_{k=0}^\infty e^{-q k} \xi_n \leq C_0 \xi_n.
\ee
\end{pf}

%Recall that
%\be
%\Xi_\delta(f_m,x)=\sum_{k \geq 0} \sum_{f_m^k(y)=x}
%\frac {1} {|Df_m^k(y)|^\delta}.
%\ee
Let   $\Xi^{[j]}$  be the Poincar\'e series
$\Xi$   truncated at level $j$, 
so that  summation in (\ref{Poincare series}) is taken over  $n \leq j$ instead of all $n$.
Let
\be
\omega_{m,n}(f,\delta) \equiv \omega_{m,n}(\delta)=
\frac {1} {|A^n|} \int_{A^n} \Xi_\delta(f_m,x) dx,
%\sum_k \sum_{f^k_m(y)=x} \frac {|Df^k_m(y)|}^{-\delta} dx
\quad \omega_n=\omega_{0,n}.
\ee
We will use $\omega^{[j]}$ to denote the truncation of
$\omega$ at level $j$, that is, taking the
integral of $\Xi^{[j]}$ instead of $\Xi$.

Our key estimates are presented in Lemmas~\ref {pert1} and~\ref {pert2}.

\begin{lemma} \label {pert1}

Let $f$ be a Feigenbaum map.  Then for any $0 \leq l<m<n$,
there exists $C>0$, only depending on the geometric bounds,
and $\epsilon>0$, depending on the geometric and combinatorial bounds,
and on $n-l$, $m-l$, such that if $2-\epsilon \leq \delta \leq 2$ then
\be \label {301}
\omega_{l,m}^{[j+1]}(\delta) \leq C \frac {\eta_{l,m}}
{\rho_{l,m}}+(1-C^{-1}(\eta_{m,n}+\xi_{l,m}))
\omega_{l,m}^{[j]}(\delta)+C \rho_{m,n} \xi_{m,n}
\omega_{l,n}^{[j]}(\delta),
\ee
\be \label {302}
\omega_{l,n}^{[j+1]}(\delta) \leq C \left (\frac {\eta_{l,m}}
{\rho_{l,m}}+\omega_{l,m}^{[j]}(\delta) \right )\omega_{m,n}^{[j]}(\delta).
\ee
\comm{
\be
\omega_{n,r}^{[l+1]}(\delta) \leq
C \frac {\eta_{n,r}} {\rho_{n,r}}+(1-C^{-1}(\eta_{n,r}+\xi_{n,r})
\omega_{n,r}^{[l]}(\delta)+C \eta_{n,r} \xi_{m,n+r}
\frac {\rho_{m,n+r}} {\rho_{n,r}} \omega_{m,n+r}^{[l]}(\delta)+
C \xi_{m,n+r} \rho_{m,n+r} \omega_{n,r}^{[l]}(\delta)
\omega_{m,n+r}^{[l]}(\delta).
\ee
}

\end{lemma}

Let $\tau_{m,n}$, $m<n$ be the infimum of $|Df^k_m(x)|$
over all $x \in U^m$, $k \geq 0$ such that
$f^k_m(x) \in V^n$ and $f^j_m(x) \notin V^n$,
$0 \leq j \leq k$, and let $\tau_n=\tau_{0,n}$.
Let $\upsilon_{m,n}$, $m<n$ be the infimum of
$|Df^k_m(x)|$ over all $x \in A^n$, $k \geq 1$ such that
$f^k_m(x) \in A^m$ and $f^j_m(x) \notin A^m \cup A^n$,
$1 \leq j \leq k$, and let $\upsilon_n=\upsilon_{0,n}$.
Let us show that the product $\tau \upsilon$ is exponentially large
in terms of $n-m$.

\begin{lemma} \label {tauupsilon}

Let $f$ be a Feigenbaum map.  Then for every $0 \leq m<n$, there exist $C>0$,
$\theta<1$, only depending on the geometric bounds, such that for every
$m<n$ we have
\be
\tau_{m,n} \upsilon_{m,n} \geq C^{-1} \theta^{m-n}.
\ee

\end{lemma}

\begin{pf}

We may assume that $m=0$.  Let $f^k(x)$ be the first landing of $x$ in
$V^n$.  Then there exists a domain $D \subset U$ containing $x$ such that
$f^k:D \to V^n$ is
univalent and has bounded distortion.  Let now $z \in A^n$ and let 
$y=(f^k|D)^{-1}(z)$.  We claim that if $f^s(z) \in A$ then 
$|Df^{k+s}(y)| \geq C^{-1} \theta^{-n}$, which implies the 
result (notice that $|Df^k(x)|\, |Df^s(z)| \geq C^{-1} |Df^{k+s}(y)|$ by
bounded distortion of $f^k|D$).

To see this, let $W$ be the complement of the postcritical set in some 
definite neighborhood of $V$, and consider the Poincar\'e metric $\| 
\cdot \|_W$ on $W$.  Since $A$ is away from the postcritical set, we have
$$
|Df^{k+s}(y)| \geq C^{-1} \|Df^{k+s}(y)\|_W.
$$
It follows from the Schwarz Lemma that 
$f:U \to V$ does not contract $\| \cdot \|_W$, so
$$
\|Df(f^j(y))\|_W \geq 1, \quad 0 \leq j \leq k+s-1.
$$
Closer inspection, see \cite[Prop. 6.9]{McM-towers}
%\marginpar {Give more precise reference.}
shows that whenever $f^j(y) \in A^i$ with $1 \leq i \leq n$
we have  $\|Df(f^j(y))\|_W \geq \theta^{-1}$.  Thus
$$
\|Df^{k+s}(y)\|_W \geq \theta^{-n}.
$$
\end{pf}

\begin{lemma} \label {pert2}

Let $f$ be a Feigenbaum map.  Then for every $0
\leq l<m<n$, there exist $C>0$,
only depending on the geometric bounds,
such that if $\delta \leq 2$ then
\be \label {303}
\omega_{l,m}(\delta) \geq C^{-1} \tau^{2-\delta}_{l,m}
\frac {\eta_{l,m}} {\rho_{l,m}}+\max \{1-C (\eta_{m,n}+\xi_{l,m}),0\}
\omega_{l,m}(\delta)+C^{-1} \upsilon^{2-\delta}_{m,n}
\xi_{m,n} \rho_{m,n} \omega_{l,n}(\delta),
\ee
\be \label {304}
\omega_{l,n}(\delta) \geq C^{-1} \left (\tau^{2-\delta}_{l,m}
\frac {\eta_{l,m}} {\rho_{l,m}}+\omega_{l,m}(\delta) \right )
\omega_{m,n}(\delta).
\ee
\comm{
\be
\omega_{n,r}(\delta) \geq
C^{-1} \frac {\eta_{n,r}} {\rho_{n,r}}+(1-C(\eta_{n,r}+\xi_{n,r})
\omega_{n,r}(\delta)+C^{-1} \theta^{\delta-2} \eta_{n,r} \xi_{m,n+r}
\frac {\rho_{m,n+r}} {\rho_{n,r}} \omega_{m,n+r}(\delta)+
C^{-1} \theta^{\delta-2} \xi_{m,n+r} \rho_{m,n+r} \omega_{n,r}(\delta)
\omega_{m,n+r}(\delta).
\ee
}

\end{lemma}

\subsection{Proof of Lemmas \ref {pert1} and \ref {pert2}}
\label{proofofestimates}

We may assume that $l=0$.  To simplify notation, truncation of the
Poincar\'e series (labeled by $j$ or $j+1$) will be implicit in what
follows.

We define $Z_{m,n}$ as the set of points in $A^m$ that do not return to
$A^m$ before passing through $A^n$.  We let
\be
\kappa_{m,n}=\frac {|Z_{m,n}|} {|A^m|}.
\ee
Notice that
\be
\kappa_{m,n} \asymp \xi_m+\eta_{m,n}
\ee

In this section, an orbit (under iteration of $f_r$) will denote a sequence
$(x_0,...,x_k)$, $k \geq 0$, where $x_k \in V^r$ and
$f_r(x_i)=x_{i+1}$, $i<k$.  Unless explicitly stated otherwise, the orbits
considered will be orbits under iteration of $f$.

An orbit will be said to be trivial if its length is $0$.

Given a family $\FF$ of orbits $(x_0,...,x_k)$ under iteration of some $f_r$,
we define a function $\C \to [0,\infty]$
\be
\Xi_\delta(\FF)(z)=\sum_{(x_0,...,x_k=z) \in \FF} |Df^k_r(x_0)|^{-\delta}
\ee
(notice that we do not make $f_r$ explicit in this notation).

We shall need some special notation for certain families of orbits.  By $D
\xleftarrow{} E$, we will understand the family of orbits
$(x_0,...,x_k)$ with $x_0 \in E$ and $x_k \in D$.
The family of orbits $(x_0,...,x_k)$ with $x_0
\in E$, $x_k \in D$ and $x_1,...,x_{k-1} \in S$ will be denoted $D
\xleftarrow[S]{} E$.  A ``plus sign'' over the arrow will indicate that
only non-trivial orbits are considered.
The justaposition of arrows will denote composition in the
natural way.  Thus
\be
D \xleftarrow[S]{+} D \xleftarrow[S]{} E,
\ee
denotes the family of orbits $(x_0,...,x_k)$, with $x_0 \in E$, $x_k \in D$,
and such that there exists $0 \leq i<k$ such that $x_i \in D$ and
$x_1,...,x_{i-1},x_{i+1},...,x_{k-1} \in S$.

When considering orbits under some $f_r$ other than
$f$, this will be indicated above the arrow, for instance, $D
\xleftarrow[S]{f_r} E$ or $D \xleftarrow{f_r,+} E$.

Given a function $\phi:\C \to [0,\infty]$ and a measurable set $N \subset
\C$, we let
\be
I(\phi|N)=\int_N \phi(x) dx
\ee
and if $|N|>0$, we let
\be
E(\phi|N)=\frac {1} {|N|} I(\phi|N).
\ee
More generally, if $\mu$ is a measure on $\C$ with $\mu(\V)<\infty$,
we let
\be
E_\mu(\phi|N)=\frac {1} {\mu(N)} \int_N \phi(x) d\mu(x)
\ee
whenever $\mu(N)>0$.

In this notation, we have for instance
\be
\omega_m(\delta)=E(\Xi_\delta(A^m \leftarrow U)|A^m).
\ee

Let us say that a family $\FF$  of orbits is hyperbolic if $(x_0,...,x_k)
\in \FF$ implies that $|Df^k(x_0)| \geq C \lambda^k$, $C>0$, $\lambda>1$.
Notice that $\FF$ is hyperbolic if and only if there exists a neighborhood
$W$ of $0$ such that if $(x_0,...,x_k) \in \FF$ then $x_0,...,x_{k-1} \notin
W$.

Let $\delta_\FF$ be the infimum over all $\delta$ such that
$\Xi_\delta(\FF)$ is uniformly bounded.  By Corollary~\ref {Delta},
$\delta_\FF<2$ if $\FF$ is hyperbolic.

\comm{
\begin{lemma} \label {finiteness}

Let $\FF$ be a hyperbolic family of branches.  Then there exists
$\delta_\FF<\delta_\crit$ such that if $\delta \geq \delta_\FF$ then
$\Xi_\delta(\FF)$ is uniformly bounded.

\end{lemma}

\begin{pf}

It is enough to prove the result for $\FF$ of the form $U
\setminus W \xleftarrow[U \setminus W]{} U \setminus W$
for a neighborhood $W$ of $0$.  But in this case the result follows from
Lemma~\ref {continuity}.
\end{pf}
}

\begin{lemma} \label {nu_delta}

%Let $\FF$, $\delta_\FF$ be as in the previous lemma.
%Assume that
%$\Xi_\delta(\FF)$ is uniformly bounded for $\delta$ in some open
%interval $\Delta$.
With $\delta_\FF$ defined as above, we have:
\begin{enumerate}
\item
The map $\delta \mapsto \sup \Xi_\delta(\FF)$ is a continuous convex
function of $\delta>\delta_\FF$,
\item
Let $\mu$ be a finite measure on $\V$, and let
\be
\nu_\delta=\int \sum_{(x_0,...,x_k=x) \in \FF} \frac {\delta_{x_0}}
{|Df^k(x_0)|^\delta} d\mu(x)
\ee
which is a finite measure for $\delta>\delta_\FF$.
Then for any Borelian set $N$,
$\delta \mapsto \nu_\delta(N)$ is a continuous convex function
of $\delta>\delta_\FF$.  In particular $\delta \mapsto \nu_\delta$,
$\delta>\delta_\FF$ is
continuous in the weak$^*$ topology.
\end{enumerate}

\end{lemma}

\begin{pf}

Follows immediately from the convexity of $\delta \mapsto x^{-\delta}$.
%Use Lemma~\ref {finiteness} and the convexity in $\delta$ of
%$\Xi_\delta(\FF)$ for $\delta \geq \delta_\FF$.
\end{pf}

In the most common application of the previous lemma, we will take $\mu$ as
Lebesgue and $\FF$ as a set of hyperbolic first returns or landings, say
$A^m \xleftarrow[U \setminus V^n]{+} A^m$.  In this case, the measure
$\nu_2$ is the restriction of Lebesgue measure to the set of points that
land or return.

\subsubsection{Proof of (\ref {301})}

Let
\be
B \equiv B^{m,n}=U \setminus (A^m \cup A^n).
\ee

We clearly have
\be \label {221}
\omega_m(\delta)=
E(\Xi_\delta(A^m \leftarrow U)|A^m)=1+E(\Xi_\delta(A^m \xleftarrow{+}
U)|A^m),
\ee
where the $1$ accounts for the trivial orbits.

Since $U$ is the disjoint union of $A^n$, $A^m$ and $B$, we can split
$A^m \xleftarrow{+} U$,
the set of non-trivial orbits that terminate in $A^m$ (under iteration by
$f$) in three groups:
\begin{enumerate}
\item $A^m \xleftarrow[B]{} B$ which consists of orbits
$(x_0,...,x_k)$ terminating in $A^m$
such that $x_0,...,x_{k-1}$ do not intersect $A^m \cup A^n$,
\item $A^m \xleftarrow[B]{} A^n \leftarrow
U$ which consists of orbits $(x_0,...,x_k)$ such that
$x_0,...,x_{k-1}$ intersects $A^n \cup A^m$ and the last such
intersection happens in $A^n$,
\item $A^m \xleftarrow[B]{+} A^m \leftarrow
U$ which consists of non-trivial orbits $(x_0,...,x_k)$ such that
$x_0,...,x_{k-1}$ intersects $A^m \cup A^n$ and the last
such intersection happens in $A^m$.
\end{enumerate}

Together with (\ref {221}), this decomposition gives
\be \label {81}
\omega_m(\delta)=1+E(\Xi_\delta(A^m \xleftarrow[B]{} B)|A^m)+
E(\Xi_\delta(A^m \xleftarrow[B]{} A^n
\leftarrow U)|A^m)+E(\Xi_\delta(A^m \xleftarrow[B]{+} A^m \leftarrow
U)|A^m).
\ee

Notice that
\be
I(\Xi_2(A^m \xleftarrow[B]{} B)|A^m),
\ee
is the area of the set of points in $B$
that land in $A^m$ before passing through $V^n$, which is
certainly smaller than $\eta_m |U|$.  Thus we have
\be \label {401}
E(\Xi_2(A^m \xleftarrow[B]{} B)|A^m) \leq \eta_m \frac {|U|}
{|A^m|}=\frac {\eta_m} {\rho_m} \frac {|U^m|} {|A^m|}.
\ee
By hyperbolicity of the set of orbits that avoid $V^n$, we have
\be \label {71}
E(\Xi_\delta(A^m \xleftarrow[B]{} B)|A^m)
\leq C \frac {\eta_m} {\rho_m}, \quad \delta \approx 2,
\ee
by continuity of the left hand side in $\delta \approx 2$.

We have the obvious estimate
\be \label {51}
E(\Xi_\delta(A^m \xleftarrow[B]{} A^n
\leftarrow U)|A^m) \leq
E(\Xi_\delta(A^m \xleftarrow[B]{} A^n)|A^m)
\sup_{A^n} \Xi_\delta(A^n \leftarrow U),
\ee
By bounded oscillation of the Poincar\'e series in $A^n$, we have
\be \label {52}
\inf_{A^n}\, \Xi_\delta(A^n \leftarrow U) \asymp
\sup_{A^n}\, \Xi_\delta(A^n \leftarrow U) \asymp E(\Xi_\delta(A^n
\leftarrow U)|A^n)=\omega_n(\delta).
\ee
Moreover,
\be
I(\Xi_2(A^m \xleftarrow[B]{} A^n)|A^m)
\ee
is the area of the set of points in $A^n$ that land in $A^m$ before
passing through $V^n$.  This is precisely $|A^n|\, \xi_{m,n}$, so
\be \label {402}
E(\Xi_2(A^m \xleftarrow[B]{}
A^n)|A^m)=\xi_{m,n} \frac {|A^n|} {|A^m|}=\xi_{m,n} \rho_{m,n} \frac
{|U^m|} {|A^m|} \frac {|A^n|} {|U^n|} \leq C \xi_{m,n} \rho_{m,n}.
\ee
By hyperbolicity of the set of orbits avoiding $U^n$, we get
\be \label {53}
E(\Xi_\delta(A^m \xleftarrow[B]{} A^n)|A^m)
\leq C \xi_{m,n} \rho_{m,n}, \quad \delta \approx 2.
\ee
Estimates (\ref {51}), (\ref {52}) and (\ref {53}) imply
\be \label {55}
E(\Xi_\delta(A^m \xleftarrow[B]{} A^n
\leftarrow U)|A^m) \leq C \xi_{m,n} \rho_{m,n} \omega_n(\delta), \quad
\delta \approx 2.
\ee

Let us estimate the last term of (\ref {81}).  We have the identity
\be \label {92}
\frac {E(\Xi_2(A^m \xleftarrow[B]{+} A^m
\leftarrow U)|A^m)} {E(\Xi_2(A^m \leftarrow U)|A^m)}=
\frac {E(\Xi_2(A^m \setminus Z_{m,n} \leftarrow U)|A^m)}
{E(\Xi_2(A^m \leftarrow U)|A^m)}=
1-\frac {E(\Xi_2(Z_{m,n} \leftarrow U)|A^m)} {E(\Xi_2(A^m \leftarrow U)|A^m)}.
\ee
By bounded oscillation of Poincar\'e series, we have
\be \label {93}
\frac {E(\Xi_\delta(Z_{m,n} \leftarrow U)|A^m)}
{E(\Xi_\delta(A^m \leftarrow U)|A^m)} \asymp \kappa_{m,n}.
\ee
By hyperbolicity of the set of orbits that do not pass through $V^n$, we
have
\be \label {91}
E(\Xi_\delta(A^m \xleftarrow[B]{+} A^m
\leftarrow U)|A^m) \leq (1-C^{-1} \kappa_{m,n}) E(\Xi_\delta(A^m \leftarrow
U)|A^m)=(1-C^{-1} \kappa_{m,n}) \omega_m(\delta), \quad \delta \approx 2.
\ee

Plugging (\ref {71}), (\ref {55}) and (\ref {91}) in (\ref {81}),
we get
\be
\omega_m(\delta)=1+C \frac {\eta_m} {\rho_m}+(1-C^{-1} \kappa_{m,n})
\omega_m(\delta)+C \xi_{m,n} \rho_{m,n} \omega_n(\delta), \quad \delta
\approx 2
\ee
which implies (\ref {301}) (notice that $\eta_m \geq \rho_m$, so the
$1$ can be put inside $C \eta_m/\rho_m$).

\subsubsection{Proof of (\ref {302})}

Let us split $A^n \leftarrow U$ into two groups of orbits:
\begin{enumerate}
\item $A^n \xleftarrow[U \setminus A^m]{} U \setminus A^m$,
\item $A^n \xleftarrow[U \setminus A^m]{} A^m \leftarrow U$.
\end{enumerate}

We further split $A^n \xleftarrow[U \setminus A^m]{} U \setminus A^m$
into two groups:
\begin{enumerate}
\item[(1a)] $A^n \xleftarrow[U \setminus A^m]{} U^m$,
\item[(1b)] $A^n \xleftarrow[U \setminus A^m]{} U^m \xleftarrow[U \setminus
V^m]{} U \setminus V^m$.
\end{enumerate}

Notice that we can write
\be
\left (A^n \xleftarrow[U \setminus A^m]{} A^m \right )=\left (A^n
\xleftarrow[U \setminus A^m]{} U^m \xleftarrow[U \setminus V^m]{} A^m
\right ).
\ee

Thus we have
\begin{align} \label {115}
\omega_n(\delta)=E(\Xi_\delta(A^n \leftarrow U)|A^n)=
&E(\Xi_\delta(A^n \xleftarrow[U \setminus A^m]{} U^m)|A^n)\\
\nonumber
&+E(\Xi_\delta(A^n \xleftarrow[U \setminus A^m]{} U^m
\xleftarrow[U \setminus V^m]{} U \setminus V^m)|A^n)\\
\nonumber
&+E(\Xi_\delta(A^n\xleftarrow[U \setminus A^m]{} U^m
\xleftarrow[U \setminus V^m]{} A^m \leftarrow U)|A^n).
\end{align}

This gives the estimate
\be \label {111}
\frac {E(\Xi_\delta(A^n \leftarrow U)|A^n)}
{E(\Xi_\delta(A^n \xleftarrow[U \setminus A^m]{} U^m)|A^n)} \leq
1+\sup_{U^m}\,
\Xi_\delta(U^m \xleftarrow[U \setminus V^m]{} U \setminus V^m)+
\sup_{U^m}\, \Xi_\delta(U^m \xleftarrow[U \setminus V^m]{} A^m
\leftarrow U).
\ee

Since $f_m$ is the first return map from $U^m$ to $V^m$,
there is a natural correspondence between the orbits of
$A^n \xleftarrow[U \setminus A^m]{}
U^m$ and $A^n \xleftarrow{f_m} U^m$.  It follows that
\be
\Xi_\delta(A^n \xleftarrow[U \setminus A^m]{} U^m)=\Xi_\delta(A^n
\xleftarrow{f_m} U^m),
\ee
thus
\be \label {191}
E(\Xi_\delta(A^n \xleftarrow[U \setminus A^m]{} U^m)|A^n)=
E(\Xi_\delta(A^n \xleftarrow{f_m} U^m)|A^n)=\omega_{m,n}(\delta).
\ee

By the Markov property\footnote{Let us remark that one does not need to use
the Markov property to conclude the upper bound in estimates (\ref {801})
and (\ref {101}).} of the first landing map to $V^m$, we have
\be \label {801}
\sup_{U^m} \Xi_2(U^m \xleftarrow[U \setminus V^m]{} U \setminus V^m)
\asymp
\inf_{U^m} \Xi_2(U^m \xleftarrow[U \setminus V^m]{} U \setminus V^m)
\asymp
E(\Xi_2(U^m \xleftarrow[U \setminus V^m]{} U \setminus V^m)|U^m) \asymp
\frac {\eta_m} {\rho_m}.
\ee
By hyperbolicity of the set of orbits avoiding $U^m$,
\be \label {121}
\sup_{U^m} \Xi_\delta(U^m \xleftarrow[U \setminus V^m]{} U \setminus V^m)
\leq C \frac {\eta_m} {\rho_m}, \quad \delta \approx 2.
\ee

Similarly,
\be \label {101}
\sup_{U^m} \Xi_2(U^m \xleftarrow[U \setminus V^m]{} A^m) \asymp
\inf_{U^m} \Xi_2(U^m \xleftarrow[U \setminus V^m]{} A^m) \asymp
E(\Xi_2(U^m \xleftarrow[U \setminus V^m]{} A^m)|U^m) \asymp 1,
\ee
and by hyperbolicity of the set of orbits avoiding $U^m$, we get
\be \label {102}
E(\Xi_\delta(U^m \xleftarrow[U \setminus U^m]{} A^m)|U^m) \leq C, \quad
\delta \approx 2.
\ee

\comm{
Notice that any pullback of $U^m$ along an orbit $(x_0,...,x_k)$
of $U^m \xleftarrow[U \setminus V^m]{} U \setminus V^m$ is univalent and
can be extended to a univalent pullback of $V^m$.  Thus, for every $z \in
U^m$, if $z \in D$ is a disk of diameter $C^{-1}$
times the diameter of $U^m$, any pullback of $D$ along an orbit
$(x_0,...,x_k=z)$ of $U^m \xleftarrow[U \setminus V^m]{} U \setminus V^m$ is
univalent with bounded distortion, does not intersect $U^m$, and moreover
such pullbacks form a disjoint family.  It follows that
\be
\sup_{U^m} \Xi_\delta(U^m \xleftarrow[U \setminus V^m]{} U \setminus V^m)
\leq C \frac {\eta_m} {\rho_m}.
\ee
By hyperbolicity of the set of orbits avoiding $U^m$,
\be \label {121}
\sup_{U^m} \Xi_\delta(U^m \xleftarrow[U \setminus V^m]{} U \setminus V^m)
\leq C \frac {\eta_m} {\rho_m}, \quad \delta \approx 2.
\ee
}

\comm{
Notice that any pullback of $U^m$ along an orbit of $U^m \xleftarrow[U
\setminus V^m]{} U \setminus V^m$ is univalent and can be extended to a
univalent pullback of $V^m$.  It follows that such univalent pullback has
bounded distortion and is contained in $U \setminus U^m$.  Thus all such
pullbacks are disjoint.  This implies that
\be \label {801}
\sup_{U^m} \Xi_\delta(U^m \xleftarrow[U \setminus V^m]{} U \setminus V^m)
\leq C E(\Xi_\delta(U^m \xleftarrow[U \setminus U^m]{} U \setminus
U^m)|U^m).
\ee
Moreover,
\be
I(\Xi_2(U^m \xleftarrow[U \setminus U^m]{} U \setminus
U^m)|U^m)
\ee
is the area of the set of points in $U \setminus U^m$ that land in $U^m$
which is smaller than $\eta_m |U|$.
Thus
\be \label {121}
\sup_{U^m} \Xi_\delta(U^m \xleftarrow[U \setminus V^m]{} U \setminus V^m)
\leq C \eta_m \frac {|U|} {|U^m|}=C \frac {\eta_m} {\rho_m}, \quad \delta
\approx 2.
\ee

Similarly,
\be \label {101}
\sup_{U^m} \Xi_2(U^m \xleftarrow[U \setminus V^m]{} A^m) \leq C,
\ee
and by hyperbolicity of the set of orbits avoiding $U^m$, we get
\be \label {102}
E(\Xi_\delta(U^m \xleftarrow[U \setminus U^m]{} A^m)|U^m) \leq C, \quad
\delta \approx 2.
\ee
}

\comm{
Similarly, any pullback of $U^m$ along an orbit of
$U^m \xleftarrow[U \setminus  V^m]{} A^m$ is univalent, has bounded
distortion and is contained in a neighborhood of $A^m$ of size comparable to
the diameter of $U^m$.  The family of such pullbacks is not necessarily
disjoint, but has intersection multiplicity $2$, that is, each
$z \in \C$ belongs to at most two such pullbacks.   It follows that
\be \label {101}
\sup_{U^m} \Xi_2(U^m \xleftarrow[U \setminus V^m]{} A^m) \leq C,
\ee
and by hyperbolicity of the set of orbits avoiding $U^m$, we get
\be \label {102}
E(\Xi_\delta(U^m \xleftarrow[U \setminus U^m]{} A^m)|U^m) \leq C, \quad
\delta \approx 2.
\ee
}
Estimate (\ref {102}) implies
\be
\sup_{U^m} \Xi_\delta(U^m \xleftarrow[U \setminus V^m]{} A^m \leftarrow U)
\leq C \sup_{A^m} \Xi_\delta(A^m \leftarrow U), \quad \delta \approx 2,
\ee
and by bounded oscillation of the Poincar\'e series in $A^m$, we get
\be \label {105}
\sup_{U^m} \Xi_\delta(U^m \xleftarrow[U \setminus V^m]{} A^m \leftarrow U)
\leq C \omega_m(\delta), \quad \delta \approx 2.
\ee

Plugging (\ref {121}) and (\ref {105})
in (\ref {111}) we get
\be
\frac {E(\Xi_\delta(A^n \leftarrow U)|A^n)}
{E(\Xi_\delta(A^n \xleftarrow[U \setminus A^m]{} U^m)|A^n)} \leq
1+C \frac {\eta_m} {\rho_m}+C \omega_m(\delta), \quad \delta \approx 2.
\ee
This and (\ref {191}) imply (\ref {302}).

\subsubsection{Proof of (\ref {303})}

It is easy to see that if $x \in A^n$ and $f^k(x) \in A^m$ then $|Df^k(x)|
\geq C^{-1} \theta^{m-n}$ for some $\theta<1$.  This follows from the
Schwarz Lemma and the Koebe Distortion Lemma, see Lemma~\ref {exp} for an
analogous estimate.

For simplicity, we will assume that if $x \in A^m$ and $f^k(x) \in A^m$ then
$|Df^k(x)| \geq 1$ (a slightly different argument allows to deal with the
general case\footnote{
Let $W$ be the complement of the postcritical set of $f$ in some definite
neighborhood of $V$ and let $\| \cdot \|_W$ denote the Poincar\'e metric on
$W$.  By the Schwarz Lemma, $f:U \to V$ expands $\| \cdot \|_W$.
Let $\| \cdot \|_*$ denote the metric on $\C$ which is equal to the usual
Euclidean metric outside $\cup A^r$ and in each $A^r$ is given
by $\| \cdot \|$, suitably rescaled so that the $\|\cdot\|_*$ area of $A^r$
coincides with the usual Euclidean area of $A^r$.  Since  the $A^r$ are away
from the postcritical set, the metric $\|\cdot\|_*$ is
comparable to the usual Euclidean metric everywhere.  Moreover, we get
improved inequality that if $x \in A^m$ and $f^k(x) \in A^m$ then
$\|Df^k(x)\|_* \geq 1$ (instead of $|Df^k(x)|>C^{-1}$).
If one defines Poincar\'e series, $\omega_{m,n}$, $\eta_{m,n}$,
$\xi_{m,n}$,... with respect to the metric $\|\cdot\|_*$, the resulting
objects are comparable to the Euclidean ones.  Thus, it is enough to
prove (\ref {303}) for those modified objects, which can be done with the
argument given below.}).
This assumption can be made without loss of generality,
see Remark~\ref {small domains}.

\comm{
We will first prove (\ref {303}) assuming additionally an extra condition
\begin{enumerate}
\item[*] If $x \in A^m$ and $f^k(x) \in A^m$ then $|Df^k(x)| \geq 1$.
\end{enumerate}
}

As in the proof of (\ref {301}), we shall use the identity (\ref {81}).

Notice that the argument for the upper bound on $E(\Xi_2(A^m
\xleftarrow[B]{} B)|A^m)$, (\ref {401}) also implies a lower bound of the
same order,
\be
E(\Xi_2(A^m \xleftarrow[B]{} B)|A^m) \geq C^{-1} \frac {\eta_m} {\rho_m}.
\ee
Similarly, the argument of (\ref {402}) gives the bound
\be
E(\Xi_2(A^m \xleftarrow[B]{}
A^n)|A^m) \geq C^{-1} \xi_{m,n} \rho_{m,n}.
\ee

Thus
\be \label {701}
E(\Xi_\delta(A^m \xleftarrow[B]{} B)|A^m) \geq C^{-1}
\frac {\eta_m} {\rho_m} \tau_m^{2-\delta},
\ee
\be \label {702}
E(\Xi_\delta(A^m \xleftarrow[B]{}
A^n)|A^m) \geq C^{-1} \xi_{m,n} \rho_{m,n} \upsilon_{m,n}^{2-\delta}.
\ee

Recall equality (\ref {92}).  The expansion of the Euclidean metric
along orbits of $A^m \xleftarrow[B]{+} A^m$ implies that
\be
\frac {E(\Xi_\delta(A^m \xleftarrow[B]{+} A^m
\leftarrow U)|A^m)} {E(\Xi_\delta(A^m \leftarrow U)|A^m)} \geq
\frac {E(\Xi_\delta(A^m \setminus Z_{m,n} \leftarrow U)|A^m)}
{E(\Xi_\delta(A^m \leftarrow U)|A^m)}=
1-\frac {E(\Xi_\delta(Z_{m,n} \leftarrow U)|A^m)}
{E(\Xi_\delta(A^m \leftarrow
U)|A^m)},
\ee
for $\delta \leq 2$.  By (\ref {93}), we get
\be \label {704}
E(\Xi_\delta(A^m \xleftarrow[B]{+} A^m
\leftarrow U)|A^m) \geq (1-C \kappa_{m,n}) E(\Xi_\delta(A^m \leftarrow
U)|A^m)=(1-C \kappa_{m,n}) \omega_m(\delta).
\ee

Plugging (\ref {701}), (\ref {702}) and (\ref {704}) in (\ref {81}), we
obtain (\ref {303}).

\subsubsection{Proof of (\ref {304})}

We will use the identity (\ref {115}), which leads to the lower bound
(analogous to (\ref {111}))
\be \label {2111}
\frac {E(\Xi_\delta(A^n \leftarrow U)|A^n)}
{E(\Xi_\delta(A^n \xleftarrow[U \setminus A^m]{} U^m)|A^n)} \geq
1+\inf_{U^m}\,
\Xi_\delta(U^m \xleftarrow[U \setminus V^m]{} U \setminus V^m)+
\inf_{U^m}\, \Xi_\delta(U^m \xleftarrow[U \setminus V^m]{} A^m
\leftarrow U).
\ee

By (\ref {801}) and (\ref {101}),
\be \label {2801}
\inf_{U^m} \Xi_\delta(U^m \xleftarrow[U \setminus V^m]{} U \setminus V^m)
\geq C^{-1} \tau_m^{2-\delta}
\frac {\eta_m} {\rho_m},
\ee
\be \label {2101}
\inf_{U^m} \Xi_\delta(U^m \xleftarrow[U \setminus V^m]{} A^m) \geq C^{-1}.
\ee
By bounded oscillation of Poincar\'e series in $A^m$ and (\ref {2101}),
we have
\be \label {2302}
\inf_{U^m} \Xi_\delta(U^m \xleftarrow[U \setminus V^m]{} A^m
\leftarrow U) \geq C^{-1} \omega_m(\delta).
\ee

Plugging (\ref {2801}) and (\ref {2302}) in (\ref {2111}), and using
(\ref {191}) one gets (\ref {304}).

\section{Trichotomy}\label{trichotomy}

Let $f$ be a periodic point of renormalization of period $p$.  In this case,
we may assume that the fundamental annulus $A^n$ form a periodic sequence
up to rescaling by $\rho_n^{-1/2}$ (see Remark~\ref {permarkimpro}).  Then
$\eta_{m,n}$, $\xi_{m,n}$,... only depend on $n-m$ and the congruence
class of $m$ modulo $p$.

The results in this section will be obtained from some particular cases of
the estimates in Lemmas~\ref {pert1} and~\ref {pert2}.
Setting $n=2m$ an integer multiple of $p$,
the estimates can be rewritten as
\be \label {2n11}
\omega_m^{[j+1]}(\delta) \leq C \frac {\eta_m}
{\rho_m}+(1-C^{-1}(\eta_m+\xi_m)\omega_m^{[j]}(\delta)+C \xi_m \rho_m
\omega_{2m}^{[j]}(\delta), \quad \delta \approx 2,
\ee
\be \label {2n1}
\omega_{2m}^{[j+1]}(\delta) \leq C \left (\frac {\eta_m}
{\rho_m}+\omega_m^{[j]}(\delta) \right )\omega_m^{[j]}(\delta),
\quad \delta \approx 2,
\ee
\be \label {2n12}
\omega_m(\delta) \geq C^{-1} \tau_m^{2-\delta}
\frac {\eta_m} {\rho_m}+\max \{1-C (\eta_m+\xi_m),0\}
\omega_m(\delta)+C^{-1} \upsilon^{2-\delta}_m \xi_m \rho_m
\omega_{2m}(\delta),
\ee
\be \label {2n2}
\omega_{2m}(\delta) \geq C^{-1} \left (\tau^{2-\delta}_m
\frac {\eta_m} {\rho_m}+\omega_m(\delta) \right )\omega_m(\delta).
\ee

Putting the respective estimates  together, we get the following
quadratic estimates: 
\be \label {above}
\omega_m^{[j+1]}(\delta) \leq
C \frac {\eta_m} {\rho_m}+(1-C^{-1}(\eta_m+\xi_m)+C \eta_m \xi_m)
\omega_m^{[j]}(\delta)+
C \xi_m \rho_m \omega_m^{[j]}(\delta)^2, \quad \delta \approx 2,
\ee
\be \label {below}
\omega_m(\delta) \geq
C^{-1} \tau_m^{2-\delta} \frac {\eta_m} {\rho_m}+(\max
\{1-C(\eta_m+\xi_m),0\}+
C^{-1} (\tau_m \upsilon_m)^{2-\delta} \eta_m \xi_m) \omega_m(\delta)+
C \upsilon_m^{2-\delta} \xi_m \rho_m \omega_m(\delta)^2.
\ee

Let us now explore those relations.

\begin{lemma}\label{bounds on om}

Let $f$ be a periodic point of renormalization of period $p$.  If $m=0 \mod
p$ then there exists $C>0$ (only depending on the geometric
bounds) such that
\begin{enumerate}
\item If the polynomial
\be \label {up}
   P:  x \mapsto C \frac {\eta_m} {\rho_m}+(1-C^{-1}(\eta_m+\xi_m)+C \eta_m
\xi_m)x+C \xi_m \rho_m x^2
\ee
has a real positive fixed point then $\delta_\crit<2$ and for $\delta$
close to $2$, the smallest such fixed point gives an upper bound on
$\omega_m(\delta)$.
\item If $\delta \geq \delta_\crit$, the polynomial
\be \label {down}
     Q_\de : x \mapsto C^{-1} \tau_m^{2-\delta} \frac {\eta_m} {\rho_m}+
(\max \{1-C(\eta_m+\xi_m),0\}+
C^{-1} (\tau_m \upsilon_m)^{2-\delta} \eta_m \xi_m) x+
C \upsilon_m^{2-\delta} \xi_m \rho_m x^2
\ee
has a fixed point and the smallest (respectively, largest)
such fixed point gives a lower bound (respectively, upper bound)
to $\omega_m(\delta)$.
\end{enumerate}

\end{lemma}

\begin{pf}

By repeated application of (\ref {above}), for $\delta \approx 2$ and for
every $j \geq 0$, we have
$\omega^{[j]}_m(\delta) \leq P^{j+1}(0)$.  If there exists a positive
fixed point for $P$ then the smallest such fixed point
attracts the orbit of $0$.  This gives the first statement.

By (\ref {below}), for $\delta>\delta_\crit$, $Q_\delta(\omega_m(\delta))
\leq \omega_m(\delta)$.  Thus $\omega_m(\delta)$ is contained in the
interval bounded by the fixed points of $Q_\delta$.  This still holds for
$\delta=\delta_\crit$ by continuity of the coefficients of $Q_\delta$.  This
implies the second statement.
\end{pf}

\begin{thm}[Balanced case]\label{balanced case}

Let $f$ be a periodic point of renormalization of period $p$.
Then there exists $C>0$,
only depending on the geometric bounds, such that if $\delta_\crit=2$ then
\be \label {eta xi}
C^{-1} \leq \frac {\eta_m} {\xi_m} \leq C, \quad m\equiv 0\ \mod p.
\ee
Moreover,
\be\label{1/n}
      \xi_m \asymp \eta_m \asymp \frac {1} {m}
\ee
for all $m\in \N$ with the constants depending also on $p$.

\end{thm}

\begin{pf}

If (\ref {eta xi}) does not hold then 
\be\label{means}
  \xi_m \eta_m << (\xi_m +\eta_m)^2,
\ee
and then the quadratic polynomial  (\ref {up})  has two positive fixed points.
By the first part of Lemma \ref{bounds on om},  $\delta_\crit<2$ contradicting the assumption.  
The same argument applied to the $k$-fold renormalization of $f$, $k=0,1,\dots, p-1$,  gives
$$
\eta_m \asymp \xi_m, \quad m\equiv k\ \mod p,
$$
with the constant depending on $k$. 
Thus, $\eta_m \asymp \xi_m$ for all $m$, with the constant depending on $p$.
By Lemma~\ref {exp}, this implies that
$$
\eta_m-C \eta_m^2 \leq \eta_{m+1} \leq \eta_m-C^{-1} \eta_m^2.
$$
Hence
$$
C^{-1} \leq \frac{1}{\eta_{m+1}}-\frac{1}{\eta_m}  \leq C, 
$$
and thus $1/ \eta_m \asymp m$. 
\end{pf}

The following result gives Theorem~D stated in the Introduction.

\begin{cor}\label{area and de}

Let $f$ be a periodic point of renormalization.  If $|J(f)|>0$ then
$\delta_\crit<2$.

\end{cor}

\begin{pf}
Since $f$ is a renormalization periodic point, $f_m: U^m \ra V^m$ is a rescaling of $f: U\ra V$ for $m\equiv 0\ \mod p$.
Hence  
$$
    \eta_m = p(X_m|U^m)  \geq p(J(f_m) | U^m) = p( J(f)| U)  >0, \quad m\equiv 0 \ \mod p . 
$$ 
Then by Lemma \ref{exp}, $\xi_m\to 0$, so that (\ref{eta xi}) is violated -
this contradicts $\delta_\crit=2$.
\end{pf}

\begin{lemma} \label {lemmabla}

Let $f$ be a periodic point of renormalization of period $p$ and let $m\equiv 0 \ \mod p$.  
Then there exists $C>0$ depending only on the geometric bounds such
that
\be \label {n2}
\omega_m(2) \geq C^{-1} \frac {\eta_m} {(\xi_m+\eta_m) \rho_m} \quad
\text {if }\quad  \eta_m\xi_m<C^{-1},%
\footnote{Note that this assumption is always satisfied on high levels  since  $\eta_m\xi_m\to 0$ by Lemma \ref{exp}}
%\leq \omega_n(2) \leq C \frac {\eta_n} {(\xi_n+\eta_n) \rho_n}.
\ee
\be \label {n7}
\omega_m(2) \leq C \frac {\eta_m} {(\xi_m+\eta_m) \rho_m} \quad 
  \text {if }\quad  \left| \ln \frac{\eta_m}{\xi_m} \right|>C,
\ee
\be \label {8n1}
\omega_{2m}(2) \asymp \left (\frac {\eta_m} {\rho_m}+\omega_m(2) \right
)\omega_m(2) \asymp \omega_m(2)^2.
\ee
In particular, if $\eta_m \xi_m<C^{-1}$ and if either
$\eta_m/\xi_m$ or $\omega_m \rho_m$ is small,  then both are small and indeed
they are comparable.

\end{lemma}

\begin{pf}

If $\eta_m \xi_m$ is small then the linear coefficient of $Q_2$ is comparable with $\eta_m + \xi_m$, 
and a direct calculation of the smallest fixed point of $Q_2$ 
gives the lower bound (\ref {n2}) on $\omega_m(2)$. 

 If $|\ln (\eta_m/ \xi_m) |$ is large, then (\ref{means}) holds,
and the polynomial $P$ has two real fixed points. 
A direct calculation  of  the smallest one gives the upper bound (\ref {n7}) on $\omega_m(2)$.

By (\ref {2n1}), (\ref {2n2}),
\be \label {2n47}
\omega_{2m}(2) \asymp \left (\frac {\eta_m} {\rho_m}+\omega_m(2) \right )
\omega_m(2).
\ee
By (\ref {2n12}), $\omega_m(2)>C^{-1} \eta_m/\rho_m$ which
shows that $(\eta_m/\rho_m)+\omega_m(2) \asymp \omega_m(2)$ which together
with (\ref {2n47}) implies (\ref {8n1}).

The last assertion is a straightforward consequence of (\ref {n2}) and (\ref {n7}).
\end{pf}

\comm{
If $\eta_m\xi_m$ is small then it is possible to estimate the smallest
fixed point of (\ref {down}) when $\delta=2$.  We have
\be \label {n2}
\omega_m(2) \geq C^{-1} \frac {\eta_m} {(\xi_m+\eta_m) \rho_m}, \quad
\eta_m\xi_m<C^{-1}.
%\leq \omega_n(2) \leq C \frac {\eta_n} {(\xi_n+\eta_n) \rho_n}.
\ee

If moreover $|\ln \eta_m-\ln \xi_m|$ is large, then the smallest fixed point
of (\ref {up}) exists and is comparable to the smallest fixed point of (\ref
{down}), so we can write
\be \label {n7}
\omega_m(2) \asymp \frac {\eta_m} {(\xi_m+\eta_m) \rho_m}, \quad |\ln
\eta_m-\ln \xi_m|>C.
\ee

It follows that
\be \label {n8}
\omega_m(2) \asymp \frac {\eta_m} {\xi_m \rho_m}, \quad \eta_m<C^{-1} \xi_m.
\ee
%By (\ref {2n1}), (\ref {2n2}),
%\be \label {3n1}
%\omega_{2n}(2) \approx \left (\frac {\eta_n} {\rho_n}+\omega_n(2) \right
%)\omega_n(2),
%\ee

By (\ref {2n1}), (\ref {2n2}),
\be \label {8n1}
\omega_{2m}(2) \asymp \left (\frac {\eta_m} {\rho_m}+\omega_m(2) \right
)\omega_m(2) \asymp \omega_m(2)^2,
\ee
where the last estimate comes from $\omega_m(2)>C^{-1} \eta_m/\rho_m$, which
follows from (\ref {2n12}).
}

\begin{thm}[Lean case] \label {Shallow case}

Let $f$ be a periodic point of renormalization of period $p$.
Then there exists $C>0$, only depending
on the geometric bounds, such that if for some $m\equiv 0 \ \mod p$ we have
\be \label {433}
\eta_m<C^{-1} \xi_m
\ee
then $\inf \xi_n > 0$.\footnote {An explicit  lower bound for $\inf \xi_n$ in terms of
the geometric and combinatorial bounds and the period of $f$ can be worked out.}  
Moreover, $\eta_n \to 0$ exponentially fast.

\end{thm}

\begin{pf}
Since $\eta_n$ monotonically decrease, 
assumption (\ref {433}) implies that $\eta_n$ (and hence  $\eta_n \xi_n$) are small for $n \geq m$.
It also implies by Lemma~\ref {lemmabla},  that $\omega_m \rho_m$ is small. 

By (\ref{8n1}), $\om_{2m} \leq K \om_m^2$, so that
$  K \om_{2m} \rho_{2m} \leq  (K \om_m \rho_m)^2.$
Iterating this estimate $k$ times, we conclude that 
$$
   K \omega_{2^k m} \rho_{2^k m} \leq ( K \omega_m \rho_m )^{2^k}\equiv q^{2^k} 
$$
where $q<1$ provided $\omega_m \rho_m $ is small enough.  
By Lemma~\ref {lemmabla},
$$
   \frac{\eta_{2^k m}}{\xi_{2^k m}} = O(q^{2^k}).
$$ 
Since $\eta_n$ monotonically decrease, 
we conclude that  they decay exponentially:  $\eta_n  = O (q^{n/2})$. 

Furthermore, it follows from the renormalization periodicity that the $\xi_{kp}$ monotonically decrease in $k$.
It is also easy to check that $\xi_{n+1}\asymp \xi_n$.
%\note{OK?}
It follows that  
$$
  \inf \xi_n=0 \implies \xi_n\to 0. 
$$
Since $\eta_n$ decay exponentially, $\lim\sup \xi_n>0$ by Lemma~\ref {exp}. 
Hence $\xi_n$ are bounded away from $0$ .
\comm{
By (\ref {n2})$\eta_m/\xi_m$ is
small implies thus by induction that
The assumption implies that $\eta_r$ (and hence $\eta_r \xi_r$) is
small for $r \geq m$, and allows us to use estimate (\ref {n2}).
By (\ref {2n1}), (\ref {2n2}),
\be \label {8n1}
\omega_{2m}(2) \asymp \left (\frac {\eta_m} {\rho_m}+\omega_m(2) \right
)\omega_m(2) \asymp \omega_m(2)^2
\ee
(the last estimate follows from (\ref {n2})).
By (\ref {n8}), we have $\omega_m(2) \rho_m \leq
C^{-1}$.
Applying (\ref {8n1})$ we get $\omega_{2m}(2)
\rho_{2m}=\omega_{2m}(2) \rho^2_m \asymp (\omega_m(2) \rho_m)^2 \leq
C^{-2}$.  By (\ref {n2})
By induction we get $\omega_{2^k m} \rho_{2^k m} \leq C^{-2^k}$.
Thus $\eta_r \to 0$ exponentially fast.  By Lemma \ref {exp}, this implies
that $\xi_r$ is bounded away from $0$.
}
\end{pf}

\comm{
If $|\ln \eta_m-\ln \xi_m|$ is large, the smallest fixed points of (\ref
{above}) and (\ref {below}) are given (up to multiplicative constant) by
$\eta_n/(\eta_n+\xi_n)\rho_n$.  Thus if $\eta_n/\xi_n$ is small then
\be
\omega_n(2) \approx \frac {\eta_n} {\xi_n \rho_n}.
\ee
By (\ref {2n1}), (\ref {2n2}),
\be
\omega_{2n}(2) \approx \left (\frac {\eta_n} {\rho_n}+\omega_n(2) \right
)\omega_n(2),
\ee
}

\begin{lemma} \label {measurezero}

Let $f$ be a periodic point of renormalization of period $p$ and let
$m\equiv 0\ \mod p$.  If $\area (J(f))=0$ then there exists $C>0$, only
depending on the geometric bounds, such that
\be \label {977}
C^{-1} \frac {\eta_m} {\xi_m \rho_m} \leq \omega_m(2) \leq
C \frac {\eta_m} {\xi_m \rho_m}.
\ee

\end{lemma}

\begin{pf}
Since the set $Y_m$ is wandering,
$$
   \int_{Y_m} \Xi_2(f,x) dx  = \area \left(\bigcup_{n=0}^\infty f^{-n}(Y_m)\right) \equiv |Z_m|.
$$
The set  $Z_m$ consists of points that land at $A^m$ at least once but at most finitely many times.
Hence 
$$ 
   X_m\sm J(f)\subset Z_m \subset X_m. 
$$
Since area of the Julia set is zero, $| Z_m| = | X_m| $ and hence 
\be
\int_{Y_m} \Xi_2(f,x) dx =| X_m|.
\ee
By bounded oscillation of the Poincar\'e series, 
\be
\omega_m(2) \asymp  \frac{1}{|Y_m|} \int_{Y_m} \Xi_2(f,x) dx = \frac {|X_m|} {|Y_m|} \asymp \frac {\eta_m} {\xi_m \rho_m}.
\ee
\end{pf}

\begin{lemma} \label {10.9}

Let $f$ be a periodic point of renormalization of period $p$.  If $\delta_\crit<2$ then  there
exists $\theta<1$, only depending on $\delta_\crit$ and the geometric bounds
such that
\be\label{one more}
\frac {\eta_m \xi_m} {(\eta_m+\xi_m)^2} \leq \theta^m, \quad m\equiv 0 \ \mod p.
\ee

\end{lemma}

\begin{pf}
By Lemma \ref{bounds on om},
the quadratic polynomial $Q_{\de_\crit}$ has a real  fixed point and hence has the positive
discriminant. Together with  Lemma~\ref {tauupsilon} this implies (\ref{one more}).
\end{pf}

\begin{thm}[Black hole case] \label {Black hole case}

Let $f$ be a periodic point of renormalization of period $p$.  Then there exists $C>0$,
only depending on the geometric bounds, such that if for some $m\equiv 0 \ \mod p$
we have
\be \label {937}
\eta_m>C \xi_m
\ee
then $\lim \eta_n >0$.\footnote {A lower bound for $\lim \eta_n $ in terms of
the geometric and combinatorial bounds and the period of $f$ can be worked out.} 
 Moreover, $\xi_n \to 0$ exponentially
fast.

\end{thm}

\begin{pf}

%Using the upper estimate, we get that if $x \in A^n$ then
%$\Xi_2(f,x)$ is bounded by the smallest fixed point of
%\be
%P(x)=C \frac {\eta_n} {\rho_n}+(1-C^{-1}(\xi_n+\eta_n)+C \xi_n \eta_n)x+C
%\xi_n \rho_n x^2.
%\ee
%Thus
%\be
%\Xi_2(f,x) \leq C \frac {\eta_n} {(\eta_n+\xi_n) \rho_n}.
%\ee
%so that there exists $x \in Y_{n_0,n}$ such that
%\be
%\Xi_2(f,x) \geq \frac {\eta_n |U|} {|Y_n|} \geq C \frac {\eta_n} {\xi_n
%\rho_n}
%\ee
Assumption (\ref {937}) allows us to apply (\ref {n7}).
Since (\ref {937}) and  (\ref {n7}) are not compatible with (\ref {977}),
area of the Julia set must be positive and hence $\lim \eta_n>0$.\footnote
{Notice that this conclusion gives  a {\it criterion for positive
measure} which might be checkable numerically.
Let us remark that this conclusion is based on an {\it upper bound} on
$\omega_m(2)$, which is derived from Lemma~\ref {pert1}.} 

Moreover, by Corollary \ref{area and de}, $\delta_\crit<2$ which implies
by Lemma \ref {10.9} that $\eta_n \xi_n \to 0$ exponentially fast, so
$\xi_n \to 0$ exponentially fast as well.
\end{pf}

\comm{
\begin{thm}

Let $f$ be a periodic point of renormalization.  Then there exists $C>0$,
only depending on the geometric bounds, such that if for some $n$ we have
\be
\eta_n<C^{-1} \xi_n
\ee
then $\liminf \xi_n>0$.\footnote {A lower bound for $\lim \xi_n$ involving
the geometric and combinatorial bounds and the
period of $f$ can be worked out.}  Moreover, $\eta_n \to 0$ exponentially
fast.

\end{thm}

\begin{pf}

It is easy to see that if $n=0 \mod m$ then
\be
1+\frac {\xi_{2n}} {\eta_{2n}} \geq C \left (1+\frac {\xi_n} {\eta_n}
\right )^2.
\ee
\end{pf}

\begin{lemma}

Let $f$ be a periodic point of renormalization of period $m$.
There exists $C>0$ (only depending on the geometric bounds)
with the following property.  If
\be
|\ln \eta_n-\ln \kappa_n|>C
\ee
for some $n=0 \mod m$ then $\delta_\crit<2$.

\end{lemma}

\begin{cor}

Let $f$ be a periodic point of renormalization
whose critical exponent is $2$.  Then
\be
|\ln \eta_n-\ln \kappa_n|
\ee
and
\be
|\ln n+\ln \eta_n|
\ee
are bounded.\footnote{One can work out explicit bounds (which depend on the
geometric bounds and on the period of $f$).}

\end{cor}

\begin{pf}

The first estimate follows immediately from the previous lemma.
From
\be
(1-C \eta_n) \eta_n<\eta_{n+1}<(1-C^{-1} \eta_n) \eta_n
\ee
we get
\be
C^{-1}<\eta_{n+1}^{-1}-\eta_n^{-1}<C,
\ee
which implies the second estimate.
\end{pf}

\begin{cor}

Let $f$ be a periodic point of renormalization whose Julia set has positive
Lebesgue measure.  Then $\delta_\crit<2$.

\end{cor}

\begin{pf}

Indeed in this case $\eta_n \geq |J(f)|$ is bounded away from zero.
\end{pf}

\begin{lemma}

Let $f$ be a periodic point of renormalization of period $m$.  There exists
$C>0$, $\theta<1$ (only depending on the geometric bounds) such that if $n=0
\mod m$ then
\be
\frac{\kappa_n \eta_n} {(\kappa_n+\eta_n)^2} \leq C
\theta^{(2-\delta_\crit)n}.
\ee

\end{lemma}

\begin{cor}

Let $f$ be a periodic point of renormalization such that $\delta_\crit<2$.
Then one of $\eta_n$, $\kappa_n$ is bounded from below and the other decays
exponentially fast.

\end{cor}

Given a holomorphic map $f:U \to V$, we let $\LL_f$ (respectively,
$\LL_f^+$) denote the set of all $(x,y,m)$ such that
$m \geq 0$ (respectively, $m>0$), and $f^m(y)=x$.

Given to subsets $\GG, \HH \subset \LL_f$, we let $\GG \circ \HH$ (the
composition of $\GG$ and $\HH$) be the set
of all $(x,z,m) \in \LL_f$ such that there exists $0 \leq k \leq m$ and $y$
such that $(x,y,m-k) \in \GG$ and $(y,z,k) \in \HH$.

Given sets $X$, $Y$, $Z$, we let
$\LL_f(X,Y,Z)$ (respectively, $\LL_f^+(X,Y,Z)$
denote the set of $(x,y,m)$ in $\LL_f$ (respectively,
$\LL_f^+$) where $x \in X$, $y \in Y$ and $f^k(y) \in Z$ for $0<k<m$.

\newcommand{\gra}[2]{\xleftarrow[{#2}]{{#1}}}
\newcommand{\grap}[2]{\xleftarrow[{#2}]{{#1,+}}}

We shall use the following more suggestive notation with arrows:
\be
X \gra{f}{Z} Y, \quad \text {for } \LL_f(X,Y,Z),
\ee
\be
X \grap{f}{Z} Y, \quad \text {for } \LL_f^+(X,Y,Z).
\ee
We shall also denote composition by concatenation of arrows, for instance:
\be
X \gra{f}{Z} Y \gra{f}{T} S, \quad \text {for the composition of }
X \gra{f}{Z} Y \text { and } Y \gra{f}{T} S.
\ee

If $X$ has positive Lebesgue measure, we let
\be
\Xi^{ave}_\delta(X \gra{f}{Z} Y)=\int_X \sum_{(x,y,z) \in \LL_f(X,Y)}

\begin{thm}

If $f$ is a periodic point of renormalization and its Julia set has positive
Lebesgue measure then its critical exponent is less than $2$.

\end{thm}

It will soon be clear that it is preferrable to state the previous theorem
as a tricotomy as follows:

\begin{thm}

Let $f:U \to V$ be a periodic point of renormalization.  Then either:
\begin{enumerate}
\item The Lebesgue measure of $J(f)$ is positive.  In this case its critical
exponent is less than $2$.  This happens if and only if $\eta_n$ is bounded
away from zero and $\kappa_n \to 0$ exponentially fast.
\item The Lebesgue measure of $J(f)$ is zero and the critical exponent is
$2$.  This happens if and only if $\eta_n/\kappa_n$ is bounded away from
zero and infinity, and we have $\lim \ln \eta_n/\ln n=\lim \ln \kappa_n/\ln
n=1$.
\item The Lebesgue measure of $J(f)$ is zero and the critical exponent is
less than $2$.  In this case $\kappa_n$ is bounded away from zero and
$\eta_n \to 0$ exponentially fast.
\end{enumerate}

\end{thm}

The three possibilities in the tricotomy can be distinguished qualitatively
in terms of the sequences $\eta_n$, $\kappa_n$ as follows:

\begin{thm}

In the setting of the previous theorem we have
\begin{enumerate}
\item Holds if and only if $\eta_n$ is bounded away from zero (by $|J(f)|$)
and $\kappa_n$ converges to zero exponentially fast.
\item Holds if and only if the sequence $|\ln \eta_n-\ln \kappa_n|$ is
bounded and in this case $\ln \eta_n \asymp \ln n$.
\item Holds if and only if $\kappa_n$ is bounded away from zero
and $\eta_n$ converges to zero exponentially fast.
\end{enumerate}
\end{thm}

A more precise ``algorithmic'' statement is as follows:

\begin{thm}

In the same setting of the previous theorem, there exists constants $C$
(only depending on the degree and combinatorial bounds of $f$) and $c$
(depending also on the critical exponent), such that
\begin{enumerate}
\item $\ln \eta_{mn}-\ln \kappa_{mn}>C$, $\eta_n \geq |J(f)|
\geq c^{-1}$, $\kappa_n \leq e^{-cn}$.
\item $|\ln \eta_{mn}-\ln \kappa_{mn}|<C$.
\item $\ln \kappa_{mn}-\ln \eta_{mn}>C$, $\kappa_n \geq c^{-1}$, $\eta_n
\leq e^{-cn}$.
\end{enumerate}

\end{thm}

We call the above statement algorithmic because the constant $C$ is
effective at least in the real case (it depends on the complex bounds),
hence if (2) does not hold then this can be in principle decided in
finite time from the computation of the sequences
$\eta_{mn}$, $\kappa_{mn}$ (and in this case one also
decides whether (1) or (3) holds).

Let $A \ntop {\underleftarrow {\scriptscriptstyle {+}}}
{\overleftarrow {\scriptscriptstyle {C}}} B$.
Let $A \scriptscriptstyle {\ntop {\underleftarrow {\scriptscriptstyle {+}}}
{\ntop {\not \leftarrow} {\scriptscriptstyle {C}}}} B$.
Let $\ntop {\ntop {\scriptscriptstyle +}
{\scriptscriptstyle {\overleftarrow{\not
\longleftarrow}}}} {\scriptscriptstyle C}$.

\comm{
\section{Combinatorial breakdown}

Let us consider an infinitely renormalizable map $f:U \to V$ with
prerenormalization $f_n:U^n \to V^n$ and let $A^n=V^n \setminus U^n$.  We
will assume that $U^n,V^n,A^n$ have bounded shape, and that $A^n$ is away
from the postcritical set (that is, it has bounded hyperbolic diameter).
We will also assume for simplicity that $V^{n+1} \subset U^n$ and that
$f_{n+1}$ is the first return map of $U^{n+1}$ to $V^{n+1}$ under $f_n$.

A $4$-tuples of the type $(f_n,x,m,y)$ will be called fine if $f_n^m(y)=x$,
$m \geq 0$.
A $4$-tuple is said to start in $y$ and end in $x$ via $f_n$.  A fine
$4$-tuple is said to pass through points $f(x),...,f^{m-1}(x)$.
We shall denote by caligraphic letters, as $\FF$, to denote sets of fine
$4$-tuples.  We shall denote by
the corresponding roman letter, as $F$, the image of the projection on the
first coordinate.  We shall call $\FF$ injective if the projection on the
fourth coordinate is injective.

Given sets $B,C \subset U^n$, $A \subset V^n$ we will denote by
$\LL(B,C,D,f_n)$ the set of all fine $4$-tuples starting at $B$ and ending
in $C$ via $f_n$, but which do not pass through $D$.  By $\LL^+(B,C,D,f_n)
\subset \LL(B,C,D,f_n)$ we shall denote the subset corrsponding to $m>0$.

Let us denote by $\FF^n_0=\LL(U,V,\emptyset,f)$,
$\FF^n_1=\LL^+(A^n,A^n,A^n \cup A^{2n},f)$,
$\FF^n_2=\LL(U,A^n,A^n \cup A^{2n},f)$,
$\FF^n_3=\LL(A^{2n},A^n,A^{2n},f)$,
$\FF^n_4=\LL(U^n,A^{2n},A^n,f)$,
$\FF^n_5=\LL(U,U^n,A^n,f)$,
$\FF^n_7=\LL(A^n,U^n,A^n,f)$,
$\FF^n_5=\LL(A^n,A,A^n,f)$,
$\FF^n_7=\LL(U^n,A^{2n},A^n,f)$, $\FF^n_8=\LL(A^n,U^n,U^n,f)$.
All of them are injective except $\FF^n_5$.

Let $\eta_n$ be the probability of a point in $U$ to have a non-negative
iterate in $V^n$.  Let $\delta_n$ be the probability of
a point in $A^n$ to have no positive iterate in $V^n$.  Let $\xi_n$ be the
probability of a point in $A^n$ to not return to $A^n$ before passing by
$V^{2n}$.  We let $\rho_n$ be the probability of $V^n$ in $V$.

Let $X \subset \C$ be a measurable set of positive Lebesgue measure.  We let
\be
\Xi^\sup_\delta(X,\FF)=\esssup_{x \in X} \sum_{(f_n,x,m,y) \in \FF}
|Df_n^m(y)|^{-\delta},
\ee
\be
\Xi^\inf_\delta(X,\FF)=\essinf_{x \in X} \sum_{(f_n,x,m,y) \in \FF}
|Df_n^m(y)|^{-\delta},
\ee
\be
\Xi^\ave_\delta(X,\FF)=\frac {1} {|X|} \int_X \sum_{(f_n,x,m,y) \in \FF}
|Df_n^m(y)|^{-\delta} dx.
\ee

We first observe that
\be
\Xi^\ave_2(A^n,\FF^n_1) \leq \Xi^\ave_2(\FF^n_1) \Xi^\sup(A^n,\FF^n_0)$.
\ee

Clearly
\be
\Xi^\ave_2(A^n,\FF^n_1) \leq (1-C^{-1}\xi_n) \Xi^\ave_2(A^n,\FF^n_0).
\ee

Moreover,
\be
\Xi^\sup_2(A^n,\FF^n_2) \leq C \frac {\eta_n} {\rho_n}.
\ee
and
\be
\Xi^\ave_2(A^n,\FF^n_3) \leq C \delta_n.
\ee

Thus we have
\be
\Xi^\ave_2(A^n,\FF^n_0) \leq (1-C^{-1} \xi_n) \Xi^\ave_2(A^n\FF^n_0)+
C \frac {\eta_n} {\rho_n}+\Xi^\ave_2(A^n,\FF^n_3)
\Xi^\sup_2(A^{2n},\FF^n_0),
\ee
which gives
\be
\Xi^\ave_2(A^n,\FF^n_0) \leq (1-C^{-1} \xi_n) \Xi^\ave_2(A^n\FF^n_0)+
C \frac {\eta_n} {\rho_n}+C \delta_n \rho_n
\Xi^\sup_2(A^{2n},\FF^n_0).
\ee

We now estimate
\be
\Xi^\sup_2(A^{2n},\FF^n_0) \leq
\Xi^\sup_2(A^{2n},\FF^n_4)(\Xi^\sup_2(U^n,\FF^n_5)+\Xi^\sup_2(U^n,\FF^n_7)
\Xi^\sup_2(A^n,\FF^n_0)).
\ee

Notice that
\be
\Xi^\sup_2(A^{2n},\FF^n_4)=\Xi^\sup_2(A^n,\FF^n_0),
\ee
\be
\Xi^\sup_2(U^n,\FF^n_5) \leq C \frac {\eta_n} {\rho_n},
\ee
\be
\Xi^\sup_2(U^n,\FF^n_7) \leq C.
\ee

We get
\be
\Xi^\sup_2(A^{2n},\FF^n_0) \leq \Xi^\sup_2(A^n,\FF^n_0)(C \frac {\eta_n}
{\rho_n}+C\Xi^\sup_2(A^n,\FF^n_0)).
\ee

Using that
\be
\Xi^\sup_2(A^n,\FF^n_0) \leq C \Xi^\ave_2(A^n,\FF^n_0),
\ee
we get
\be
\Xi^\ave_2(A^n,\FF^n_0) \leq C \frac {\eta_n} {\rho_n}+(1-C^{-1}
\xi_n+C\delta_n\eta_n) \Xi^\ave_2(A^n,\FF^n_0)+C \delta_n \rho_n
(\Xi^\ave_2(A^n,\FF^n_0))^2.
\ee

This estimate is the basis of the following:

\begin{thm}

Let $f$ be a fixed point of renormalization.  There exists $C>0$ (only
depending on the combinatorial bounds of $f$), such that
for every $n$, if the quadratic polynomial
\be
C \delta_n \eta_n+(1-C^{-1} \xi_n+C \delta_n \eta_n)x+x^2
\ee
has a fixed point then the critical exponent is less than $2$.

\end{thm}

\begin{cor}

If the critical exponent is $2$ then $\delta_n/\eta_n$ is bounded away from
$0$ and $\infty$.

\end{cor}

\begin{pf}

Notice that
\be
\frac {\xi_n} {\delta_n+\eta_n}
\ee
is bounded away from zero and infinity.
\end{pf}

\begin{cor}

If the Lebesgue measure is positive then the critical exponent is less than
$2$.

\end{cor}

\begin{pf}

Indeed, in this case $\lim \eta_n>0$ but $\lim \delta_n=0$.
\end{pf}

We will also assume for simplicity that $V^{n+1} \subset U^n$ and that
$R^{n+1}(f)$ is the first return map of $U^{n+1}$ to $V^{n+1}$ under
$R^n(f)$.

When considering a periodic point of renormalization (of period $n$)
we will assume that $U^{i+n}=\rho^{1/2}_n U^i$.

In what follows, $C$ will denote a constant depending only on the bounds on
the combinatorial type.

\begin{lemma}

We have that $\xi_n \leq C(\eta_n+\delta_n)$.

\end{lemma}

\begin{lemma}

If $\lim \eta_n>0$ then $\lim \delta_n=0$.

\end{lemma}

\begin{pf}

We are going to show that $\ln \eta_n \geq -C^{-1}\sum_{m<n} \delta_m$.
Notice that the inequality $\ln \eta_n \leq -C\sum_{m<n} \delta_m$ also
holds but will not be needed.
\end{pf}

\begin{lemma}

Let $f$ be a periodic point of renormalization of period $n$.  Then the
quadratic polynomial
\be
P(x)=C^{-1} \eta+(1-C \xi_n)x+x^2
\ee
has a fixed point.

\end{lemma}

\begin{pf}

It actually holds that
\be
\hat P(x)=C^{-1} \eta+(1-C \xi_n+C^{-1}\eta_n)x+x^2
\ee
has a fixed point, but we will not need it.

Let $X$ be the set of points in $A^n$ that return to $A^n$ without passing by
$V^{2n}$.  We let $\XX$ be the set of $(m,y)$ such that $y \in X$ and
$f^m(y)$ is the first return map to $A^n$.

Let $Y$ be the set of points that land in $A^n$ before passing by $A^{2n}$.
Let $\YY$ be the set of $(m,y)$ such that $y \in Y$ and $f^m(y)$ is the
first landing into $A^n$.

Let $Z$ be the set of points in $A^{2n}$ that land in $A^n$ before returning
to $A^{2n}$ (or $V^{2n}$).  Let $\ZZ$ be
the set of $(m,y)$ such that $y \in Z$ and $f^m(y)$ is the first landing
into $A^n$.

Let $W$ be the set of points in $A^n$ that have a positive iterate in $V^n$. 

Notice that
\be
\inf_{x \in A^n} \sum_{(m,y) \in \YY_x}
|Df^m(y)|^{-2} \geq C^{-1} \frac
{\eta_n} {\rho_n}.
\ee

Moreover,
\be
\int_X \Xi_2(f,x) dx \geq (1-C\xi_n) \int_{A^n} \Xi^2(f,x) dx.
\ee

Notice that
\be
\inf_{x \in A^n} \sum_{(m,y) \in \ZZ_x} |Df^m(y)|^{-2} \geq C^{-1} \rho

We also have
\be
\inf_{x \in A^n} \sum_{(m,y) \in \ZZ_x} |Df^m(y)|^{-2} \Xi_2(f,y) \leq
\inf\sum_{(R^n(f))^{m'}(y')=y} |D(R^n(f))^{m'}(y')|^{-2} \sum_{(m'',y'') \in
\WW} |Df^{m''}(y'')|^{-2} \Xi_2(f,y'') \geq \inf_{x \in A^n}

We now consider the set $Y \subset A^{2n}$ that land in $A^n$ before passing
through $A^n$

\end{pf}

\begin{cor}

We have that $\eta_n \leq C \xi^2_n$.

\end{cor}

\begin{lemma}

We have that $\xi_n \leq C(\eta_n+\delta_n)$.

\end{lemma}

\begin{lemma}

If $\lim \eta_n>0$ then $\lim \delta_n=0$.

\end{lemma}

\begin{pf}

We are going to show that $\ln \eta_n \geq -C^{-1}\sum_{m<n} \delta_m$.
Notice that the inequality $\ln \eta_n \leq -C\sum_{m<n} \delta_m$ also
holds but will not be needed.
\end{pf}

\begin{lemma}

Let us consider a sequence of periodic points of renormalization $f_n$ with
combinatorics $[c_0,...,c_n]$ which converge to $f$.  Then
\be
\frac {\delta_n(f_n)} {\delta_n(f)},
\ee
\be
\frac {\eta_n(f_n)} {\eta_n(f)}
\ee
are bounded away from zero and infinity.

\end{lemma}

%\end{lemma}

\begin{thm}

Let $f$ be a Feigenbaum Julia set.  Then $|J(f)|=0$.

\end{thm}

\begin{pf}

Assume that there exists a Feigenbaum Julia set $f$ of positive Lebesgue
measure with combinatorics $[c_0,...]$.  Let $g$ be a fixed point of
renormalization with combinatorics $[c_*]$ with zero Lebesgue measure.
Let $f_{n,m}$ be a periodic point of renormalization of period
$n+m$ and combinatorics $[c_*,...,c_*,c_0,...,c_{m-1}]$.  We have:
\be
\liminf_{m \to \infty} \eta_{n+m}>0,
\ee
\be
\lim_{n \to \infty} \sup_{m} \eta_{n+m}=0,
\ee
\be
\lim_{m \to \infty} \delta_{n+m}=0.
\ee

Let $C>1$ be a constant that work for lemmas.  There exists $n$, $m$ so that
\be
\delta_{n+m}<\eta_{n+m}<\frac {1} {4 C^3},
\ee
and thus
\be
\eta_{n+m} \leq C \xi^2_{n+m} \leq C(C(\eta_{n+m}+\delta_{n+m}))^2 \leq
4 C^3 \eta_{n+m}^2.
\ee
\end{pf}
}
}

\section{Conservativity}\label{conservativity}

In order to complete the proof of Theorem~E we must show that some
Feigenbaum maps have a conservative $\delta_\crit$ conformal measure.
In this section we will provide a criterium for conservativity,
and in the next section we will show that such criterium is satisfied
by some Feigenbaum maps

%In this section we will provide a criterium for conservativity of
%$\delta_\crit$ conformal measures:
%complete the proof of Theorem E.

\begin{thm} \label {conserve}

Let $f$ be a Feigenbaum map\footnote {It is enough to assume that $f$ is $n$
times renormalizable.}.  If $\delta_\crit(f)>\delta_\crit(R^n(f))$ then
the $\delta_\crit(f)$ conformal measure on $J(f)$ is conservative.

\end{thm}

\begin{pf}

Let $\delta_\crit \equiv \delta_\crit(f)$ and let $\mu$ be the
$\delta_\crit$ conformal measure on $J(f)$.
Let $Q \subset A^n$ be the set of points that return to $A^n$.
By Lemma \ref{diss criterion},  $\mu$ is conservative if and only if $\mu(Q)=\mu(A^n)$.
We will show that if $\mu(Q)<\mu(A^n)$ then
\be \label {bound31}
E_\mu(\Xi_\delta^{[j+1]}(A^n \leftarrow U)|A^n) \leq C+\kappa
E_\mu(\Xi_\delta^{[j]}(A^n \leftarrow U)|A^n), \quad \delta \approx
\delta_\crit,
\ee
for some $C>0$, $\kappa<1$.  This implies that $E_\mu(\Xi_\delta(A^n
\leftarrow U)|A^n) \leq C/(1-\kappa)$ for $\delta \approx \delta_\crit$ which
contradicts the definition of $\delta_\crit$ and concludes the proof.
 
Through the argument we will skip the notation for truncation.
We start with the estimate
\be \label {bound32}
E_\mu(\Xi_\delta(A^n \leftarrow U)|A^n)=
1+E_\mu(\Xi_\delta(A^n \xleftarrow[U \setminus A^n]{} U
\setminus A^n)|A^n)+E_\mu(\Xi_\delta(A^n \xleftarrow[U
\setminus A^n]{+} A^n \leftarrow U)|A^n).
\ee

Let
\be
\nu_\delta=\int_{A^n} \sum_{(x_0,...,x_k=x) \in \FF} \frac {\delta_{x_0}}
{|Df^k(x_0)|^\delta} d\mu(x), \quad \FF=
\left (A^n \xleftarrow[U \setminus A^n]{+} A^n \right ).
\ee
Then $\nu_{\delta_\crit}=\mu|Q$.
Notice that
\be
\Xi_\delta(A^n \xleftarrow[U \setminus A^n]{+} A^n)=\Xi_\delta
(A^n \xleftarrow[U \setminus A^n]{} U^n \xleftarrow[U \setminus V^n]{}
A^n)+\Xi_\delta(A^n \xleftarrow[U \setminus V^n]{+} A^n).
\ee
The assumption of dissipativity of $\mu$ implies in particular that
$\sup\, \Xi_{\delta_\crit}(V^n \xleftarrow[U \setminus V^n]{+} A^n)<\infty$.
By hyperbolicity of the set of orbits avoiding $V^n$, and by Lemma~\ref
{continuity}, we have that
$\Xi_\delta(V^n \xleftarrow[U \setminus V^n]{+} A^n)$ is uniformly
bounded in a neighborhood of $\delta_\crit$.  Notice that
\be
\Xi_\delta(A^n \xleftarrow[U \setminus A^n]{} U^n)=
\Xi_\delta(A^n \xleftarrow{f_n} U^n),
\ee
which by hypothesis is uniformly bounded in a neighborhood of
$\delta_\crit$.  Thus
\be \label {bound2}
\Xi_\delta(A^n \xleftarrow[U \setminus A^n]{+} A^n)<C, \quad \delta \approx
\delta_\crit
\ee
and it follows
immediately (see Lemma~\ref {nu_delta})
that $\delta \mapsto \nu_\delta(N)$ is continuous
in a neighborhood of $\delta_\crit$ for any Borelian set $N$.
If $\mu(Q)<\mu(A^n)$ it then follows that
\be \label {bound21}
\frac {E_\mu(\Xi_\delta(A^n \xleftarrow[U \setminus A^n]{+} A^n
\leftarrow U)|A^n)} {E_\mu(\Xi_\delta(A^n \leftarrow
U)|A^n)}=\frac {E_{\nu_\delta}(\Xi_\delta(Q \leftarrow U)|A^n)} {E_\mu(\Xi_\delta(A^n \leftarrow
U)|A^n)} \leq \kappa<1, \quad \delta \approx \delta_\crit.
\ee

Similarly to (\ref {bound2}) we conclude that that
\begin{align} \label {bound1}
E_\mu(\Xi_\delta(A^n \xleftarrow[U \setminus A^n]{} &U
\setminus A^n)|A^n)=
E_\mu(\Xi_\delta(A^n \xleftarrow[U \setminus A^n]{} U^n)|A^n)+
E_\mu(\Xi_\delta(A^n \xleftarrow[U \setminus V^n]{} U \setminus V^n)|A^n)\\
\nonumber
&+E_\mu(\Xi_\delta(A^n \xleftarrow[U \setminus A^n]{} U^n
\xleftarrow[U \setminus V^n]{} U \setminus V^n)|A^n)<C,
\quad \delta \approx \delta_\crit.
\end{align}
Plugging (\ref {bound1}) and (\ref {bound21}) into (\ref {bound32}) we get
(\ref {bound31}).
\end{pf}

\section{Varying the dimension}\label{varying dim}

In this section we will give two proofs that for some Feigenbaum maps,
$\de_\crit(Rf)< \de_\crit (f)$. Together with Theorem~\ref{conserve},
it will justify the assertion of Theorem E that the $\de_\crit$-conformal
measure can be conservative (up to the proof of the existence of Lean
Feigenbaum Julia sets, which will be completed in \S~\ref {Yarr}).

  \subsection{Varying the dimension within the horseshoe}

% The space of Feigenbaum maps is forward invariant under renormalization.
Below  we will make use of the  continuity properties of the critical exponent 
(Lemmas \ref{lower} and \ref{continuity1}) and an observation that it is non-increasing
under the renormalization: 
\be\label{non-decreasing}
\de_\crit(Rf)\leq \de_\crit (f)
\ee

  \begin{lemma}

  Let $f$ be a Feigenbaum map.  If $\delta_\crit(f)=\delta_\crit(R(f))$ then
  $\delta_\crit(f)=\delta_\crit(R^n(f))$ for all $n \geq 0$.

  \end{lemma}

  \begin{pf}

  Assume that $\delta_\crit(f)=\delta_\crit(R(f))>\delta_\crit(R^2(f))$.
  Then the $\delta_\crit(R(f))$
  conformal measure on $J(R(f))$ is dissipative by the first relation (Proposition \ref{diss3})
   and conservative by the second one (Theorem \ref{conserve}) - contradiction.
  \end{pf}

Let us consider the horseshoe $\AAA\equiv \AAA_\MM$ associated to some fine collection $\MM$ 
of little Mandelbrot copies (see \S \ref{horseshoe}). 
Notice that  under the renormalization operator $R: \AAA\ra \AAA$
 the orbit  of a generic point $f\in \AAA$ is dense in $\AAA$. 

\begin{lemma} \label {hete}
Assume that $\orb_R (f)$ is dense in $\AAA$.
 Then 
$$ 
  \delta_\crit(f)=\max_{g \in \AAA} \delta_\crit(g).
$$

\end{lemma}

\begin{pf}
This follows from the lower semi-continuity of the critical exponent and (\ref{non-decreasing}) .
\end{pf}

  \begin{lemma}

  If the critical exponent is not constant over a renormalization horseshoe $\AAA$
  (associated with  some fine collection of Mandelbrot copies) then there exists $f \in \AAA$
  such that $\delta_\crit(f)>\delta_\crit(R f)$.

  \end{lemma}

  \begin{pf}
 
 Take two maps $f_+,f_- \in \AAA$ such that $\delta_\crit(f_+)>\delta_\crit(f_-)$.
  Since the critical exponent is lower semi-continuous,
   $f_+$ can be selected to be a periodic point of some period $p$ of $R$. 
  Let $f$ belong to the intersection of the unstable manifold of $f_+$ with the
  stable manifold of $f_-$.  Using semi-continuity again, we conclude that
$$
   \lim \delta_\crit(R^{-n p}(f)) \geq  \delta_\crit( f_+ ).
$$ 
On the other hand, by monotonicity (\ref{non-decreasing}) and uniform continuity on  hybrid classes
(Lemma \ref{continuity1}), 
 $$
 \lim  \delta_\crit(R^n f)= \lim \delta_\crit (R^n f_- ) \leq \delta_\crit(f_-).
$$
  \end{pf}

  \begin{example}

  Let $f_p$ be the real fixed point of renormalization of period $p$ with
  combinatorics closest to Chebyshev.
  By Lemma~\ref{de ge 1}, $\Xi_1(f_p,z)=\infty$. On the other hand by
Lemma~\ref{diss4} the $\delta_\crit$-conformal measure on
of $f_p$ is dissipative; hence by Lemma~\ref{diss1},
$\Xi_{\de_\crit}(f_p,z) < \infty$.
Hence  $\delta_\crit(f_p)>1$.\footnote{This also follows, from
Theorem~\ref{Zdunik} and Corollary \ref {HD = hyp}, in the
case of $p$ large which we use here (since in this case
the Julia set has zero Lebesgue measure).}
  
 Moreover, it is shown in \cite{AL} that
%we will show in Theorem~\ref {chebthe} that
$\lim  \delta_\crit(f_p)=1$. %\note{ref} 
 Thus we can choose $p_1$ and $p_2$ such that
  $\delta_\crit(f_{p_1})>\delta_\crit(f_{p_2})$, so the horseshoe associated
  to the corresponding Mandelbrot copies satisfies the hypothesis of the
  previous lemma.

  \end{example}

\subsection{Varying the dimension within a hybrid class}\label{varying dim1}

Recall that $\HH_f$ stands for the hybrid class of a quadratic-like map $f$ (see \S \ref{q-l maps}). 
We will study now how the dimension  $g\mapsto \HD(J(g))$ vary over  $\HH_f$. 

\begin{lemma}\label{sup HD} 

Let $f$ be a Feigenbaum map.  Then 
$$ 
          \sup_{g\in \HH_f} \HD(J(g)) = 2.
$$

\end{lemma}

\begin{pf}
Let us take a little square $\De\Subset \V \sm \OO$ intersecting the Julia set.
It is easy to find   four branches $g_i\equiv f^{-n_i}_i | \De $ of inverse iterates such that 
$D_i\equiv g_i(\De)\Subset \De$ and $f^{n_i}$ is the first return map from $D_i$ to $\De$. 
Let $S=\cup D_i$ and let $G: S \ra \De$ be the associated Bernoulli map, $G|D_i = f^{n_i}$. 
Let $Q$ be the maximal $G$-invariant set contained in  $S$.
By the Schwarz Lemma, $Q$ is hyperbolic for $G$, so that the forward saturation $ \orb (Q)$ is 
a hyperbolic set for $f$. 

Let us consider a diffeomorphism  $H: \C  \to \C$ such that 
 $\tl \De\equiv h(\De)$ is the standard square $ [0,1]\times [0,1]$ of size 1
 and $\tl D_i \equiv H(D_i)$ are rescaled standard squares  of size close to $1/4$
(and we can select  $H$ to be affine  on $\C\sm \De$). 
 Let $\tl S=\cup \tl D_i$. 
 Consider a piecewise affine  $\tl G : \tl S \to \tl \De $ 
 that affinely maps each $\tl D_i$ onto $\tl \De$. 
 Then the maximal invariant set $\tl Q$ of $\tl G$ has Hausdorff dimension  $2-\eps$.

By the standard pullback argument, the maps  $G: S \to \De $ and $\tl G : \tl S  \to \tl \De $
are conjugate by some qc homeomorphism $h$ which is affine on $\C\sm \De$. 
Let us consider the corresponding $G$-invariant Beltrami differential  $\mu = \dibar h / \di h$. 
%We let $g=h_\mu \circ f \circ h_\mu^{-1}$.

Let $f: \u \ra \V$ and let $W^m = f^{-m} \u$.
For $m \geq 0$, define a Beltrami differential $\mu_m$ as follows.  
Let $\mu_m = \mu$  on $ \C \setminus W^m$. Pull it back by the iterates of $f$
to the domain $W^m \setminus J(f)$.  Finally, let  $\mu_m |J(f)=0$.
Since $G$ is composed from iterates of $f$, this Beltrami differential
 is simultaneously invariant with respect to 
$f:W^m \to W^{m-1}$ and $G: S\ra \De$.   
Since conformal pullbacks preserve the norm of  Beltrami differentials,
$\|\mu_m\|_\infty \leq \|\mu\|_\infty$. 
 
We assume now that $\area (J(f))=0$, for otherwise there is nothing to prove.
In this case $\mu_m \to \mu$ almost everywhere.

Let $h_m$ be solution of the Beltrami equation $\dibar h_m / \di h_m = \mu_m$ coinciding
with $h$ at three points. Then $h_m \to h$ uniformly on $\bar\C$.

Consider a quadratic-like map 
$$
  g_m=h_m \circ f \circ h_m^{-1}: h_m(W^m) \ra h_m(W^{m-1}). 
$$
and the corresponding expanding map 
$$
   G_m= h_m \circ G \circ h_m^{-1}: h_m(S) \ra h_m (\De)  .
$$  
Since the latter maps uniformly converge to $\tl G$, $\HD(h_m(Q)) \to \HD (\tl Q)= 2-\eps$. 
The conclusion follows since the forward saturation of $h_m(Q)$ under $g_m$ is a hyperbolic subset of $J(g_m)$.
\end{pf}

\begin{rem}
In case of the Julia set of positive area, it would be interesting to
obtain the same result for the hyperbolic dimension  $\HD_\hyp (g)$ rather than $\HD(J(g))$. 
\end{rem}

Let us  consider a fine family $\MM$ of little Mandelbrot copies and  
a Feigenbaum map $f$ with combinatorics $\bar M_+ = \{M_0, M_1, \dots\}$, $M_i\in \MM$. 
We say that $f$ has  {\it recurrent combinatorics}  
if the sequence $\bar M_+$ is recurrent. 
Such a sequence  can be extended to  the ``past''  as a bi-infinite recurrent sequence $\bar M = \{M_k\}_{k\in Z}$.
Take the map $f_*$ in the renormalization horseshoe $\AAA_\MM$ with combinatorics $\bar M $.
It  is recurrent under  the renormalization  dynamics  and the $R$-orbits of $f$ and $f_*$ are forward asymptotic.
We say that $f_*$ is ``associated'' with $f$. 

\begin{thm}
   Let $f$ be a Feigenbaum map with recurrent combinatorics,
and let $f_*$ be an associated map in the horseshoe $\AAA_\MM$. 
Then the function $h: g\mapsto \HD(J(g))$, $g\in \HH_f$, is either identically equal to 2, 
or 
$$
   \operatorname{Im} h = [\HD(J(f_*)), 2).
$$ 
\end{thm} 

\begin{pf}
Since qc maps preserve the property $\HD(J)=2$, the function $h$ is either identically equal to 2,
or always less than 2. Assume that the latter happens.

Since hybrid classes are connected and the Hausdorff dimension varies continuously under qc deformations,
$\operatorname{Im} (h)$ is an interval.  By Lemma \ref{sup HD}, it must be either  $[a, 2)$ or $(a,2)$.

Since the function $h$ is a decaying Lyapunov function for the renormalization (i.e., $h(Rg) \leq h(g)$),
 and by the recurrence assumption the orbit $R^n f$ accumulates on $f_*$, we conclude that 
$h(f)\geq h(f_*)$. Hence 
$$
   h(f_*)=\min_{g\in \HH_f} h(g),
$$
 and we are done.   
\end{pf}

Take, for instance, a renormalization fixed point $f_*$  of lean type.  Then $\HD(J(f_*))<2$
and there is a map $f\in \HH_{f_*}$ with $\HD(J(f))> \HD(J(f_*))$. Since $R^n f\to f_*$ and the critical exponent
is continuous on hybrid classes, there is an $n$ such that $\HD(R^{n+1} f) < \HD (R^n f)$,
so that the $\de_\crit$-conformal measure of $R^n f$ is conservative.

\comm{
\comm{
\begin{pf}

There exists $m \geq 4$ and a sequence of disjoint Jordan domains
$D_1,...,D_m$ such that $f|D_i$ is univalentf(D_i)=D_{i+1}$,
$1 \leq i \leq m-1$, $f(D_m)$ is a Jo

A Markov system is a pair $(g,S)$ where $S \subset \C$ is a
finite union of Jordan domains $S_i$, $1 \leq i \leq m$ with disjoint closure and
analytic boundary, $g:S \to \C$ is a holomorphic map such that $g(S)$ is a
finite union of Jordan domains with disjoint closure and analytic boundary,
$\overline S \subset g(S)$ and for every $1 \leq i \leq m$, $g|S_i$ is a
univalent map onto some component of $g(S)$.

Let $Q=\cap g^{-n}(S)$.  A Markov system is said to be non-trivial if
$Q$ is a Cantor set.  In this case we have $0<\HD(Q)<2$.

The combinatorics of a Markov system is a matrix $a_{ij}$,
$1 \leq i,j \leq m$ where $a_{ij}=1$ if $S_j \subset f(S_i)$ and
$a_{ij}=0$ otherwise.  Such a matrix will be considered up to permutation
of $\{1,...,m\}$.

Let $m \geq 4$, and let $b_{ij}(m)=1$ if $j-1be the matrixsymbols if
$m \geq n$

A Markov system is said to be essentially a horseshoe with $n$ symbols if
$f(S^

\begin{lemma}

Let $Q \subset \C$ be a hyperbolic Cantor set for some holomorphic

Let $m \geq 4$, and let $a_{It is easy to construct examples of Markov systems with Hausdorff

Two Markov systems $(g_1,S^1,T^1)$ and
$(g_2,S^2,T^2)$ are said to be combinatorially equivalent if there is a
homeomorphism $h:(T^1,S^1) \to (T^2,S^2)$ such that $g_1(S^1_i)=T^1_j$
if and only if $g_2(h(S^1_i))=h(T^1_j)$.

Two Markov systems are said to be
qc-equivalent if $h$ can be chosen quasiconformal and $h \circ g_1|S^1=
g_2 \circ h$.

A Markov system is said to be non-trivial if $\cap g^{-n}(S)$ is a Cantor
set.

\begin{lemma}

Any non-trivial Markov system is qc-equivalent to a

$\overline S
and $g$ is a holomorphic map defined on a neighborhood of $M$
such that $g(\partial M) \cap M=emptyset$, $g|M_i$ is a
diffeomorphism for $1 \leq i \leq m$, and the maximal invariant set
$Q_{g,M}$ of $g|M$ is a Cantor set which intersects each $M_i$.
If $(g,M)$ is a Markov system, we let $M^k=(g|M)^{-k}(M)$.  Notice that
$(g,M^k)$ is also a Markov system for $k \geq 0$.  Two Markov systems
$(g,M)$ and $(g',M')$ are said to be quasiconformally equivalent if there
exists a quasiconfomal map $h:\C \to \C$ such that $h(M)=M'$ and $h \circ
g|M=g' \circ h$.
}}

\begin{lemma}

Let $(g,M)$ be a Markov system.  There exists $k>0$ such that for every
$\epsilon>0$, there exists a Markov system $(g',M')$ which is
quasiconformally equivalent to $(g,M^k)$
such that $HD(Q_{g,M'}) \geq 2-\epsilon$.

\end{lemma}

\begin{lemma}

Let $f:U \to V$ be a quadratic-like map.  Then there exists $M \subset U$
such that $(f,M)$ is a Markov system.

\end{lemma}

\begin{thm}

Let $f:U \to V$ be a quadratic-like map with connected Julia set.
For every $\epsilon>0$, there exists a quadratic-like map $g$ which is
hybrid conjugate to $f$ and such that
$HD(K(g))>2-\epsilon$.\footnote {We believe
we can actually get $HD(J(g))>2-\epsilon$.}

\end{thm}

\begin{pf}

We may assume that $K(f)$ has zero Lebesgue measure, otherwise $HD(K(f))=2$.
Let $M \subset U$ be such that $(f,M)$ is a Markov system.  We may assume
that $(f,M)$ is equivalent to a Markov system $(g,M')$ such that
$HD(Q_{g,M'})>2-\epsilon$ via a quasiconformal map $h:\C \to \C$ fixing
$\{0,1,\infty\}$.  Let
$U^l=f^{-l}(U)$.  Let $D^l=U^l \cup \inter M$.  Let $\mu_l:\C \to \C$ be the
Beltrami differential such that
\begin{enumerate}
\item $\mu_l|\C \setminus D^l=\frac {\op h} {\partial h}$,
\item If $x \in D^l \setminus J(f)$
and $s \geq 1$ is minimal with $f^s(x) \in \C \setminus
D^l$ then $\mu_l(x)$ is the pullback of $\frac {\op h} {\partial h}(f^s(x))$
by $f^s$,
\item If $x \in J(f)$ then $\mu_l(x)=0$.
\end{enumerate}

Let $h_l:\C \to \C$ be such that $\frac {\op h_l} {\partial h_l}=\mu_l$ and
fixing $\{0,1,\infty\}$.  Then
$g_l=h_l \circ f \circ h_l^{-1}:h_l(D^l) \to \C$
is holomorphic.  Moreover, $h_l \to h$ uniformly, so $g_l$ form a sequence
of holomorphic maps converging uniformly to $g$ on a neighborhood of
$Q_{g,M'}$.  Thus $HD(h_l(Q_{f,M})) \to HD(Q_{g,M'})$, and for $l$ big,
$HD(h_l(J(f))) \geq HD(h_l(Q_{f,M}))>2-\epsilon$.  Notice that $g_l:U^l \to
U^{l-1}$ is a quadratic-like map and
$h_l(J(f))=J(g_l:U^l \to U^{l-1})$ so the result follows.
\end{pf}
}

\section{Interpolation argument}\label{interpolation}

Let us now prove Theorem B.

\begin{thm}

Let $\AAA$ be the horseshoe associated to some fine collection of Mandelbrot copies. 
If there exist $g_+,g_0 \in \AAA$ such that $\area (J(g_+))>0$ and
$\area (J(g_0))=0$, then for a generic $f\in \AAA$, $\delta_\crit(f)=2$ and $\area (J(f))=0$.

\end{thm}

\begin{pf}
Any map in the stable manifold of $g_0$ has a Julia set of zero area. 
Since the stable manifold of $g_0$ is dense in $\AAA$ and $\area J(g)$
depends upper semi-continuously  on $g\in \AAA$, 
it follows that a generic map in $\AAA$ has a Julia set of zero area.

Fix $s>0$ and let $f_{[s]} \in \AAA$ be such that $R^s(f_{[s]})$ belongs
to the intersection of the stable
manifold of $g_0$ and the local unstable manifold of $R^s(g_+)$,
so that the combinatorics of $f_{[s]}$ coincides with that of $g_+$ 
on many top renormalization levels 
and coincides with that of $g_0$ on all sufficiently deep renormalization levels.
Let us denote by $\phi_{n,s}$ the first ``univalent''
return map to the disk $V^{n,s}\equiv V^n(f_{[s]})$ (see \S \ref{nfd}), that
is, $\phi_{n,s}(z)=f^k_{[s]}(z)$ if $k$ is smallest possible such that
there is a univalent branch of $f^{-k}_{[s]}|V^{n,s}$ onto a neighborhood of
$z$.  Since the maps  $\phi_{n,s}$ are Markov  and have bounded distortion  
(see Remark \ref{bounded distortion}),
the critical exponent $\delta_\crit(\phi_{n,s})$ makes sense 
(i.e., it is independent of the base point $z\in V^{n,s}$ in the Poincar\'e series $\Xi_\de(\phi_{n,s},z))$. 
Obviously $\delta_\crit(f_{[s]}) \geq \delta_\crit(\phi_{n,s})$.
Let us show that we can choose $n$, $s$ so that $\delta_\crit(\phi_{n,s})$
is close to $2$,
as this will imply the result by Lemma \ref {hete}.

Let $D^{n,s}_k$ be the domain of $\phi_{n,s}^k$.  We may assume that $|D\phi_{n,s}(x)|>2$ for
$x \in D^{n,s}_1$, see Remark~\ref {small domains}.
%By construction of $V^n$, \note{should be listed}
%each component of $D^n_1$ has diameter much smaller than $\diam(V^n)$, 
%and by bounded distortion we get $|D\phi_n(x)|>2$ for every $x \in D^n_1$.

Let us show that we can choose $n$ and $s$ so that $p(D^{n,s}_1|V^{n,s})$ is close
to $1$.  Since $J(g_+)$ has positive Lebesgue measure, its critical point is a density point
of $J(g_+)$.  Given $\xi>0$, we can choose $\rho>0$ so that $p(J(g_+)|\D_\rho)>1-\xi$.
Take the minimal  $n>0$ such that $\diam(V^{n,s}) \leq \rho/2$ for all $s$ sufficiently big.
Then for $s$ big, $\diam(V^{n,s}) \asymp \rho$ and so
$V^{n,s} \supset \D_{q \rho}$ for some absolute $q>0$.  

Let $A^{n,s}$ be defined in the natural way (see \S~\ref {nfd}).
Let $L^s \subset \D_\rho$ be the set
of points whose positive orbit under $f_s$ intersects $\D_{q \rho}$.  Notice that
$L^s\cap A^{n,s} \subset D^{n,s}_1 \cap A^{n,s}$.  Since the positive orbit of almost every
$x \in J(g_+)$ (under iterates of $g_+$) accumulates on $0$,
$\liminf L^s \supset J(g_+) \cap \D_\rho$ up to a set of zero area.  Thus we can choose
$s$ so that $p(L^s|\D_\rho)>1-\xi$.  Hence $p(D^{n,s}_1|A^{n,s})>1-C\xi$.
By bounded distortion,
\be \label {Cxi}
p(D^{n,s}_1|B^{n,s})>1-C\xi \quad \text{where}\quad  B^{n,s}= \bigcup_{k \geq 0} f^{-k}_{[s]}(A^{n,s}).
\ee
Since $f_{[s]}$ is hybrid equivalent to  $g_0$, $J(f_{[s]})$ has zero area. 
 Thus $B^{n,s} =V^{n,s}$ modulo a set of zero area,
 so that  (\ref {Cxi}) is equivalent to $p(D^{n,s}_1|V^{n,s})>1-C\xi$ as desired.

Fix $\epsilon>0$.  Notice that if $p(D^{n,s}_1|V^{n,s})$ is close to $1$ then
$p(D^{n,s}_k|V^{n,s})>(1-\epsilon)^k$.  We have
\begin{align} \nonumber
\int_{V^{n,s}} \Xi_{2-\gamma}(\phi_{n,s},x) dx&=
\int_{V^{n,s}} \sum_{k \geq 0} \sum_{\phi^k_{n,s}(y)=x} |D\phi^k_{n,s}(y)|^{\gamma-2} dx\\
\nonumber
&\geq \int_{V^{n,s}} \sum_{k \geq 0}
\sum_{\phi^k_{n,s}(y)=x} |D\phi^k_{n,s}(y)|^{-2}
2^{k \gamma} dx=\sum_{k \geq 0} |D^{n,s}_k| 2^{k \gamma}\\
\nonumber
&\geq |V^{n,s}| \sum_{k \geq 0} 2^{k\gamma} (1-\epsilon)^k.
\end{align}
Thus the critical exponent $2-\gamma$ of $\phi_{n,s}$ satisfies
$2^\gamma (1-\epsilon) \leq 1$, so by choosing $\epsilon$ small we can get
$\delta_\crit(\phi_{n,s})$ arbitrarily close to $2$.
\end{pf}

\begin{rem}

Another interpolation argument allows one to show that,
under the hypothesis of the
previous theorem, there exists $f \in \AAA$ with a dense orbit (or which is
periodic) and with $\area J(f)>0$.  Notice that if $f$ has a dense orbit
then it has automatically $\delta_\crit(f)=2$
by the previous theorem and Lemma~\ref {hete}, and if $f$ is periodic then
$\delta_\crit(f)<2$ (since we are in the Lean case).

\end{rem}

\comm{
\begin{lemma}

Let $f$ be a recurrent

A subhorseshoe is a subshift of finite type which is transitive.

\begin{lemma}

Let $f$ be a Feigenbaum map.  Then the stable manifold of $f$
intersects a subhorseshoe.

\end{lemma}

\begin{pf}

It is enough to show that the $\omega$ limit of $f$ is contained in
a subhorseshoe, since the latter is automatically a locally maximal invariant
set.
\end{pf}

\begin{lemma}

Let $f$, $g$ be two Feigenbaum maps.  Then there exists a subhorseshoe
$\Sigma$ which intersects the stable manifolds of $f$ and $g$.

\end{lemma}

\begin{pf}

The stable manifold of any
Feigenbaum map is foliated by hybrid classes and is dense in the boundary
of the connectedness locus.

Let $\Sigma_f$, $\Sigma_g$ be subhorseshoes containing the $\omega$ limit
set of $f$ and $g$.  Let $f_0$ be a map in the stable manifold of $f$ which
is contained in the unstable manifold of $\Sigma_g$, and define $g_0$
analogously.

Given any neighborhood of $\Sigma_f$, $\Sigma_g$, and the full orbits of
$f_0$ and $g_0$, it is straightforward to find an orbit whose $\omega$ limit
contains $f_0$ and $g_0$.

Let us consider a neighborhood of $\Sigma_f$, $\Sigma_g$, and the full
orbits of $f_0$, and $g_0$.

\begin{rem}

Another (somewhat more technical) interpolation argument shows that
in the setting of the previous lemma, there is a dense set
of $f \in \Sigma$ such that $|J(f)|>0$ and $\delta_\crit(f)=2$.

\end{rem}

\comm{
\begin{pf}

It is enough to find a single $f$ with the required properties, as those
will also hold in the stable manifold of $f$.

We will actually show the following: given $\epsilon>0$ and $f \
}
}

\section{Lean critical points}\label{Yarr}

In this section we will show that there are many Feigenbaum
Julia sets with zero area (Yarrington's thesis \cite{Y}).  We will
restrict ourselves to a particular class of combinatorics (to be defined
below), for which the estimates are simpler.  As a consequence of our
estimates, we will also conclude that the Lean case exists, completing the
proof of the first assertion of Theorem~A.

\comm{
The specific measure estimates
obtained will be used further for perturbative estimates of Poincar\'e
series for $\delta$ close to $2$.  Notice that this class of
combinatorics does not formally include combinatorics close to Chebyshev,
but the argument can easily be adapted to this case.
}

Let $F:\cup \U_j \to \U$ be a holomorphic map with the following
properties:
\begin{enumerate}
\item $\cup \U_j \subset \U$,
\item For every $j \neq 0$, $F:\U_j \to \U$ is univalent,
\item $f \equiv F:\U_0 \to \U$ is a quadratic-like map with
connected Julia set.
\end{enumerate}
We let $\U^k=f^{-k}(\U_0)$.

%We let $\U_k=(f_0|\U^0_0)^{-k}(\U_0)$.  Let $f_k:\cup \U^j_k \to \U_k$ be
%the first return map to $\U_k$, where we let $\U_{k+1}=\U^0_k$.
We say that $F$ is $\kappa$-good if
\begin{enumerate}
\item $\mod(\U \setminus \overline \U_0) \geq \kappa$,
\item The domain of the first return map to $\U_0$ has density at most
$\kappa^{-1}$.
\end{enumerate}

Let $P$ be the (monic)
quadratic polynomial which is hybrid equivalent to $f:\U_0
\to \U$.  Denote by $D_r$ the domain bounded by the equipotential of
height $r$ for $P$.  We construct the standard Yoccoz puzzle of
$P$ with the help of the equipotential at height $4$, see for
instance \cite {puzzle}.  Let $\hat T_1 \supset \hat T_2 \supset ...$ be
the principal nest of $P$ obtained in the standard way (see \cite {puzzle}).

We will make the following combinatorial assumptions:
\begin{enumerate}
\item $P$ belongs to a selected finite family of Misiurewicz wakes,
   or, equivalently,  $\hat T_1$ is mapped 
   to the top level puzzle piece by an iterate $F^n$ with a bounded $n$
   (see \cite{parapuzzle,Sc});   
\item $P$ is renormalizable with no central cascades: there exists
$\tau(P)>0$ (the height of $P$) such that the first return of $0$
to $\hat T_k$ belongs to $\hat T_{k+1}$ if and only if $k \geq \tau(P)$.
% \note{Misha: could you please give suitable definitions and references for
% those concepts?}
\end{enumerate}
In what follows, the constants involved will depend on the choice of the
Misiurewicz wake.

We will need the following simple plane topology argument.

\begin{lemma}

Let $W_i \subset \C$, $i \geq 0$, be Jordan domains such that $\diam W_j \to 0$.
Assume that $W_0\cap W_i \not=\emptyset$, $i>0$
while $\overline W_i\cap\overline W_j = \emptyset $ for $i>j>0$. 
Consider a compact set  $W \subset \C$  obtained by taking the closure of  $\cup W_i$  and then
filling in the holes.  Then $W$ is a Jordan disk.

\end{lemma}

\begin{pf}

We may assume that $W_0$ is not contained in any $W_i$, for otherwise $W=W_i$.
We can also discard the $W_i$, $i>0$ that are contained in $W_0$.
It follows that $\#(\partial W_i \cap \partial W_0) \geq 2$ for $i \geq 0$.

%Since $W$ is a full compact set, $\partial W$ is connected.
Note that  $\partial W \subset \overline {\cup \partial W_i} $. 
Moreover,  if a  sequence of domains $W_i$ accumulates on  some $x \in \partial W$  then $x \in \partial W_0$.  
Hence  $\partial W \subset \cup \partial W_i$. 
%Moreover, if $x \in \partial W \setminus \partial W_0$, say $x \in \partial W_i$,
%then there is a neighborhood of $x$ in $\partial W_i$ which is contained in $\partial W$.

It follows that each connected component of $\partial W \setminus \partial W_0$ is an open arc $\gamma_k$
contained in some $\partial W_i$, whose endpoints are contained in $\partial
W_0$. (Warning: there can be infinitely many $\gamma_k$'s in one $W_i$).
   Notice that the endpoints of $\gamma_k$ cannot coincide (otherwise it
would bound some $W_i$ such that $\partial W_i \cap \partial W_0$ is a
point), so $Q_k\equiv \overline \gamma_k \setminus \gamma_k$ consists of
two points.  Since $\diam W_i \to 0$ we also have $\diam \gamma_k \to 0$.

Let $\delta_k$ and $\delta_k'$ denote the two components of $\partial W_0
\setminus Q_k$.  One of them, say, $\delta_k$, is such that $\gamma_k
\cup Q_k \cup \delta_k$ bounds a domain contained in $W$
(since $W$ is a full compact set) which is disjoint from
$W_0$.  This shows that each $\delta_k$ is contained in $\inter W \cap \partial W_0$.
Hence $\partial \delta_k\cap \delta_i= \emptyset$ for $k\not=i$, so that all the arcs $\delta_k$ are disjoint.
Since they  lie on the Jordan curve $\partial W_0$,  $\diam \delta_k \ra 0$.

Let $h_k:\delta_k \to \gamma_k$ be a homeomorphism fixing the endpoints,
and let $h:\partial W_0 \to \C$ be defined by $h|\delta_k =h_k$
and $h|\partial W_0 \setminus \cup \delta_k=\id$.  Since $\diam \gamma_k
\to 0$ and $\diam \delta_k\to 0$, $h$ is continuous.  It is also
injective by definition, so it is a homeomorphism onto its image.  Notice
that $h(\partial W_0) \supset \partial W$ by construction.  On the other
hand, $\partial W$ can not be a proper subset of
$h(\partial W_0)$ (otherwise it would not disconnect the plane).  Thus
$\partial W= h(\partial W_0)$ is a Jordan curve as was asserted.
\end{pf}

\begin{rem} \label {exten}

Notice that the proof constructs an explicit homeomorphism between $\partial
W^0$ and $\partial W$, which is close to the identity when
$\sup \diam(W_i)$ is small, and which can be extended to a homeomorphism of
$\overline W^0$ to $\overline W$ which is also close to the identity.

\end{rem}

We will say that $F:\cup \U_j \to \U$ is normalized if
$F(z)=F(0)+z^2+O(z^3)$ near $0$.

\begin{lemma} \label {eta}

For every $\eta>0$, there exists $\kappa=\kappa(\eta)>0$ with the following
property.  Let $F$ be $\kappa$-good and normalized.
There exists a homeomorphism $h:D_4 \to \C$ which
conjugates $P:D_2 \to D_4$ to $F:h(D_2) \to h(D_4)$ such that
\begin{enumerate}
\item $h$ is $\eta$ close to the identity,
\item If $S$ is a non-central component of the domain of the first return
map to $\U_0$, $S'$ is a component of a preimage of $S$ by some iterate
of $f$, and $Q$ is a Yoccoz puzzle piece for $P$
intersecting $h^{-1}(S')$, then $h^{-1}(S') \subset Q$.
\end{enumerate}

\end{lemma}

\begin{pf}

Let $M^0$ be the domain bounded by the equator of $\U \setminus \U_0$, and
$M^k=f^{-k}(M^0)$.
If $\kappa$ is big, then $\partial M^0$ is close to a circle centered at
$0$, and $f|M^1$ is essentially quadratic.
It follows that there exists $r$ and a quasiconformal map
$\hat h:\overline D_r \setminus D_{r^{1/2}} \to \overline
M^0 \setminus M^1$ with dilatation close to $1$ which conjugates
$f|\partial M^1$ and $P|\partial
D_{r^{1/2}}$ (notice that $r \geq \kappa$ and indeed
$\ln r$ is close to $2\mod(M^0 \setminus M^1)$).  By the pullback
argument, we obtain an extension $\hat h:\overline D_r \to \overline
M^0$ with the same dilatation which is a hybrid conjugacy.  By the
normalization at the origin, $\hat h|D_{256}$ is uniformly close to the
identity in the Euclidean metric, provided $\kappa$ is big.

Let $\SSS$ be the collection of non-central components of the domain of the
first return map to $\U^0$, and let $\SSS'$ be the collection of components
of preimages of some $S \in \SSS$ by some iterate of $f$.
Notice that if $S \in \SSS$ then there exists
$\tilde S \subset \U^0 \setminus \U^1$ such that the first
return map $\phi:S \to \U^0$ extends to a univalent map $\phi:\tilde S \to
\U$.  If $S' \in \SSS$ is such that $f^k(S')=S$
then $S' \subset \U^k \setminus \U^{k+1}$ and the map
$f^k:S' \to S$ extends to a bigger domain
$\tilde S' \subset \U^k \setminus \U^{k+1}$ so that
$f^k:\tilde S' \to \tilde S$ is univalent.
Notice that $\mod(\tilde S' \setminus S')=
\mod(\U \setminus \U^0) \geq \kappa$, so the hyperbolic diameter of $S'$
in $\tilde S'$ is very small.  By the Schwarz Lemma, since $\tilde S'
\subset \U^k \setminus \U^{k+1} \subset M^k \setminus M^{k+2}$,
the hyperbolic diameter of $S'$ in $M^k \setminus
\overline M^{k+2}$ is also very small.

Notice that if
$S' \in \SSS'$ intersects $\partial \hat h(D_4)$ then
$M^k \subset \C \setminus
\overline {\hat h(D_4)}$ and $M^{k+2} \subset \hat h(D_4)$,
so $M^k \setminus M^{k+2} \subset \hat h(D_{256}) \setminus
\hat h(D_{\sqrt 2})$, and the hyperbolic diameter of
$S'$ in $\hat h(D_{256}) \setminus \hat h(D_{\sqrt 2})$
is small.  Let $B$ be the union of $\hat h(D_4)$ with all $S'
\subset \SSS'$ which intersect $\hat h(D_4)$, and let $B^0$ be obtained from
$\overline B$ by filling in the holes and taking the interior.
%Let us consider the closure $B$ of
%the union of $\hat h(D_4)$ with all $S' \in \SSS'$ which intersect
%$\hat h(D_4)$, and let $B^0$ be the complement of the closure of the
%unbounded component of $\overline B$.  Then $B^0$ is a simply
%connected domain, and indeed $\partial B^0$ is a Jordan curve
%(it is enough to use that the $S' \in \SSS'$ have bounded shape).
By the previous lemma, $B^0$ is a Jordan domain.
Moreover, $\partial B^0$ is close to
$\partial \hat h(D_4)$ in the hyperbolic metric of
$\hat h(D_{256}) \setminus \hat h(D_{\sqrt 2})$.
Let $B^k=f^{-k}(B^0)$.  Then (using Remark~\ref {exten}) there exists
a homeomorphism $\tilde
h:\overline D_4 \setminus D_2 \to \overline B^0 \setminus B^1$
that conjugates $P|\partial D_2$ to $f|\partial B^1$, such that $\tilde
h \circ \hat h^{-1}:\overline {\hat h(D_4)} \setminus \hat h(D_2) \to
\overline B^0 \setminus B^1$ is
uniformly close to the identity in the hyperbolic metric of $\hat h(D_{256})
\setminus \hat h(D_{\sqrt 2})$.

The image by $\tilde h$ of the
external rays $\gamma$ landing at the $\alpha$ fixed point of $P$
divide $\overline B^0 \setminus B^1$ into finitely many (say, $k$)
topological rectangles.  Applying the previous lemma $k-1$ times,
we obtain a new homeomorphism
$h:\overline D_4 \setminus D_2 \to \overline B^0 \setminus B^1$, such that
$h(\gamma \cap \overline D_4 \setminus D_2)$ is disjoint from $\cup_{S' \in
\SSS'} S'$ for all external rays $\gamma$.
We may still require (see Remark~\ref {exten}) that
$h \circ \hat h^{-1}:\overline {\hat h(D_4)} \setminus
\hat h(D_2)$ is uniformly close to the identity in the hyperbolic
metric of $\hat h(D_{256}) \setminus \hat h(D_{\sqrt 2})$.

\comm{
For each external ray $\gamma$ landing at
the dividing repelling fixed point of $P$, let us modify $\tilde h$ in a
small neighborhood of $\gamma \cap (\overline D_4 \setminus D_2)$ to obtain
a new homeomorphism
$h:\overline D_4 \setminus D_2 \to \overline B^0 \setminus B^1$, such that
$h(\gamma \cap \overline D_4 \setminus D_2)$ is disjoint from $\cup_{S' \in
\SSS'} S'$ (one modifies the curve $h(\gamma \cap \overline D_4
\setminus D_2)$ to make it ``go around'' components $S' \in \SSS'$
in the same way we changed $\partial \hat h(D_4)$ to $\partial B^0$).
We may still require that
$h \circ \hat h^{-1}:\overline {\hat h(D_4)} \setminus
\hat h(D_2)$ is uniformly close to the identity in the hyperbolic
metric of $\hat h(D_{256}) \setminus \hat h(D_{\sqrt 2})$.
}

Using the pullback argument, we can extend $h$ to a homeomorphism
$h:\overline D_4 \setminus K(P) \to \overline B^0 \setminus K(f)$.  Since
$P|D_2 \setminus K(P)$, $f|B^1 \setminus K(f)$ are unbranched coverings, we
conclude that $h \circ \hat h^{-1}:\overline {\hat h(D_4)} \setminus
K(f) \to \overline B^0 \setminus K(f)$ is uniformly close to the identity
in the hyperbolic metric of $\hat h(D_{256}) \setminus K(f)$.  In particular,
$h \circ \hat h^{-1}$ extends through $K(f)$ as the identity, so
$h$ extends to a
homeomorphism $h:\overline D_4 \to \overline B^0$ which is uniformly close
to $\hat h$ in the hyperbolic metric of $\hat h(D_{256})$.  The map $h$ has the
desired properties by construction.
\end{pf}

\comm{
Let us Let $\hat B$ be the union of $h(D_4)$ and
all $S'$ that intersect $h(D_4)$.  It may happen that $\hat B$ is not simply
connected, so let $B$ be obtained from $B'$ by filling in any holes.  It
follows that $B$ is a Jordan domain containing $h(D_4)$ whose boundary
avoids all puzzle pieces $S'$.e obtaineLet us approximate $h(\partial D_4)$ by a
Jordan curve $\partial B$ that goes around It follows that we can approximate $\partial
h(D_4)$ by a Jordan curve that goes aroundIt follows that there exists a Jordan
disk $B$ close to $h(D_4)$, and a homeomorphism
$h:\overline D_4 \setminus D_2 \to \overline B \setminus (f|U_0)^{-1}(B)$,
such that
\begin{enumerate}
\item $B$ contains $D_4$ and $\partial B$ avoids all $S'$,
\item $h$ conjugates $p|\partial D_2$ and $f|\partial (f|U^0)^{-1}(B)$,
\item If $\gamma$ is an external ray landing at the dividing
repelling fixed point of $p$ then $h(\gamma \cap \overline D_4 \setminus
D_2)$ avoids all $S'$,
\item $\hat h \circ h^{-1}$ does not move points much in $D_{256}$,
that is, for any $x \in \overline B \setminus (f|U^0)^{-1}(B)$,
$\hat h(h^{-1}(x))$ is close to $x$ in the hyperbolic metric of $D_{256}$.
\end{enumerate}

Since the hyperbolic metric of $D_{256}$ is comparable with the hyperbolic
metric of $D_{256} \setminus K(f)$ in a neighborhood of $\overline D_4
\setminus D_2$, the last property implies that $\hat h \circ h^{-1}$
does not move points much in the hyperbolic metric of
$D_{256} \setminus K(f)$

Let us extend $h$ to a homeomorphism conjugating $p|\overline D_4 \setminus
K(p)$ and $f|\overline
B \setminus K(f|U^0)$.  It follows from the Schwarz Lemma
that $\hat h \circ h^{-1}$ does not move points much in the hyperbolic
metric of $D_{256} \setminus K(f)$.  In particular, $h$ extends to a
homeomorphism conjugating $p|\overline D^4$ and $f|\overline B$, such that
$\hat h \circ h^{-1}$ is the identity on $K(f)$.  This homeomorphism
satisfies the required properties by construction.
\end{pf}
}

%Let $\hat T_1 \supset \hat T_2 \supset ...$ be the principal nest of $P$
%obtained in the standard way, see \cite {puzzle}.
We can use $h$ to transfer the principal nest $\hat T_1 \supset \hat T_2
\supset...$ for $P$ to the principal nest
$T_1 \supset T_2 \supset ...$ for $f$, $T_k=h(\hat T_k)$.  Let us say that
$F$ has $(N,\tau_0)$ {\it bounded
combinatorics} if the first return time of $T_{\tau_0+1}$ to $T_{\tau_0}$
is bounded by $N$.

\begin{lemma} \label {chi}

For every $N>0$, $\tau_0>0$, there exists $\kappa>0$, $\xi>0$, with the
following property.  If $F$ is $\kappa$-good and has $(N,\tau_0)$ bounded
combinatorics then
\begin{enumerate}
\item $|T_{\tau_0}|>\xi$ if $F$ is normalized,
\item The density of the domain of the
first return map to $T_{\tau_0}$ is at most $1-\xi$.
\end{enumerate}

\end{lemma}

\begin{pf}

Let us split the returning set of $T_{\tau_0}$ under the dynamics
of $F$ into two parts: $X$ consists of points that return under the dynamics
of $f$, and $Y$ is the complement.

Let us first show that $T_{\tau_0} \setminus X$ contains a disk of definite
radius.  Parameters $N$, $\tau_0$ specify a compact set of polynomials
with the property that the domain of the first return map to
$\hat T_{\tau_0}$ is not dense.  In particular, there
exists $\chi>0$ such that for any such polynomial,
the non-returning set of $\hat T_{\tau_0}$ contains a disk of radius
$\chi$.  Taking $\eta<\chi/10$ and $\kappa=\kappa(\eta)$ given by
Lemma~\ref {eta}, we conclude that $T_{\tau_0} \setminus X$ contains a
disk of radius $\chi/2$.  In particular this implies the first
statement, since $|T_{\tau_0}| \geq \chi^2/4$.

Since the area of $T_{\tau_0}$ is bounded by
$100$, the density of $X$ in $T_{\tau_0}$ is at most $1-\chi^2/400$.
If $x \in Y$, there exists $k>0$ minimal such that $f^k(x)$ belongs to some
non-central component $S$ of the domain of the first return map to $\U^0$. 
The first return map $\phi:S \to \U^0$ extends to a
univalent map $\phi:\tilde S \to \U$, where $\tilde S \subset \U^0 \setminus
\U^1$.  Let $S'$, $\tilde S'$ be the components of
$f^{-k}(S)$, $f^{-k}(\tilde S)$ containing $x$.  Notice that $f^k:\tilde S'
\to \tilde S$ is univalent, so $f^k:S' \to S$ has small distortion.
It follows that the density of $Y \cap S'$ in
$S'$ is at most $1-\sigma$ for some absolute $\sigma$.
The domains $S'(x)$ obtained in this way form a disjoint cover of $Y$ which
does not intersect $X$.  It follows that the density of $X \cup Y$ in
$T_{\tau_0}$ is at most $1-(\sigma \chi^2/400)$.
\end{pf}

Let $\tau(F)=\tau(P)$ be the height of $F$.

\begin{lemma}[see \cite {puzzle}] \label {rho}

There exists $\kappa>0$, $\rho>0$ such that if $F$ is
$\kappa$-good then $\mod(T_k \setminus \overline T_{k+1})>\rho k$,
$k \leq \tau(F)+2$.

\end{lemma}

For the proof of the next lemma, we will need some combinatorial preparation
for dealing with first return maps.
Let us denote the first return map to $T_k$ by $R_k:\cup T^j_k \to T_k$,
where $j$ takes values in $\Z$ and we set $T_{k+1}=T^0_k$.  Let $\Omega$ be
the set of finite sequences of non-zero integers $\d$.  The length of $\d$
will be denoted $|\d|$ (we allow $|\d|=0$).  If $\d=(j_1,...,j_m)$,
let $T^\d_k=\{x \in T_k,\, R^{i-1}_k(x) \in T^{j_i}_k, 1\leq i\leq m\}$, and
let $R^\d_k=R^m|T^\d_k$, so that $R^\d_k:T^\d_k \to T_k$ is univalent.
Let $W^\d_k=(R^\d_k)^{-1}(T_{k+1})$.  The map
$L_k:\cup W^\d_k \to T_{k+1}$, $L_k|W^\d_k=R^\d_k$
is the first landing map from $T_k$ to $T_{k+1}$, and its distortion
is small provided $\mod(T_k \setminus
T_{k+1})$ is big.

\begin{lemma} \label {Yk}

There exists $\kappa>0$, $\tau>0$, $\lambda_1<1$ such that for
$\tau \leq k \leq \tau(F)+2$ we have

\begin{enumerate}

\item Let $Y_k$ be the union of $T^j_k$ such that the distortion of
$R_k|T^j_k$ is at least $1+\lambda^k_1$.  Then $|Y_k| \leq \lambda^k_1
|T_k|$.

\item Let $Z_{k+1}$ be the union of connected components of
$(R_k|T_{k+1})^{-1}(\cup
T^j_k)$ where the distortion of $R^2_k$ is at least $1+\lambda^k_1$.  Then
$|Z_{k+1}| \leq \lambda^k_1 |T_{k+1}|$.

\end{enumerate}

\end{lemma}

\begin{pf}

Let $R_k|T^j_k=F^m$.  It is easy to see that $F^{m-1}|F(T^j_k)$ extends to a
univalent map onto $T_{k-1}$.  By the previous lemma, $F^{m-1}|F(T^j_k)$
has exponentially small distortion in $k$, $\tau \leq k \leq \tau(F)+2$.

By the previous lemma, if $\tau \leq k \leq \tau(F)$, the hyperbolic
diameter of any $T^j_k$ in $T_k$ is exponentially small.  Thus, the
distortion of $F|T^j_k$ is also exponentially small in $k$, provided
$T^j_k$ is not exponentially close to $0$ in the hyperbolic metric of $T_k$.
We conclude that $Y_k$ is contained in an exponentially small hyperbolic
neighborhood of $0$ in $T_k$.  The first assertion follows.

Let $D$ be a component of $(R_k|T_{k+1})^{-1}(\cup T^j_k)$.  By the previous
lemma, the hyperbolic diameter of $D$ in $T_{k+1}$ is exponentially small.
If the distortion of $R_k|D$ is not exponentially small, then the
distortion of either $R_k|(R_k(D))$ or of $R_k|D$ is
not exponentially small.  In the
first case, $D \subset (R_k|T_{k+1})^{-1}(Y_k)$.  In the second case, we
conclude as before that $D$ is contained in an exponentially small
hyperbolic neighborhood of $0$ in $T_{k+1}$.  Thus $Z_{k+1}$ is contained in
the union of $3$ sets with exponentially small hyperbolic diameter in
$T_{k+1}$.  The second assertion follows.
\end{pf}

\begin{lemma} \label {lambda}

There exists $\tau_0>0$ such that for every $N>0$, there exists $\kappa>0$,
$\lambda<1$ with the following property.  If $F$ is $\kappa$-good, has
$(N,\tau_0)$ bounded combinatorics, and $\tau_0 \leq k \leq \tau(F)+2$
then the density of the domain of the first return map to
$T_k$ is at most $\lambda^k$.

\end{lemma}

\begin{pf}

In the proof, $\lambda$'s will denote constants which are
definitely smaller than one,
depending only on $\rho$ from Lemma~\ref {rho}.
\comm{
Let us consider $\tau_0$ big and $\tau_0 \leq k \leq \tau(f)+2$.  By the
previous lemma, $\mod(T_{k-1} \setminus T_k) \geq \rho (k-1)$ is big.  It
follows that the distortion of $L_{k-1}$ is at most $1+\lambda^k_1$, and the
hyperbolic diameter of any $T^j_k$ in $T_k$ is at most $\lambda^k_1$.

Let $Y_k$ be the union of all $T^j_k$ where the distortion of $R_k$ is at
least $1+\lambda^k_2$ (in particular, $T^0_k \subset Y_k$).  Here we choose
$\lambda_2$ depending on $\lambda_1$, so that $Y_k$ is contained in an
exponentially small hyperbolic neighborhood of $0$ in $T_k$, hence $|Y_k|
\leq \lambda^k_3 |T_k|$.

Let $Z_{k+1}$ be the union of the components of
$(R_k|T_{k+1})^{-1}(\cup T^j_k)$ where the distortion of
$R_k^2|Z'$ is at least $1+\lambda^k_4$.  Here we choose $\lambda_4$
sufficiently close to $1$ so that $Z_{k+1}$ is the union of three sets
with exponentially small hyperbolic diameter, and we get
$|Z_{k+1}| \leq \lambda^k_5 |T_{k+1}|$.
}

Let $Y_k$, $Z_k$, and $\lambda_1$ be as in the previous lemma.
If $D \not\subset Z_{k+1}$ is a pullback of  some $T^j_k$ under  $  R_k | T_{k+1}$
then the distortion of the map $R_k^2: D\ra T_k$ is at most $1+\la_1^k$. Hence 
$$ 
    1- p(\cup T^j_{k+1} | D ) \geq (1+\lambda_1^k)^{-1}
(1- p (\cup W_k^\d| T_k)).
$$
It follows that 
\be \label {eq}
1-p(\cup T^j_{k+1}|T_{k+1}) \geq
(1-p(Z_{k+1}|T_{k+1})) (1+\lambda_1^k)^{-1}
(1-p(\cup W^\d_k|T_k)).
\ee
By Lemma \ref{Yk}, $1-p(\cup T^j_{k+1}|T_{k+1}) \geq
(1+\lambda^k_1)^{-1} (1-p(\cup T^j_k|T_k))$.
By Lemma \ref {chi}, there exists $\xi>0$ such that
$1-p(\cup T^j_{\tau_0}|T_{\tau_0})
\geq \xi$, so there exists $\sigma>0$ such that
\be \label {tilde chi}
1-p(\cup T^j_k|T_k) \geq \sigma, \quad \tau_0 \leq k \leq
\tau(F)+2.
\ee
Notice that for $j \neq 0$ the  map 
$ R_k:  T^j_k \ra T_k$ carries the set $  T^j_k  \cap \bigcup W^\d_k$  to $\bigcup W^\d_k$.
Hence, if $T^j_k \subset T_k \setminus Y_k$ then by the first part of Lemma \ref{Yk}, 
$p(\cup  W^\d_k\,  |\, T^j_k) \leq (1+\lambda^k_1)\,
p(\cup W^\d_k\, |\, T_k)$.  Hence
\be \label {Y_k}
   p( \cup W^\d_k\,  |\, T_k ) \leq p( Y_k\, |\,  T_k)  +(1-\sigma)
(1+\lambda^k_1)\, p( \cup W^\d_k\, |\, T_k ),
\ee
so if $k$ is big so that $(1-\sigma) (1+\lambda^k_1) \leq 1-(\sigma/2)$,
(\ref {tilde chi}) and (\ref {Y_k}) imply
\be
  p( \cup W^\d_k\, |\, T_k )  \leq \min   \{ 1-\sigma, 2 \sigma^{-1}
p(Y_k\, |\,   T_k)  \}.
\ee
Hence $|\cup W^\d_k| \leq \lambda^k_2 |T_k|$.  Putting this into
(\ref {eq}), we get $|\cup T^j_k| \leq \lambda^k_3 |T_k|$ as desired.
\comm{
For further use, let us notice that we may now write
\be
\frac {|W^\d_k|} {|T_k|} \leq |Y_k|+(1+\lambda^k_2)
\frac {|\cup T^j_k|} {|T_k|} \frac {|\cup W^\d_k|} {|T_k|},
\ee
that is
\be
(1-\lambda^k_8(1+\lambda^k_2)) \frac {|\cup W^\d_k|} {|W_k|} \leq |Y_k|.
\ee
}
\end{pf}

\comm{
\be
1-\hat \Delta_k \leq |Y_k|+\min \{1-|Y_k|,1-\tilde \chi\} \min \{1,

Let $\JJ_k$ be the set of $j$ such that the distortion of $R_k|U^j_k$ is
bigger than $\lambda^k_2$ (in particular $0 \in \JJ_k$).  Up to increasing
$\lambda_2$, we may assume that $\cup_{j \in \JJ_k} U^j_k$ is contained in a
$\lambda_3^k$ neighborhood of $0$.

Let $Y_k$ be the union of all $U^j_k$ where the distortion of $R_k$ is at
least $1+\lambda^k_4$ (in particular, $U^0_k \subset Y_k$).  Increasing
$\lambda_4$ if necessary, we may assume that the hyperbolic diameter of
$Y_k$ in $U_k$ is at most $\lambda^k_5$, so $|Y_k| \leq \lambda^k_5$.
Let $Z_{k+1}$ be the union of the components of
$(R_k|U_{k+1})^{-1}(\cup U^j_k)$ where the distortion of
$R_k^2|Z'$ is at least $1+\lambda^k_4$.  we
have $|Z_{k+1}| \leq \lambda_5^k |U_{k+1}|$.

Let
\be
\Delta_k=|U_k \setminus \cup U^j_k|,
\ee
so that $\Delta_{\tau_0} \geq \xi>0$ by the previous lemma.
Let
\be
\hat \Delta_k=|U_k \setminus \cup W^\d_k|,
\ee
so that $\hat \Delta_k \geq \Delta_k$.  Clearly
\be
\Delta_{k+1} \geq (1-|Z_{k+1}|)(1-\lambda_4^k) \hat \Delta_k \geq
(1-\lambda^k_6) \Delta_k,
\ee
so there exists $\sigma>0$ such that $\Delta_k \geq \sigma$ for
$\tau_0 \leq k \leq \tau(f)+2$.

Let $j^k$ be such that $R_k(0) \in U^{j_k}_k$.
Let $Z_{k+1}$ be the union of the components of
$(R_k|U_{k+1})^{-1}(\cup U^j_k)$ where the distortion of
$R_k^2|Z'$ is at least $1+\lambda_4^k$.  Up to incrasing $\lambda_4$, we
have $|Z_{k+1}| \leq \lambda_5^k$.

preimages of $U^j_k$
\be
Z_{k+1}=(R_k|U_{k+1})^{-1}(\cup_{j \in \JJ_k \cup \{j^k\}} U^j_k),
\ee
so that $|Z_{k+1}| \leq \lambda^k_4$.  Then
\be
\cup U^\d_{k+1} \subset Z_{k+1} \cup_{j \in \JJ_k \cup \{j^k\}}
(R_k|U_{k+1})^{-1}U^j_k

It follows that
\be
\Delta_{k+1} \geq (1-\lambda^k_4) |\Delta_k|,
\ee

Notice that
\be
\cup W^\d_k \subset \cup_{j \in \JJ_k} U^j_k \cup
\cup_{j \notin \JJ_k} (R_k|U^j_k)^{-1}(\cup W^\d_k),
\ee
so that
\be
|U_k|-|\cup W^\d_k| \geq \chi+(1-\lambda^k_4)|U_k|-|\cup W^\d_k|+
|\cup_{j \notin \JJ_k} W^\d_k(1+\lambda^k_4)|\cup

Since the
hyperbolic diameter of any $U^j_k$ is at most $\lambda^k_2$, we conclude
that for any $j
\in \JJ_k$
}

\begin{lemma} \label {4.19}

There exists $\tau_0>0$, such that for every $N>0$, there exists $\kappa>0$,
$\tau>\tau_0$, with the following property.  If $F$ is
$\kappa$-good, has $(N,\tau_0)$ bounded combinatorics and
$\tau(F) \geq \tau$ then the first return map to
$T_{\tau(F)}$ is $\kappa$-good.

\end{lemma}

\begin{pf}

Let $\kappa$, $\lambda$ be as in Lemma~\ref {lambda}, let
$\rho$ as in Lemma~\ref {rho},
and let $\tau$ satisfy $\lambda^\tau<\kappa^{-1}$, $\tau \rho \geq
\kappa$.
\end{pf}

\begin{thm}[Lean case exists]

There exists $\tau_0>0$, such that for every $N>0$, there exists
$\tau>\tau_0$ with the following property.  Let $f$ is an infinitely
renormalizable quadratic-like map such that all renormalizations of $f$ have
$(N,\tau_0)$ bounded combinatorics and height at least $\tau$.  Then
the Julia set of $f$ has Lebesgue measure zero.  Moreover, if $f$ is
a Feigenbaum map then $\HD(J(f))<2$.
%Moreover, there exists a
%sequence $r_i \to 0$ such that the density of the domain of the first return
%map to $\D_{r_i}$ is bounded away from $1$.

\end{thm}

\begin{pf}

Let $\tau_0$, N, $\kappa$, and $\tau$ be as in the previous lemma.
It is enough to prove the result for quadratic maps.  In this case, we can
consider a quadratic-like restriction $f=F:\U^0 \to \U$
which is $\kappa$-good.  By the previous lemma, there exists a domain
$\U_{(1)}=T_{\tau(F)}$ such that the first return map to $\U_{(1)}$ is
$\kappa$-good.  Repeating this procedure, we obtain a sequence of domains
$\U_{(n)}$ around $0$ such that the first return map $F_n$ to $\U_{(n)}$ is
$\kappa$-good.  We also let $\U_{(0)}=\U$.  Let $\U^0_{(n)}$ be the central
domain of the first return map to $\U_{(n)}$.  Let $f_n=F_n:\U^0_{(n)} \to
\U_{(n)}$.

Let $M_{(n)}^0$ be the domain bounded by the equator of $\U_{(n)} \setminus
\U_{(n)}^0$.  Let $M_{(n)}^k=f_n^{-k}(M_{(n)}^0)$.  Let $\V^n=M_{(n)}^{k_n}$
where $k_n$ is minimal with $\diam(M_{(n)}^{k_n})<100 \diam(J^n)$, and let
$\u^n=f_n^{-1}(\V^n)$.  Then $\u^n$ and $\V^n$ have the Properties (P1-4) of
Lemma~\ref {shapes} (see the proof of Lemma~\ref {eta}).
Let $V^n=f_n^{-k_n-2}(\U^0_{(n)})$,
$U^n=f_n^{-1}(V^n)$, and $A^n=V^n \setminus U^n$.
Then $U^n$ and $V^n$ satisfy the conditions (C1-3) and (G1-2) of \S~\ref
{nfd}.

We claim that the density of the domain of the first return map to $V^n$ is
bounded by some $\lambda<1$.  The argument is similar to the proof of
assertion (2) of Lemma~\ref {chi}.
If $x \in A^n$ is such that $F^k_n(x) \in V^n$ for some $k$,
then there exists $k \geq k_n+4$, and there exists $k_n+3 \leq r_x \leq k-1$
such that $F^l_n(x) \in \U_{(n)} \setminus \U_{(n)}^0$, $k_n+3 \leq l \leq
r_x$ and $F^{r_x+1}_n(x)$ belongs to the domain of the first return map
(under $F_n$) to $\U_{(n)}^0$.  Moreover, there exists
$\tilde S_x \ni x$ such that $F^{r_x+1}_n|\tilde S_x$ is a
univalent map onto $\U_{(n)}$.  Let
$S_x=(F^{r_x+1}_n|\tilde S_x)^{-1}(\U_{(n)}^0)$.  Then $S_x \subset V^n
\setminus U^n$ and the $S_x$ form a disjoint cover of
the domain of the first return map to $V^n \setminus
U^n$.  Since $F^{r_x+1}_n:S_x \to \U^0_{(n)}$ has
bounded distortion, one sees that there exists a definite (at least $1-C
\kappa^{-1}>0$) probability that a point in $A^n$ never returns to $V^n$.
Since $\area(A^n) \asymp \area(V^n)$, this concludes the claim.

Let $\Delta_n$ be the set of points in $\V^0$ that visit $V^n$.  Each branch
of the first landing map $\phi_n:\Delta_n \to V^n$ has a univalent extension
to $\V^n$ and thus has bounded distortion.  Notice that if $x \in
\Delta_{n+1}$ then $\phi_n(x)$ belongs to the domain of the first return map
to $V^n$.  This implies that $p(\Delta_{n+1}|\Delta_n)<\tilde \lambda<1$, so
$\area(\Delta_n)$ decreases exponentially fast.  By \cite {old}, almost
every $x \in J(f)$ is such that $0 \in \omega(x)$, so $\area(J)=\area(\cap
\Delta_n)=0$.

Assume now that $f$ is Feigenbaum.  In the notation of \S \ref {pert est},
we have shown
that $\xi_n \geq 1-C \kappa^{-1}$.  The same argument applied to $F_m$
shows that $\xi_{m,n} \geq 1-C \kappa^{-1}$.
By Lemma~\ref {exp}, $\eta_{m,n}$ is exponentially small in $n-m$.
Fix some $\epsilon>0$ small.  For $0<\rho<1$ small, let
\be
\omega^{[j]}_\rho(\delta)=\sup \omega^{[j]}_{m,n}(\delta), \quad
\xi_\rho=\sup \xi_{m,n}, \quad \eta_\rho=\sup \eta_{m,n},
\ee
where we take the supremums over pairs $m<n$ such that
$\rho^{1+\epsilon}<\rho_{m,n}<\rho^{1-\epsilon}$ (for every $m$, there
exists such an $n$ provided $\rho$ is sufficiently small).
Notice that such pairs satisfy $-C^{-1} \ln \rho \leq n-m \leq -C \ln \rho$
(where $C$ depends on $f$ through the geometric and combinatorial bounds),
so $\eta_\rho<\rho^\upsilon$ for some $\upsilon>0$.  Moreover, $\xi_\rho
\geq 1-C \kappa^{-1}$.
Putting together (\ref {301}) and (\ref {302}), we get that
for every $\delta$ close to $2$ (depending on $\rho$), we have
$\omega^{[j+1]}_\rho(\delta) \leq P_\rho(\omega^{[j]}_\rho(\delta))$, where
\be
P_\rho(x)=C \frac {\eta_\rho} {\rho^{1+\epsilon}}+
(1-C^{-1}(\eta_\rho+\xi_\rho)+C\eta_\rho \xi_\rho \rho^{-2\epsilon}) x+
C \xi_\rho \rho^{1-\epsilon} x^2.
\ee
If $\epsilon>0$ is small, then for every $\rho$ sufficiently small
$P_\rho$ has a positive fixed point $s$ that attracts $0$, and we have
$\lim_{j \to \infty} \omega^{[j]}_\rho(\delta) \leq s$, $\delta \approx 2$.
This implies that $\delta_\crit(f)<2$, so $\HD(J(f))<2$ by Theorem~C.
\end{pf}

\comm{
component, so $
By Lemma~\ref {old}, almost every $x \in J(f)$ satisfies $0 \in \omega(x)$. 
Thus we conclude that $p(J(f)|V^n) \leq \lambda<1$.  This implies that
$f$
In other words, in
the notation of \S, $\xi_n$ is bounded from below.
It follows that the density of the domain of the
first return map (under $f$).

Let us first show that $\area(J(f))=0$.  By \cite {old}, almost every $x \in
J(f)$ accumulates at $0$.  We will show that no such $x$ can be a density
point of $J(f)$.

hypothesis of \S. as
the central domain of the first return map to $M_{(k)}^j$.  As in Lemma~\ref

Clearly, $\mod(\U_k \setminus
\U_{k+1}) \geq \kappa$, so those domains shrink.  Let $\U^0_k$ be the
central component of the first return map to $\U_k$, and
$\U^1_k$ be the
central component of the first return map to $\U^0_k$.
Let $\Gamma^i_k=\{x \in \U^0,\, f^k(x) \in \U^i_k, \text { for some } k \geq
0\}$, $i=0,1$, and let $L^0_k:\Gamma^0_k \to \U^0_k$, be the
first landing map from $\U^0$ to $\U^0_k$.
All branches of $L^0_k$ are
univalent maps onto $\U^0_k$ that extend to univalent maps onto $\U_k$ and
hence have bounded distortion.
Since the domain of the first return map to $\U^0_k$ has density
$\kappa^{-1}$ in $\U^0_k$, we conclude that the density of $\Gamma^1_k$
in $\Gamma^0_k$ is bounded away from $1$.  Since $\Gamma^0_{k+1} \subset
\Gamma^1_k$, we conclude that $\cap \Gamma^0_k$ has Lebesgue measure
zero.  In particular, $\{x \in \U^0,\, 0 \in \omega(x)\}$ has Lebesgue
measure zero.  But it is well known, see \cite {old}
that $\{x \in J(f),\,
0 \notin \omega(x)\}$ has Lebesgue measure zero as well, so $|J(f)|=0$.
\end{pf}
}

\section{Various remarks and open problems}\label{remarks}

In this section we will formulate some related results, possible
further generalizations, and open problems. Some of these issues will
be elaborated in forthcoming notes.  

\subsection{Outstanding Problems}

The following  problems  still remain open:

\begin{problem} 
 Are there Feigenbaum Julia sets with positive area?
\end{problem}

We are inclined to believe that the answer is affirmative. 

\begin{problem}
Are there Feigenbaum Julia sets with zero area but Hausdorff dimension 2?
\end{problem}

According to Theorem~B, the affirmative answer to the first question
(together with certain {\it a priori} bounds)
implies the affirmative answer to the second.
%\marginpar{Not completely justified}

\begin{problem}\label{doubling}   
What is the area and dimension of the Julia set of the classical 
Feigenbaum quadratic polynomial $F_2$ corresponding to the doubling
renormalizations? 
\end{problem}

It is unlikely that Problem \ref{doubling} can be resolved
without computer assistance.

\ssk

All results of this paper extend naturally (with the same proofs) to the
setting of unicritical polynomial-like maps of some fixed degree $d \geq 2$
(specific constants may depend on the degree).
Let $F_d: z\mapsto z^d+c_d$ (where $d \geq 2$ is  an even integer) be the
doubling Feigenbaum map of degree $d$, and let $J_d=J(F_d)$.

\begin{problem}
  What is $\area(J_d)$ for big $d$?
\end{problem}
 
In case $\area(J_d)=0$,
it would  be interesting to find out whether  $\HD(J_d)=2$ for big $d$. 
(According to \cite{LS}, $\HD_\hyp(J_d)\to 2$ as $d\to\infty$.)

\subsection{Hausdorff dimension in the parameter space}

By our method one can produce a fair amount of Feigenbaum parameters
with small Julia set:
% (compare \cite{par-dim}): 

\begin{thm}
  The set of parameters $c\in M$ such that $P_c: z\mapsto z^2+c$ is a
Feigenbaum map with $\HD(J(f))<2$ has Hausdorff dimension at least 1.  
\end{thm} 

To see this, it is enough to use the method of \cite {par-dim} to obtain
large sets of Feigenbaum maps of high combinatorial type.
%\note{Proof?}

\subsection{Hyperbolic dimension and critical exponents}

Our methods of proving that $\HD_\hyp(J) = \HD(J)$ and
$\de_\crit = \de_*$  can be applied to much broader class of
rational maps: it will be elaborated in forthcoming notes. 
%What is needed to carry it through is that the
%dynamics is hyperbolic on the set of points that stay away from the
%critical points. 

\subsection{Absolutely continuous invariant measures}

Construction of invariant geometric measures  is probably the most important
traditional problem of ergodic theory. Here is an answer in our situation:  

\begin{thm}
   For any Feigenbaum map, there is a $\si$-finite measure
$\la$ equivalent to the $\de_\crit$-conformal measure $\mu$. 
\end{thm}

To construct this measure, one can use the method of \cite{KL}.
%First construct $\de_\crit$-conformal transverse and leafwise
%measures on the natural extension $\NN_f$ (the latter can be
%obtained just by lifting the conformal measure on $J$).  The product
%$\hat\la$ of these measures is an invariant measure on $\NN_f$ whose
%projection to $\C$ is $\la$.
%\marginpar{Is it possible to avoid laminations?}

\begin{problem} 
  Can the invariant geometric measure of a Feigenbaum map be finite? 
(Of course, it may happen only in the conservative case.)  
\end{problem}

\subsection{Hausdorff measure}

Conformal measures appear as a useful substitution of
Hausdorff measures $h_\de$ since in the non-hyperbolic dynamical situations the latter usually
vanish. The following statement (which easily follows from Lemma \ref{de and si} from Appendix B)
confirms this principle: 
%\note{how about packing measure?}

\begin{prop}
 Let $\de= \HD(J)$ for a Feigenbaum Julia set $J$ with periodic combinatorics. 
  If $\area(J)= 0$ then $h_\de(J)=0$.
\end{prop}

\subsection{Conservativity/dissipativity of towers}

The Trichotomy nicely translates to measure-theoretic properties
of towers $\bar f$
(see Appendix B): 

\begin{itemize}
\item {\it Lean case}: the Lebesgue measure is dissipative, supported on $\C\sm J(\bar f)$,
     and almost all orbits escape to infinity under the tower dynamics;

\item {\it Balanced case}: the Lebesgue measure is conservative under the tower dynamics,
               and supported on $\C\sm J(\bar f)$;

\item {\it Black hole case}:  the Lebesgue measure is dissipative, supported on $J(\bar f)$, 
     and almost all orbits are attracted to the postcritical set $\OO(\bar f)$.  
\end{itemize}  

\subsection{Real Fibonacci maps and wild attractors}

In this section, we will consider a class $\AAA$ of 
``real generalized polynomial-like maps''
 $$
     f:\bigcup_{j=0}^s I_j \to I, \quad s\geq 1,
$$ where
\begin{enumerate}
\item $I_j$ are disjoint closed intervals
      contained in a $0$-symmetric closed interval $I$, and moreover $I_j\subset \inter I$;
\item One of the $I_j$ (conventionally $I_0$) is a $0$-symmetric interval
     and $f:(I_0, \di I_0) \to (I, \di I)$ is a unimodal map which can be written as
     $\phi(|x|^l)$ where $\phi$ is a diffeomorphism with negative
     Schwarzian derivative onto a neighborhood of $I$, and $l>1$ is a real number 
     (the {\it criticality} of $f$);
\item If $j \neq 0$ then $f: (I_j, \di I_j)  \to (I, \di I)$ extends to  a diffeomorphism with negative
      Schwarzian derivative onto a neighborhood of $I$;
\item The critical point is recurrent and $\omega(0)$ is a minimal set;
\item All intervals $I_j$ intersect $\omega(0)$;
\item $f$ is not renormalizable in the classical sense: 
    there exist no $0$-symmetric intervals $T \subset I_0$ such that the first return map to $T$
    is a unimodal map.
\end{enumerate}
We shall consider such maps up to rescaling.  The {\it real Julia set} is defined as
$J^r(f)=\cap f^{-n}(I)$.  If it has positive length ($\equiv$ one-dimensional Lebesgue measure)
 then almost all orbits are attracted to $\omega(0)$, which is called a {\it wild attractor}.

The principal nest is an infinite sequence of intervals defined inductively
by $I^0=I$, $I^1=I_0$, and $I^{n+1}$ as the domain of the first return to
$I^n$ containing $0$.  Let $I_j^n$ be the connected components of the first
return map to $I^n$ which intersect $\omega(0)$ (conventionally, we let
$0 \in I_0^n$ so that $I_0^n=I^{n+1}$).  We let $f_n:\cup I_j^n \to I^n$ be
the first return map, and we call $f_{n+1}$ the {\it generalized renormalization}
of $f_n$.

We say that $f$ has a {\it bounded geometry} if all the intervals $I^j_n$
 and all the gaps in between are commensurable with $I^n$ 
(with constants independent of $j$ and $n$).  In this case
$0<\HD(\omega(0))<1$, and one is interested in $\HD(J^r(f))$.
%$\HD(J^r(f))>0$.

Maps of class $\AAA$ with bounded geometry  are quite similar to Feigenbaum maps, 
%(in fact, the combinatorics is sometimes simpler in the former case),
and our methods apply to  this setting as well.  There is one (temporary?)
advantage however: by \cite {BKNS}, we know that wild attractors do exist for some
maps $f\in \AAA$,
so that our results concerning Julia sets with positive measure become unconditional
in this context.
%Thus it is worth to translate a couple of results in this context.

Though we will concentrate on the real setting described above, 
maps of class $\AAA$ have natural complex counterparts: we refer the reader to \cite {LM}
and \cite {Bu}, as well as to a more recent work by Smania \cite {smania}.
Our methods apply equally well in this setting 
(where existence of Julia sets of positive Lebesgue measure is still an
open problem).

Maps of class $\AAA$ with $1<l \leq 2$ cannot have bounded geometry
(see \cite {attractors} for $l=2$ and \cite {shen} for a recent generalization).

We say that $f\in \AAA$ is a {\it  Fibonacci  map}  if for every
$n$ the domain of $f_n$ has two components, $I_0^n$ and $I_1^n$, and $f_n|I_1^n=f_{n-1}$ while
$f_n|I_0^n=f_{n-1}^2$.  For $l=2$,  geometry of Fibonacci maps was
analyzed in \cite {LM}.
A Fibonacci map with criticality $l>2$ has bounded geometry (see \cite {KN}).

For $l>2$ but close to $2$, 
the intervals $I_0^n$ and $I_1^n$ are much smaller than $I^n$  \cite {KN},
and in this case one can show that $J^r(f)$ is lean at the critical point.
% along the lines described here (the argument is actually simpler than for Feigenbaum
% maps, since the combinatorics is less complicate).  
Our methods then imply:

\begin{thm} \label {11.4}
There exists an $\eps>0$ such that if  $f$ 
is a Fibonacci map of degree  $l\in (2, 2+\epsilon)$, 
then $\HD(J^r(f))<1$.

\end{thm}

On the other hand, 
it is proved in \cite {BKNS} that if $f$ is a Fibonacci map of large degree $l$, 
then the real Julia set $J^r(f)$ has positive length (so that $f$ has a wild attractor). 
%Moreover, the estimates are uniform for the renormalizations of $f$, so
Then our methods imply:

\begin{thm} \label {11.5}

If $f$ is a Fibonacci map of sufficiently large degree $l$, 
then $\HD_\hyp(J^r(f)) < 1$.

\end{thm}

% \begin{problem}
% Analyze $\HD_\hyp(J^r(f))$ as $l \to \infty$.  Does it go to $1$?
%\end{problem} 

If $l$ is an even integer, one can prove more by methods of complex dynamics:

\begin{thm}[see \cite {Bu}]

If $l>2$ is an even integer, then  there exists a unique fixed point of the
generalized renormalization with Fibonacci combinatorics.  It attracts all
other Fibonacci  maps with the same criticality.

\end{thm}

%For fixed points of renormalization, the more complete description given by 
%the Tricothomy applies. 
It is generally believed that the same result is true with arbitrary criticalities:

\begin{conj} \label {fibo}

For every $l>2$ there exists a unique fixed point of the
generalized renormalization with the Fibonacci combinatorics.

\end{conj}

 Martens' method  \cite {martens} can probably be applied to prove existence,
 but the uniqueness part appears to be wide open at this point.

Let us denote by $f_l$ the fixed point of the Fibonacci generalized renormalization if it 
exists and is unique.
By compactness (real {\it a priori} bounds),
$f_l$ depends continuously on $l$ (over $l$'s for which it is defined). 
Thus if Conjecture~\ref {fibo} holds, we can see a ``phase transition'' in
the family $\{f_l\}_{l \in (2,\infty)}$.

\begin{thm}

Assume that Conjecture~\ref {fibo} holds.  Then all three cases of the
Trichotomy occur for $f_l$, $2<l<\infty$:
\begin{enumerate}
\item (Lean case)  For an open set of $l$'s containing $(2,2+\epsilon)$, 
                             $\HD_\hyp(J^r(f))=\HD(J^r(f))<1$;
\item (Balanced case) For a compact non-empty set of $l$'s, 
   $\HD_\hyp(J^r(f))=\HD(J^r(f))=1$ but the length  of $J^r(f)$ is zero;
\item (Black hole case) For an open set of $l$'s containing $(K,\infty)$,
     $\HD_\hyp(J^r(f))<1$ but $J^r(f)$ has positive length.
\end{enumerate}

\end{thm}

\begin{pf}
By Theorems~\ref {11.4} and~\ref {11.5}, 
the Lean and Black hole cases do occur for some sets of $l$'s
containing intervals $(2,2+\epsilon)$ and $(K,\infty)$
respectively.  To see that they occur for open sets of parameters, notice
that both are equivalent to a finite criteria
(existence of $n$ such that $\eta_n/\xi_n$ is either small or large,
with constants depending continuously on $l$, see Theorems~\ref {Shallow
case} and ~\ref {Black hole case}).  Since $(2,\infty)$ is
connected, the Balanced case must also occur.
\end{pf}

\begin{problem}
Is there a unique value  $l_0$ for which the Balanced case occurs?
% It is reasonable to expect that the Balanced case (being borderline)
%only happens at a very small set
% (of zero Lebesgue measure).  
%Is the function $l \mapsto \HD_\hyp(f_l)$ monotonically increasing up to $l=l_0$
%and decreasing after it?
\end{problem}
%\note {We used to ask whether $l \mapsto \HD_\hyp(f_l)$ is unimodal, but
%this seems incorrect after the work of Levin and Swiatek.
%Should we make reference to the work of Levin and Swiatek to discuss
%the behavior of $\HD_\hyp$?}

It seems unlikely that the Balanced case can ever happen for an integer $l$.%
\footnote {For the same reason
it is unlikely that there exists any example of a complex
Feigenbaum/Fibonacci periodic point of renormalization which falls into the
Balanced case, since there are only countably many parameters
(combinatorics and degree) to play with.}

\comm{
A generalized polynomial-like map will be called non-renormalizable if there
exists no $0$-symmetric interval $T$ such that the first return map to $T$
is unimodal.
If the critical point returns to $I^0$, then we define the generalized
renormalization of $f$ as the first return map to $I^0$.  It is also a
generalized polynomial-like map.  We will denote $f_n:\cup I^j_n \to I_n$
its generalized renormalizations.

A Fibonacci map is a non-renormalizable generalized polynomial-like map
with a recurrent critical such that:
\begin{enumerate}
\item 

A Fibonacci map is said to have bounded geometry if

In what follows, by a unimodal map we will understand an even map $f:I \to
I$ of the interval $I=[-1,1]$ which can be written as $f(x)=\phi(|x|^l)$,
where $\phi:I \to I$ is an orientation reversing diffeomorphism with
negative Schwarzian derivative, and $l>1$ is a real number
(the {\it criticality} of $f$).

The classical Fibonacci map is a unimodal map of the real interval
which is characterized by the property that the closest returns of the
critical point happen precisely at Fibonacci numbers $1,2,3,5,8,13,21,...$
(here we say that $s$ is a closest return of the critical point if
$|f^s(0)|<|f^k(0)|$ for $1 \leq k \leq s-1$).

We refer the reader to \cite {LMil}
}

% \subsection{Three dimensional laminations}

\appendix
\section{Nice fundamental domains} \label {mark}

Let $f:\u \to \V$ be a Feigenbaum map, and  let $f_n:\u^n \to \V^n$ be its
pre-renormalizations satisfying properties (P1)-(P4) of  Lemma~\ref {shapes}.
\comm{
We may assume that $\V^n \setminus \u^n$ does not
intersect the postcritical set of $f$ and that the modulus of
$\V^n \setminus \overline \u^n$ is bounded from below (such a lower bound
is called an unbranched {\it a priori} bound for $f$).
}
In this section we will construct a sequence of simply connected domains
$V^n$ satisfying properties (C1-3) and (G1-2) from \S \ref{nfd}.  
Our construction is reminiscent of Sinai's construction of Markov partitions \cite{Si}.

%The constants appearing in (G1-2) will only depend on the unbranched
%{\it a priori} bound of $f$.  Notice that such an {\it a priori} bound
%already imposes
%combinatorial restrictions, for instance, if some $f_n$ is immediately
%renormalizable then its renormalization period can be bounded in terms of
%the {\it a priori} bound.

\subsection{From admissible pairs to nice domains}
Take some $n\in \N$, and let $\FF^n$ be the family of inverse branches 
of $f^k$, $k \geq 0$, which are well defined  on $\V^n$. 

Let  $ X^n \subset \V^n$. 
For $g\in \FF^n$ (an inverse branch of $f^k$),
a set $g(X^n)\subset \V$ is called 
a {\it copy}  %, or a {\it univalent pullback},
of $ X^n$ of {\it depth} $k \geq 0$. 
In particular, let us consider the set $\CC^n= \{ g(0)\}_{g \in \FF^n}$
of all ``copies'' of $0$. 
There is a natural one-to-one correspondence between
$\CC^n$ and $\FF^n$ so that we can 
use $\CC^n$ to label all copies $X_x^n$  of $X^n$,  $x\in  \CC^n$.

%Let $\JJ^n=\cup_{k \geq 0} f^{-k}(J^n) \setminus \cup_{k \geq 0}
%f^{-k}(\beta_n)$.  If $x \in \JJ^n$, then there exists a minimal $k \geq 0$,
%called the {\it depth} of $x$ such that $f^k(x) \in J^n$.  Then there exists
%a domain $D \subset \V$ containing $x$ such that $f^k:D \to \V^n$ is a
%univalent map.  If $X \subset \V^n$ is any set, we shall denote
%$X_x=(f^k|D)^{-1}(X)$ which will be called a {\it copy}, or
%{\it univalent pullback} of $X$.

\begin{lemma} \label {copies1}

Let $X^n \subset \V^n$ be such that $f_n^{-1}(X^n) \subset X^n$.
Then $\bigcup_k f^{-k}(X^n)$ is the union of all copies of $X^n$.

\end{lemma}

\begin{pf}

Let $z\in \bigcup_k f^{-k}(X^n) $, and let $k$ be the first landing moment of $z$ to $X^n$.
Let $D_k\equiv \V^n$, and let us  define $D_{k-j}$ for $1 \leq j \leq k$ inductively as the connected
component of $f^{-1}(D_{k-j+1})$ containing $f^{k-j}(z)$.  Let us show that
$f^k:D_0 \to \V^n$ is univalent.  If this is not
the case, then   $0 \in D_{k-j}$ for some $1 \leq j \leq k$.
Take the  minimal such  $j$. 
Then  it is easy to see that $j$ is the renormalization period on level $n$  and $D_{k-j}=\u^n$. 
Hence $f^{k-j}(z) \in f_n^{-1}(X^n) \subset X^n$,
contradicting to the choice of  $k$.

It follows that $(f^k|D_0)^{-1}(X^n)$ is a copy of $X^n$ containing $z$.
\end{pf}

A pair $(Y^n,X^n)$ of nested sets $X^n\subset Y^n\subset \V^n$ is called 
{\it admissible} if every copy of $Y^n$ that intersects $X^n$ is contained
in $Y^n$.

\begin{lemma} \label {copies2}

Let $(Y^n,X^n)$ be an admissible pair.
\begin{enumerate}
\item Let $x,y \in \CC^n$ be such that $\depth(x) \leq \depth(y)$.  If
     $X^n_x \cap Y^n_y \neq \emptyset$ then $Y^n_y \subset Y^n_x$.
\item Assume $X^n$ is a domain. Let $W$ be a domain which is the
union of some copies of $X^n$.  If $X^n_x \subset W$ is a
copy of $X^n$ of smallest depth then $W \subset Y^n_x$.
\end{enumerate}

\end{lemma}

\begin{pf}

{\it First assertion}.  Let $k=\depth(x)$.  Then $f^k(X^n_x)=X^n$ intersects
$f^k(Y^n_y)=Y^n_{f^k(y)}$, so $f^k(Y^n_y) \subset Y^n$.  It follows that
$Y^n_y \subset Y^n_x$.
%If $X^n_x \cap Y^n_y \neq \emptyset$, but $Y^n_y$
%is not contained in $Y^n_y$ then $X^n \cap f^k(Y^n_y) \neq \emptyset$, but
%$f^k(Y^n_y) \not \subset Y^n$.  But $\depth(y) \geq k$ implies that
%$f^k(Y^n_y)=Y^n_{f^k(y)}$ is a copy of $Y^n$, which contradicts the
%hypothesis.

\ssk
{\it Second assertion}.  Notice that there is a sequence $x=x_1, x_2, \dots $ 
of points of $\CC^n$ such that $W$ is equal to the increasing union of
domains $X^n_{x_1}\cup \dots\cup X^n_{x_m}$.
Thus, we can restrict ourselves to the case where $W$ is itself a
finite union of copies of $X^n$.

Let us carry induction in the number of copies, $m$.  The assertion is obvious for $m=1$.  
For the induction step, let us consider connected components $U_l$ of $\bigcup_{i= 2}^m X^n_{x_i}$. 
By the induction hypothesis, $U_l\subset Y^n_{x_j}$ for some $j=j(l)\in [2,m]$.  
Since the union $U_l\cup X^n_x$ is connected, $U_l\cap X^n_x\not=\emptyset$;
all the more, $Y^n_{x_j}\cap X^n_x\not=\emptyset$.
By the first assertion, $Y^n_{x_j}\subset Y^n_x$. All the more, $U_l\subset Y^n_x$ and hence
$$
     W=X^n_x\cup \bigcup_l U_l \subset Y^n_x.
$$ 
 \end{pf}

Let us say that an open set $X \subset \V$ is  {\it nice} if $f^k(\partial X) \cap X=\emptyset$ for all $k \geq 0$
(compare \cite{Ma}).  Let us mention some nice properties of nice sets:   
%The nice property is very nice to manipulate:

\begin{lemma}

Let $f$ be a quadratic-like map.  Then

\begin{itemize}
\item[(N1)] For any open set $X \subset \V$, $\cup f^{-k}(X)$ is nice;

\item[(N2)] Components of a nice set are nice; 

\comm{
\item[N3] If $X$ is connected and nice and
$Y \subset \V$ is an open set such that
$\partial Y \subset \partial X$ then $Y$ is nice%(provided it is open)
\footnote {Here we use that the filled Julia set of $f$ has empty interior
(since $f$ is a Feigenbaum map).};
}

\item[(N3)] If $X$ is a nice domain then the filling of $X$ is nice,
provided $f$ has no attracting periodic orbits;
\comm{
\bignote{Misha: you had asked why the hypothesis of no bounded invariant
open sets play a role here.  First, this hypothesis makes the obvious proof
of this statement work.  Second, if say, $f$ has an attracting fixed point
$p$, then $X=\D_\epsilon(p) \setminus \overline {f(\D_\epsilon(p))}$ and
$Y=\D_\epsilon(p)$ give a counterexample.}
}

\item[(N4)] Intersection of two nice sets is  nice;

\item[(N5)] Preimages of nice sets are nice.
\end{itemize}

\end{lemma}

\begin{pf}

We will show (N3), as all the other properties are obvious.  Let $Z$ be the
filling of $X$.  Then $\partial Z \subset \partial X$.  If $z \in \partial
Z$ is such that $f^k(z) \in Z$ then $f^k(z) \in D$ where $D$ is a Jordan
disk with $\partial D \subset X$.  Since $X$ is nice and $\partial Z$ is
connected, this implies that $f^k(\partial Z) \subset D$, and hence
$f^k(Z) \subset D$ as well. 
Thus $f$ has an attracting periodic point in $D$, contradiction.
\end{pf}

Notice that property (C3) states precisely that $V^n$ should be nice.
The next lemma shows how to use certain admissible pairs to construct nice
domains.

\begin{lemma} \label {admissible}

Let $(Y^n, X^n)$ be an admissible pair, where $Y^n$ is a simply connected
domain and $X^n$ is a domain containing
$0$ satisfying $f_n^{-1}(X^n) \subset X^n$.  Then the filled connected component
$S^n$ of $\cup f^{-k}(X^n)$ containing $0$ is a nice simply connected domain
contained in $Y^n$, and such that $f_n^{-1}(S^n) \subset S^n$.

\end{lemma}

\begin{pf}

Let $Z^n$ be the connected component of $\cup f^{-k}(X^n)$ containing $0$.
$Z^n$ is  nice by (N1) and (N2). 
By Lemma~\ref {copies1}, $Z^n$ is the  union of some copies $X_{x_i}^n$ of $X^n$, $x_i\in \CC^n$.  
The smallest depth of these copies (that is, $0$) is reached on  $X^n\equiv X_0^n$ itself. 
By the second assertion of Lemma~\ref {copies2}, $Z^n\subset Y^n$.
By (N3), $S^n$ is nice as well.  Since $Y^n$ is simply connected, $S^n
\subset Y^n$, and since $X^n \subset S^n$, $f_n^{-1}(S^n) \cap S^n \neq
\emptyset$, the nice property implies that $f_n^{-1}(S^n) \subset S^n$.
\end{pf}

\comm{
\bignote{It is worthy to mention what properties of $(Y^n, X^n)$ are needed to ensure (G1) and (G2).} 
}

\subsection{Separation between some Julia copies}

Let $p_n$ be the renormalization period on level $n$, i.e.,  $f_n=f^{p_n}$.
Let $J^n=J(f_n)$, and let $J^n_i= f^i (J^n)$, $i=0,1,\dots, p_n-1$.
  Notice that the $J^n_i$ are copies of
$J^n$.  There are two cases to distinguish:
whether the $J^n_i$ are all disjoint (the {\it primitive} case) or not (the
{\it satellite} case).

To distinguish those cases, recall that $f_n$ has two fixed points,
conventionally denoted $\alpha_n$ and $\beta_n$, 
such that  $J^n \setminus \{\alpha_n\}$ is disconnected while $J^n \setminus \{\beta_n\}$ is connected.  
The satellite case happens if and only if $n>0$ and
$f_{n-1}$ is  {\it immediately renormalizable}, that is,
$\alpha_{n-1}=\beta_n$.  In this case, $J^n_i \cap J^n_j \neq
\emptyset$, $0 \leq i<j \leq p_n-1$ if and only if $j-i=0 \mod \frac {p_n}
{p_{n-1}}$, and in this case $J^n_i \cap J^n_j=f^i(\beta_n)$.

A $(\nu\diam J)$-neighborhood of a set $J$ will be called the
{\it $\nu$-enlargement} of $J$.  

\begin{lemma} \label {bouquets}

There exists $\nu>0$ such that if $J^n \cap J^n_i = \emptyset$ 
then the $ \nu $-enlargement of $J^n$ is disjoint from $J^n_i$.

\end{lemma}

\begin{pf}

Let $m \leq n$ be minimal
such that $J^n_i \subset J^m_k$ with $J^m_k \cap J^m =  \emptyset$. 
 Then $J^m_k \subset J^{m-1}$ (in the primitive case) 
or $J^m_k \subset J^{m-2}$ (in the satellite case). 
By a simple compactness argument, $J^m_k$ is disjoint from
the $\nu$-enlargement of $J^m$.  
All the more,  the smaller set $J^n_i \subset J^m_k$ is
disjoint from the $\nu$-enlargement  of  $J^n \subset J^m$.
\end{pf}

\subsection{Construction of admissible pairs}

Let $\nu$ be as in Lemma~\ref {bouquets}, and fix $\rho \leq \nu$.
\comm{
\bignote{Can we take $\rho=\nu$? Misha: Yes, but this extra freedom will be
useful later to obtain additional properties of this construction (namely,
expansion of the first return map).}
}
Let $l \geq 0$ be minimal such that $\Gamma^n\equiv f_n^{-l}(\u^n)$
is contained in the intersection of the $\rho$-enlargement of $J^n$ with
$\u^{n-1}$.  Notice that:
\begin{enumerate}
\item $\mod(\Gamma^n \setminus f_n^{-1}(\overline \Gamma^n))$ is definite
           (once $\rho$ is fixed),
\item If $J^n_i \cap \Gamma^n \neq \emptyset$ then either
        $J^n_i=J^n$ or $f_{n-1}$ is immediately renormalizable and
        $J^n_i=f_{n-1}^j(J^n)$ for some $j$,
\item If a copy $\Gamma^n_x$ of $\Gamma^n$ of depth $k>0$ intersects $J^n$
then $f_{n-1}$ is immediately renormalizable and $f^k|\Gamma^n_x$ is an
iterate of $f_{n-1}$.

Indeed, $\Gamma^n_x \cap J^n \neq \emptyset$ implies that $f^k(\Gamma^n_x)
\cap f^k(J^n)=\Gamma^n \cap f^k(J^n) \neq \emptyset$.  By Lemma~\ref
{bouquets}, $J^n$ touches $f^k(J^n)$, so $f_{n-1}$ is immediately
renormalizable, and $k=s p_{n-1}$ for some integer $s>0$.
To conclude, we must show that $f^{j p_{n-1}}(\Gamma^n_x)
\subset \u^{n-1}$, $0 \leq j<s$.  But if this is not the case then
the connected set $f^{j p_{n-1}}(\Gamma^n_x)$ intersects
both $J^{n-1} \supset f^{j p_{n-1}}(J^n)$ and $\partial
\u^{n-1}$.  This implies that further iterates $f^{r p_{n-1}}(\Gamma^n_x)$,
$r>j$, have the same property, and in particular $\Gamma^n=f^k(\Gamma^n_x)$
intersects $\partial \u^{n-1}$, contradiction.

\comm{
$f^k(J^n)$ intersects $\Gamma^n$ and (by Lemma~\ref {bouquets})
touches $J^n$, and hence $J^n_x=(f^k|\Gamma^n_x)^{-1}(J^n)$ touches $J^n$.
}
\comm{
\bignote{why? Misha: is it better now?}
}
\end{enumerate}

\comm{
\bignote{Should we also require that $\Gamma^n\cap \OO(f) \subset  J^n$? 
Misha:  This is quite immediate, since $\V^n \cap \OO(f)$ is already
contained in $J^n$ (unbranched apriori bounds).  I don't think it needs to
be stressed at this point.}
}

We will consider separately the primitive and satellite cases.

\subsubsection{Primitive case}

Assume that $f_{n-1}$ is not immediately renormalizable.  In this case, any
copy of $\Gamma^n$ of depth $k>0$ is disjoint from $J^n$.

Let $Y^n=f_n^{-1}(\Gamma^n)$.  Notice that $\partial Y^n$ has bounded
hyperbolic diameter in the hyperbolic metric of $\V^n \setminus J^n$.

If $Y^n_x$ is a copy of $Y^n$ of positive depth,
then the hyperbolic diameter of $Y^n_x$ in $\V \setminus J^n \supset
\Gamma^n_x$ is bounded 
(since $\mod(\Gamma^n_x \setminus \overline Y^n_x)$ is definite). 
 Notice that the hyperbolic metrics of
$\V \setminus J^n$ and  $\V^n \setminus J^n$ are comparable on $\u^n \supset Y^n$. 
Hence if $Y_x^n \cap \di Y^n\not=\emptyset$, then the hyperbolic diameter of
$Y^n_x \cap Y^n$ 
in the hyperbolic metric of $\V^n \setminus J^n$ is bounded by some $K$.

Let $s \geq 0$ be minimal such that the hyperbolic distance
 (in the hyperbolic metric of $\V^n \setminus J^n$)
from any $z \in X^n\equiv f_n^{-s}(Y^n)$ to $\partial Y^n$
is greater than   $K$.  Then $(Y^n,X^n)$ is admissible.

\ssk
 Notice also that the 
hyperbolic diameter of $\overline Y^n \setminus X^n$ 
(in the hyperbolic metric of $\V^n \setminus J^n$) is uniformly bounded.

\subsubsection{Satellite case}

Assume that $f_{n-1}$ is immediately renormalizable.

\begin{lemma} \label {boundedness}

There exists a simply connected domain $Y^n \subset \Gamma^n$ such that
\begin{enumerate}
\item $Y^n$ is nice for $f_{n-1}$ and $f_n^{-1}(Y^n) \subset Y^n$;
\item There exists a finite set $F$ of copies of $\beta_n$ such that $Y^n$
       is a neighborhood of $J^n \setminus F$;
\item $\mod(\Gamma^n \setminus \overline Y^n)$ is definite;
%\item $Y^n$ contains a definite neighborhood of $0$ in scale $\diam(J^n)$,
\item The hyperbolic diameter of $\partial Y^n$ in $\V^n \setminus \OO(f)$
        is uniformly bounded;
\item If $Y^n_x \cap \partial Y^n \neq \emptyset$ then
       $\Gamma^n_x \subset \Gamma^n \setminus J^n$;
\item If $Y^n_x \cap \partial Y^n \neq \emptyset$ then the diameter of
$Y^n_x$ in the hyperbolic metric of $\V^n \setminus J^n$ is bounded.
\comm{
If $Y^n_x$ is a copy of $Y^n$ of positive depth such that $Y^n_x\cap \di Y^n\not=\emptyset$;
% and if $z \in Y^n_x \cap Y^n$ then $z$ is at a bounded 
then the Hausdorff distance from $Y^n_x$ to  $\partial Y^n$  in the hyperbolic metric of $\V^n \setminus J^n$
is bounded.
}
\end{enumerate}

\end{lemma}

\begin{pf}

The properties described in Lemma~\ref {shapes} allow us to obtain a
straightening of $f_{n-1}:\u^{n-1} \to \V^{n-1}$ with bounded dilatation,
and thus to define equipotentials and external rays in $\V^{n-1}$.
Let $E^{n-1}$ be the domain bounded by the equipotential of level $2$. 
Let us cut $E^{n-1}$ by the external rays landing at
$\alpha_{n-1}$, and let $Z^n$ be the connected component of $0$. 
Let us show that $Y^n=f_n^{-l}(Z^n)$ has all the desired properties if $l$ is large enough
(independent of $n$). 

Properties (1) and (2) are obvious.
If $l$ is large enough then $Y^n \subset \Gamma^n$.
Moreover, then properties (3) and (4) are also valid for every fixed $l$.

\ssk
Let us deal with property (5).
Let us show first that if  $Y^n_x \cap\di Y^n\not= \emptyset$ then $\Gamma^n_x \cap J^n = \emptyset$. 
Indeed,  if $\Gamma^n_x \cap J^n \neq \emptyset$ then
$f^k:J^n_x \to J^n$, $k=\depth(x)$, is an iterate of $f_{n-1}$.
But since  $Y^n$ is nice for $f_{n-1}$,
$Y^n_x \cap \partial Y^n=\emptyset$ - contradiction.

Let us now show that if $Y^n_x \cap \partial Y^n\not=\emptyset$ 
then $\Gamma^n_x \cap \partial \Gamma^n = \emptyset$.
Notice that $\diam(Y^n_x) \asymp \diam(\Gamma^n_x)$,
$\mod(\Gamma_x^n \setminus \overline Y_x^n)=\mod(\Gamma^n \setminus \overline Y^n)$ 
is bounded from below (uniformly for all sufficiently big $l$),
 and $\partial Y^n$ is contained in a $\epsilon$-neighborhood of $J^n$ 
(where $\epsilon \to 0$ when $l \to \infty$). 
Those geometric properties together show that if $\Gamma^n_x$ intersects
both boundaries,
$\di \Gamma^n$ and $\partial Y^n$ (so that $Y^n_x$ has size comparable
with $J^n$), then $Y^n_x \cap \partial Y^n \neq \emptyset$ implies that
$\Gamma^n_x \cap J^n \neq \emptyset$ for $l$ big enough, contradiction.
This proves (5). 

\ssk
The last property follows from the third, the fifth, and the fact that the hyperbolic
metrics of $\V \setminus J^n$ and of $\V^n \setminus J^n$ are comparable on $\u^n$.
\end{pf}

Let $X^n$ be the connected
component of $0$ of the interior of
$$
Y^n \setminus \bigcup_{Y^n_x \cap \partial Y^n \neq \emptyset} Y^n_x.
$$

\begin{lemma}

The pair $(Y^n,X^n)$ is admissible.  Moreover:
\begin{enumerate}
\item $X^n$ is simply connected,
\item $X^n$ is a neighborhood of $J^n \setminus F$,
\item $f^{-1}_n(X^n) \subset X^n$,
\item $\overline Y^n \setminus X^n$ has bounded hyperbolic diameter in $\V^n
\setminus \OO(f)$.
\end{enumerate}

\end{lemma}

\begin{pf}

If $Y^n_x \cap X^n \neq \emptyset$ then $Y^n_x \cap \partial Y^n \neq
\emptyset$ by definition of $X^n$, and hence $Y^n_x \subset Y^n$.  Let us
now check Properties (1-4).

(1).  Notice that $X^n$ is a bounded connected
component of the complement of the closure of the connected set
$$
\bigcup_{\partial Y^n_x \cap \partial Y^n \neq \emptyset} \partial Y^n_x,
$$
so it is simply connected.

%\bignote{how?}
(2).  By Properties (2) and (6) of Lemma~\ref {boundedness}, 
$X^n$ is a neighborhood of $J^n \setminus F$.

(3).  If $Y^n_x \cap \partial Y^n \neq \emptyset$ then by Property (5) of
Lemma~\ref {boundedness}, $\Gamma^n_x \subset \Gamma^n \setminus J^n$ which
implies that $\depth(x) \geq p_n$ and in particular
$f_n(Y^n_x)$ is a copy of $Y^n$.
Thus, if $Y^n_x \cap f_n^{-1}(X^n) \neq \emptyset$ and $Y^n_x
\cap \partial Y^n \neq \emptyset$, then there is a copy of $Y^n$,
$f_n(Y^n_x)$, intersecting $X^n$ and $\partial (f_n(Y^n))$,
contradicting the definition of
$X^n$.  Since $J^n \setminus F \subset X^n$ (by (2)),
$f_n^{-1}(X^n)$ is a connected set intersecting $X^n$, and
by the definition of $X^n$ it follows that
$f_n^{-1}(X^n) \subset X^n$.

%\note{I don't understand logic here} 
(4).  By Property (6), $\overline Y^n \setminus X^n$ has bounded
hyperbolic diameter in $\V^n \setminus \OO(f)$.
\end{pf}

%\bignote{We should made more clear what these remarks are for.}

\subsection{Construction of domains $U^n$ and $V^n$}

Let $(Y^n,X^n)$, $n \geq 0$ be given by the construction of the previous
section.  Let $S^n$ be given by Lemma~\ref {admissible}.  We can not set
$V^n=S^n$ because this does not guarantee property (C2)
(although it would possess the other properties).

This can be fixed as follows.  Let $Q^n$ be the connected component of
$0$ of $\cap_{m \leq n} S^m$.
By (N2) and (N4), $Q^n$ is nice.  Let
$V^n=f_n^{-1}(Q^n)$, and $U^n=f_n^{-1}(V^n)$.  Then $V^n$ is nice by
(N5), so $V^n$ has Property (C3).
Notice that $J^n \setminus F \subset V^n$, where $F$ is a finite set
of copies of $\beta_n$ (empty if $f_{n-1}$ is not immediately
renormalizable).  Since we also have $V^n \subset \u^n$, this implies
(C1).  To check (C2), assume there is $x \in V^n \cap \partial
U^{n-1}$.  Then $f_n(x) \in Q^n \subset Q^{n-1}$ and $f_{n-1}^2(x) \in
\partial Q^{n-1}$.  Notice that $f_n(x)=f_{n-1}^k(x)$ where $k=p_n/p_{n-1}
\geq 2$.  So $f_{n-1}^{k-2}(\partial Q^{n-1}) \cap Q^{n-1} \neq \emptyset$,
contradicting that $Q^{n-1}$ is nice.

We have so far checked the combinatorial properties (C1-3), and we will now
verify the geometric properties (G1-2).

\comm{
\bignote{Artur, It is not clear what you are after in the section.
It looks like everything is already done and then you start over again
(verifying various random properties of different domains).
It think it should be organized as follows: \\
1) In A1, formulate geometric properties of admissible pairs needed for the desired
   geometric properties of $V^n, U^n$;\\
2) After the construction of admissible pairs in A2, check that these properties are satisfied.
QED
}
}

We first make a couple of remarks on the geometry of $S^n$.
We have that $\partial S^n$ has
bounded hyperbolic diameter in $\V^n \setminus \OO(f)$.  Since $\V^n \cap
\OO(f) \subset S^n$, this implies that the Euclidean distance from $\partial
S^n$ to $\OO(f) \cup \partial \V^n$ is comparable with $\diam(J^n) \asymp
\diam(\OO(f) \cap \V^n) \asymp \diam(\V^n)$.  In
particular, $S^n$ contains a neighborhood of $0$ of radius comparable with
$\diam(J^n)$.  On the other hand, $S^n \subset \V^n$
is contained in a disk of radius $\diam(\V^n) \asymp \diam(J^n)$.  Since
$\diam(J^n)$ decreases exponentially fast, we see that there exists
$\upsilon \geq 0$ such that $S^n \subset S^m$ whenever $n \geq
m+\upsilon$.
\comm{
Moreover, $S^n$ is a neighborhood of $J^n \setminus F$ where $F$ is a finite
set consisting of copies of $\beta_n$ (empty if $f_{n-1}$ is not immediately
renormalizable).
}

\comm{
\begin{lemma}

Let $r$ and $R$ be respectively  minimal and  maximal radii such that
$\D_r \subset S^n \subset \D_R$.  Then
$R \asymp r \asymp J^n$.

\end{lemma}

\begin{pf}

Notice that $S^n \subset \V^n$, so $\diam(J^n) \asymp \diam(\V^n)$
implies $R \asymp \diam(J^n)$.

Since the hyperbolic diameter of $\partial S^n$ in $\V^n \setminus \OO(f)$
is uniformly bounded, $S^n \supset \OO(f) \cap \V^n$
 and $\diam(\V^n) \asymp \diam(\OO(f) \cap \V^n)$, it follows that
$S^n$ contains a neighborhood of $\OO(f) \ni 0$ of radius comparable with
$\diam(\V^n)$.
\end{pf}

\begin{lemma} \label {upsi}

There exists $\upsilon_0 \geq 0$ such that if $\upsilon \geq \upsilon_0$
then $S^{n+\upsilon} \subset S^n$ for all $n$.

\end{lemma}

\begin{pf}

This follows from the fact that $\diam(J^{n+k}) \leq C \lambda^k
\diam(J^n)$, for some $C>0$, $\lambda<1$.
\end{pf}
}

%Let $Q^n$ be the connected component of $0$ in $\cap_{m \leq n} S^m$.

\begin{lemma} \label {shee}

The hyperbolic diameter of $\partial Q^n$ in $\V^n \setminus \OO(f)$ is
bounded.

\end{lemma}

\begin{pf}

%By Lemma~\ref {upsi}, if $x \in \partial Q^n$ then $x \in \partial S^m$
%for some $m \leq n$ and $n-m<\upsilon_0$.
%Each $\partial S^m$ has bounded
%hyperbolic diameter in $\V^m \setminus \OO(f)$.
For $m<n$, the Euclidean distance
of points in $\partial S^m$ to $\OO(f)$ is comparable with $\diam
J^m \geq \diam(J^n)$.
This implies that the hyperbolic distance in $\V^n \setminus \OO(f)$
of points in $\partial S^m \cap S^n$ to $\partial S^n$ is bounded.
Since $\partial Q^n \subset \partial S^n \cup \bigcup_{m<n} (S^n \cap
\partial S^m)$, the result follows.
%Since for $m \leq n$, the
%hyperbolic metrics of $\V^n \setminus \OO(f)$ and of $\V^m \setminus \OO(f)$
%are comparable on $S^n \subset \u^n \cap \u^m$, the result follows.
\end{pf}

\comm{
\begin{lemma} \label {C1-3}

Let $V^n=f_n^{-1}(Q^n)$.  Then
\begin{enumerate}
\item $V^n$ is a nice simply connected domain,
\item $V^n$ is a neighborhood of $J^n \setminus F$, where $F$ is a finite
set of copies of $\beta_n$ (empty if $f_{n-1}$ is not immediately
renormalizable),
\item $V^n \subset f_{n-1}^{-1}(V^{n-1})$.
\end{enumerate}

\end{lemma}

\begin{pf}

The first two properties are obvious.  For the third,
let $x \in V^n \cap \partial f_{n-1}^{-1}(V^{n-1})$.  Then $f_n(x)
\in Q^n \subset Q_{n-1}$ and $x \in \partial f_{n-1}^{-2}(Q^{n-1})$.  In
particular, $f_{n-1}^2(x) \in \partial Q^{n-1}$.  Notice that
$f_n(x)=f_{n-1}^k(x)$ where $k=p_n/p_{n-1} 
\geq 2$.  Thus $f_{n-1}^{k-2}(\partial Q^{n-1}) \cap Q^{n-1} \neq \emptyset$,
contradicting that $Q^{n-1}$ is nice.
\end{pf}

Let $U^n=f_n^{-1}(V^n)$.
}

\begin{lemma} \label {dc}

The hyperbolic diameter of $\overline V^n \setminus U^n$
in $\V^n \setminus \OO(f)$ is uniformly bounded.

\end{lemma}

\begin{pf}

Recall that
if $g:X \to Y$ is a covering map of degree $n$ between hyperbolic
Riemann surfaces and $Z \subset Y$ is
such that $g^{-1}(Z)$ is connected then the hyperbolic diameter of
$g^{-1}(Z)$ in $X$ is at most $n$ times the hyperbolic diameter $d$
of $Z$ in $Y$.
%Indeed, if $x \in Z$ then, by the Schwarz Lemma,
%$g^{-1}(Z)$ is contained in the neighborhood of $g^{-1}(x)$ of radius $d$,
%and the result follows since $g^{-1}(Z)$ is connected.

Let us show first that $\partial V^n$ and $\partial U^n$ have bounded
hyperbolic diameter in $\V^n \setminus \OO(f)$.
By the Schwarz Lemma, the hyperbolic diameter of
$\partial V^n$ in $\V^n \setminus \OO(f)$ is
bounded by its hyperbolic diameter in $f_n^{-1}(\V^n \setminus \OO(f))$.
Since $\partial V^n=f_n^{-1}(\partial Q^n)$ is
connected, its hyperbolic diameter in
$f_n^{-1}(\V^n \setminus \OO(f))$ is bounded by two times the hyperbolic
diameter of $\partial Q^n$ in $\V^n \setminus \OO(f)$, which is bounded by
Lemma~\ref {shee}.  Applying this argument to $\partial
U^n=f_n^{-1}(\partial V^n)$, we see that $\partial U^n$ also has bounded
hyperbolic diameter in $\V^n \setminus \OO(f)$.

This implies that the
Euclidean distance of every $x \in \overline V^n \setminus U^n$ to $\partial
\V^n \cup \OO(f)$ is comparable with $\diam J^n$.  So the hyperbolic
distance of every $x \in \overline V^n \setminus U^n$ to $\partial V^n \cup
\partial U^n$ in $\V^n \setminus \OO(f)$ is bounded.  The result follows.
\comm{
From now on, we only consider the hyperbolic metric in $\V^n \setminus
\OO(f)$.  Let us show that any $x \in \partial U^n \setminus J^n$
is at bounded
hyperbolic distance from $\partial V^n$, which clearly implies the result.
Let $k \geq 1$ be such that
$f^k_n(x) \in \overline \u^n \setminus f_n^{-1}(\overline \u^n)$.
Since $\overline \u^n \setminus f_n^{-2}(\overline \u^n)$ has bounded
hyperbolic diameter, we see that $f_n^k(x)$
and $f_n^{k-1}(x)$ are at bounded hyperbolic distance.  By the Schwarz
Lemma, $x$ and is at bounded hyperbolicinside $\V^n \setminus \OO(f)$
with bounded hyperbolic length.
It follows from the Schwarz Lemma (applied
$k-1$ times) that the
hyperbolic distance between $x$ and $f_n(x)$
is bounded by the maximum of the hyperbolic distances in
$\V^n \setminus \OO(f)$ between $f^k_n(x)$ and its two preimages by
$f_n$, and so it is uniformly bounded.
}
\end{pf}

Property (G1) is given by Lemma~\ref {dc}.  Let us show that (C1-3) and (G1)
imply (G2).

Since $\partial U^n$ has bounded hyperbolic diameter in $\V^n
\setminus \OO(f)$, $U^n \supset \V^n \cap \OO(f)$, and $\diam(\V^n) \asymp
\diam(\V^n \cap \OO(f))$, it follows that $U^n$ contains a disk $D$
centered around $0$ of radius comparable with $\diam(V^n)$.  Thus
$\area(U^n) \asymp (\diam(U^n))^2 \asymp (\diam(V^n))^2$.  We must
show that $\area(V^n \setminus U^n) \asymp \area(V^n)$.

\begin{lemma} \label {W and k}

There exists $k>0$ (independent of $n$) and a round disk
$W^n \subset V^n \setminus f^{-k}_n(\u^n)$ with
$\diam(W^n) \asymp \diam(J^n)$.

\end{lemma}

\begin{pf}

This follows easily by a compactness argument, but we will show how it can
be deduced from the construction, since this gives more explicit estimates.

In the primitive case, the construction implies that $V^n \supset X^n
\supset f_n^{-l}(\u^n)$ for some $l$.  The annulus $f_n^{-l-1}(\V^n
\setminus \overline \u^n)$ has a definite modulus, and encircles $J^n$,
so it contains a round disk of diameter comparable with $\diam(J^n)$.  We
can take $W^n$ as this disk and $k=l+1$.

In the satellite case, notice that for every $\epsilon>0$, there exists
$r>0$ such that $f_n^{-r}(\u^n) \setminus Y^n$ can be tiled into a bounded
number $N$ (independent of $\epsilon$) of sets of Euclidean diameter
less than $\epsilon \diam(J^n)$ (contained in $\epsilon \diam
(J^n)$-neighborhoods of points where $\partial Y^n$ touches $J^n$).
Notice that $\partial V^n$ is contained in a bounded
neighborhood of $\partial Y^n$ in the hyperbolic metric of
$\V^n \setminus J^n$.  This implies that we can choose $l$ such that
$f^{-l}(\u^n) \setminus V^n$ is also the
union of $N$ sets with Euclidean diameter less than $\diam(J^n)/(100 N)$.
Since the annulus $f_n^{-l-1}(\V^n \setminus \overline \u^n)$ has a definite
modulus and encircles $J^n$,
$f_n^{-l-1}(\V^n \setminus \overline \u^n) \cap
V^n$ contains a round disk of diameter comparable with $\diam(J^n)$.  We can
take $W^n$ as this disk and $k=l+1$.
\end{pf}

%Notice that there exists $k>0$ (independent of $n$) and
%a disk $W \subset U^n$ of radius comparable with
%$\diam(V^n)$ such that $W \cap f_n^{-k}(\u^n)$.  This follows either by a
%compactness argument, or directly from the construction.
%\note {Detail the construction}

Taking $W$ and $k$ as in Lemma~\ref {W and k}, we have
$W \subset \cup_{j \leq k} f_n^{-j}(V^n \setminus U^n)$,
so $\area(V^n) \asymp \area(W) \leq C \area(V^n \setminus U^n)$.
This concludes the proof of (G2).

\begin{rem} \label {small domains}

Notice that this construction gives $V^n$ with additional
properties:
\begin{enumerate}
\item [(A1)] $U^n$ contains a round disk around $0$ of radius comparable with
$\diam(V^n)$,
\item [(A2)] $V^n$ is a neighborhood of $J^n \setminus F$ where $F$ is a
finite (bounded) set of copies of $\beta_n$.  Moreover, if $f_{n-1}$ is
not immediately renormalizable then $F=\emptyset$.
\end{enumerate}
Moreover, in the construction we have freedom to choose $\rho$ small so
that $\Gamma^n$, and hence
$V^n \subset Y^n \subset \Gamma^n$, is contained in the
$\rho$-enlargement of $J^n$.  This implies that if $J^n_x$ is a copy
of $J^n$ of depth $k>0$ which
intersects $V^n \subset Y^n$ then $\diam(J^n_x)$ is much smaller than
$\diam(J^n)$ (this is easy in the primitive case, and follows from
bounded shape of $\Gamma^n_x$ and
assertion (5) in Lemma~\ref {boundedness} in the satellite case).
Under these circumstances, we obtain:
\begin{enumerate}
\item [(A3)] If $V^n_x$ is a univalent pullback of $V^n$ of depth $k>0$
contained in $V^n$ then $\diam(V^n_x)$ is much smaller than $\diam(V^n)$
(since $\diam(V^n_x) \asymp \diam(J^n_x)$).
\item [(A4)] If $V^n_x$ is a univalent pullback of $V^n$ of depth $k>0$
contained in $V^n$ then $|Df^k(x)|>2$ for every $x \in V^n_x$.  This follows
from the Koebe Distortion Lemma and (A3).
\end{enumerate}

\end{rem}

\begin{rem} \label {permarkimpro}

If $f$ is a periodic
point of renormalization, $f_p(x)=\lambda f(\lambda^{-1} x)$, we may assume
that $\V^{n+p}=\lambda \V^n$.  In this case we can modify the above
construction to obtain $V^n$ satisfying $V^{n+p}=\lambda V^n$ as follows.

We can construct $Y^n$ satisfying $Y^{n+p}=\lambda Y^n$.  In the primitive
case, we can easily get that $X^{n+p}=\lambda X^n$, but in the
immediately renormalizable case
we will only get $X^{n+p} \subset \lambda X^n$.
%(although we can easily get equality if $f_{n-1}$ is not
%immediately renormalizable).
To fix this, in the immediately renormalizable case we modify
the construction by taking the connected component of $0$ in the interior
of $\cap_k \lambda^{-k} X^{n+kp}$ as the new definition
of $X^n$.  Then $X^{n+p}=\lambda X^n$, but this only gives
$S^{n+p} \supset \lambda S^n$.  We take
$\cup_k \lambda^{-k} S^{n+kp}$ as the new definition of $S^n$.  Then
$S^{n+p}=\lambda S^n$, but we only have $Q^{n+p}=\lambda Q^n$ for large
$n$.  We take $Q^n=\lambda^{-r} Q^{n+rp}$ for $r$ large.  Then
$V^n=f_n^{-1}(Q^n)$ satisfies $V^{n+p}=\lambda V^n$.  All properties (C1-3),
(G1-2) and (A1-4) can be easily verified as before.

\end{rem}

\begin{rem}

The only parameters determining the constants in this construction are the
unbranched {\it a priori} bounds and the combinatorics.  Actually, one can
see that a lower bound on the unbranched {\it a priori} bounds implies an
upper bound on the period of any immediate renormalizations.  As a
consequence, the combinatorics are only used to get a lower bound on $\nu$
(the spacing between non-touching $J^n_i$).
\comm{
As such, they are uniform for
all real quadratic maps (or real unicritical polynomials of a fixed degree),
as well as for certain classes of quadratic maps of
``high combinatorial type'' considered in \cite {puzzle}.
}

\end{rem}

\comm{
We may take $W^n(\rho)$ satisfying
$W^{n+p}(\rho)=\lambda W^n(\rho)$.  We then define $Q^n$ as above.
The sequence $\lambda^{-k}Q^{n+kp}$ is then
decreasing and we obtain $V^n$ as the filled
$\inter (\cap \lambda^{-k}Q^{n+kp})$.  Properties (C1-3) and (G1-2) are
easily checked as before.

\end{rem}

It is
a nice simply connected domain, which is a neighborhood of $J^n \setminus
F$, where $F$ is a finite number of copies of $\beta_n$.  This implies also
that $f_n^{-1}(Q^n) \subset Q^n$.  By
Lemma~\ref {upsi}, $\partial Q^n$ is contained in the union of $\partial
S^{n-\upsilon}$,...,$\partial S^n$.  Let $V^n=f_n^{-1}(Q^n)$.

\begin{lemma}

We have $V^n \subset f_{n-1}^{-1}(V^{n-1})$.

\end{lemma}

\begin{pf}

Let $x \in V^n \cap \partial f_{n-1}^{-1}(V^{n-1})$.  Then $f_n(x)
\in Q^n$ and $x \in \partial f_{n-1}^{-2}(Q^{n-1})$.  In
particular, $f_{n-1}^2(x) \partial Q^{n-1}$.  Notice that
$f_n(x)=f_{n-1}^k(x)$ where $k=p_n/p_{n-1} 
geq 2$.  Thus $f_{n-1}^{k-2}(\partial Q^{n-1}) \cap Q^{n-1} \neq \emptyset$,
contradicting that $Q^{n-1}$ is nice.
\end{pf}

\begin{lemma}

The hyperbolic diameter of $

\comm{
\begin{lemma} \label {dc}

The hyperbolic diameter of $\partial f_n^{-1}(S^n)$ in $\V^n \setminus
\OO(f)$ is uniformly bounded.

\end{lemma}

\begin{pf}

Recall that
if $g:X \to Y$ is a covering map with $n$ sheets between hyperbolic
Riemann surfaces and $Z \subset Y$ is
such that $g^{-1}(Z)$ is connected then the hyperbolic diameter of
$g^{-1}(Z)$ in $X$ is at most $n$ times the hyperbolic diameter $d$
of $K$ in $Y$.  Indeed, if $x \in K$ then, by the Schwarz Lemma,
$g^{-1}(Z)$ is contained in the neighborhood of $g^{-1}(x)$ of radius $d$,
and the result follows since $g^{-1}(Z)$ is connected.

By the Schwarz Lemma, the hyperbolic diameter of $f_n^{-1}(\partial S^n)$ in
$\V^n$ is bounded by its hyperbolic diameter in $f_n^{-1}(\V^n \setminus
\OO(f))$.
Since $f_n^{-1}(\partial S^n)$ is connected, its hyperbolic diameter in
$f_n^{-1}(\V^n \setminus \OO(f))$ is bounded by twice the hyperbolic
diameter of $\partial S^n$ in $\V^n \setminus \OO(f)$.
\end{pf}

\begin{lemma} \label {shee}

Let $r$, respectively $R$, be minimal, respectively maximal, such that
$\D_r \subset f_n^{-1}(S^n)$ and $S^n \subset \D_R$.  Then
$R \asymp r \asymp J^n$.

\begin{lemma}

\begin{pf}

Notice that $S^n \subset \V^n$, so $\diam(J^n) \asymp \diam(\V^n)$
implies $R \asymp \diam(J^n)$.

By Lemma~\ref {dc}, the hyperbolic diameter of
$f_n^{-1}(\partial S^n)$ in $\V^n
\setminus \OO(f)$ is also uniformly bounded.  Since $S^n \supset \OO(f) \cap
\VV^n$ and $\diam(\V^n) \asymp \diam(\OO(f) \cap \V^n)$, it follows that
$S^n$ contains a neighborhood of $0$ of radius comparable with
$\diam(\V^n)$.
\end{pf}
}

\comm{
\begin{lemma} \label {Snk}

There exists $k>0$ such that $S^{n+k} \subset f_n^{-1}(S^n)$.

\end{lemma}

\begin{pf}

It follows from Lemma~\ref {shee} and the fact that $\diam(J^{n+k}) \leq
C \lambda^k \diam(J^n)$ for some $C>0$, $\lambda<1$.
\end{pf}
}

Let $V^0=S^0$ and define by induction $V^n$ as the connected component of
$0$ in $f_n^{-1}(V^$.  By Lemma~\ref {Snk},

Since $f_n^{-1}(S^n)$ contains $\OO(f) \cap
\V^n$, $\partial S^n$ has bounded hyperbolic diameter in
$\V^n \setminus \OO(f)$, and $\par
(by the Schwarz Lemma, it is at most twice the diameter of \partial S^n$),
it follows that $r \asymp \diam(J^n)$.
\end{pf}

\begin{lemma}

There exists $k>0$ such that $S^{n+k} \subset f_n^{-1}(S^n)$.

\end{lemma}

The pair $(X^n,Y^n)$ is admissible.  Moreover, $\overline Y^n \setminus X^n$
has bounded hyperbolic diameter in $

\end{lemma}

\begin{pf}

It is clear that $M^n$ is simply connected.  It is easy to see that if
$Y^n_x$ intersects $J^n$ then $T^n_x=T^n$, so $J^n \setminus \{\beta_n\}
\subset M^n$.  To see that $f^{-1}_n(M^n) \subset M^n$, it is enough to show
that $f_n(\partial M^n) \cap M^n=\emptyset$.  Notice that if $y \in \partial
M^n$ then either $y \in \partial T^n$ or $y \in \partial T^n_x$, where
$T^n_x \cap \partial T^n \neq \emptyset$ (this follows from the fact that
the diameters of the $T^n_x$ go to zero).  In the first case, $f_n(y) \notin
T^n \subset M^n$.  In the second case, either $f_n(T^n_x) \cap T^n=\emptyset$,
in which case $y \notin T^n \supset M^n$ again, or $f_n(T^n_x)=T^n_{f_n(x)}
\cap \partial T^n \neq \emptyset$, in which case $y \notin M^n$ by
construction.

By taking iterates, we may assume that $\depth(T^n_x)=0$, and $M^n_x=M^n$.
If $M^n \cap T^n_y \neq \emptyset$, then by construction $T^n_y \cap
\partial T^n \neq \emptyset$, so $T^n_y \subset T^n$.
\end{pf}

\begin{lemma} \label {Mn}

With respect to the hyperbolic metric in $\V^n \setminus J^n$,
$\partial M^n$ and $\partial T^n$ are at bounded distance in the Hausdorff
sense.  In particular, $M^n$ contains a definite neighborhood of $0$ in
scale $\diam(J^n)$.

\end{lemma}

\begin{pf}

Since $\partial T^n$ has uniformly bounded diameter in $\V^n \setminus
\OO(f)$, it is enough to show that if $T^n_x \cap \partial T^n \neq
\emptyset$ then $T^n_x$ has bounded hyperbolic diameter in $\V^n \setminus
\OO(f)$.  By the Schwarz Lemma, the hyperbolic diameter of $T^n_x$
inside $\Gamma^n \setminus J^n$ is smaller than the hyperbolic
diameter of $T^n_x$ in $\Gamma^n_x$ which is bounded.

The last statement follows from the behavior of the hyperbolic metric near
$0$: at every $x \in \V^n$, it is at least of order of the inverse of the
Euclidean distance to $0$.  This follows from $J^n$ connected and
$\diam(J^n) \asymp \diam(\V^n)$.
\end{pf}

Let $S^n$ be the connected component of $0$ in $\cup f^{-k}(X)$.

\begin{lemma}

$S^n \subset T^n$.

\end{lemma}

\begin{pf}

$S^n$ is an increasing union of a finite number
$M^n_{x_1}$,...,$M^n_{x_l}$ of copies of $M^n$.  One proves by
induction on $l$ that such finite union is contained in $T^n_{x_i}$ such
that $\depth(M^n_{x_i})$ is minimal.
\end{pf}

\begin{lemma}

The hyperbolic diameter of $\partial S^n$ in $\V^n \setminus \OO(f)$ is
uniformly bounded.

\end{lemma}

\begin{pf}

We have $\partial S^n \subset \overline T^n \setminus M^n$ which is a
compact subset of $\V^n \setminus \OO(f)$.  The result now
follows from Lemma~\ref {Mn}, assertion (3) of Lemma~\ref {boundedness},
and hyperbolic geometry.
\end{pf}

\begin{lemma}

In the primitive case, there exists an admissible pair $(X^n,Y^n)$ such that
$X

There exists $k \geq 0$ (independent of $n$) and $0 \leq s \leq t \leq k$
such that $(f^{-s}_n(\u^n),f^{-t}_n(\u^n))$ is admissible.

\end{lemma}

\begin{pf}

Let $\nu$ be as in Lemma~\ref {bouquets} and take $s$ minimal such that if
$Y^n=f_n^{-s}(\u^n)$ then every copy $Y^n_x$ is contained in a
$\nu$-neighborhood of $J^n_x$.  It follows that $s$ is uniformly bounded,
and thus $\mod(Y^n \setminus J^n)$ is bounded from below.  By the choice of
$s$, no copy of $Y^n$ intersects $J^n$ except for $Y^n$ itself.

Consider a copy $Y^n_x$ that interects $\partial Y^n$.

In particular, each connected component of a finite 
Since the connected component of $\cup f^{-k}(X^n)$ containing $X^n$ is the
union of all 

Let us say that $X \subset \V^n$ is admissible if it is a simply connected
domain containing $0$ such that $f_n^{-1}(X) \subset X$.

If $f^k(x)$ is the first landing moment of some
$x \in U$ to $J^n$ then there exists a domain $D \subset V$ containing $x$
such that $f^k|D$ is a univalent map onto $V^n$.  We call the sets
$J^n_x=(f^k|D)^{-1}(J^n)$  a {\it copy} or {\it univalent pullbacks}
of $J^n$ of {\it depth} $k$.

Let us first consider small (but definite) fundamental annuli of the maps
$f_n$:

\begin{lemma} \label {nrho}

For every $\rho>0$, there exists $\tau>0$ and a domain $W^n(\rho) \subset
U^n$ containing $J^n$, such that
\begin{itemize}

\item   $f_n^{-1}(W^n(\rho)) \Subset W^n(\rho)$;
\item  $\mod(W^n(\rho) \setminus \overline {f_n^{-1}(W^n(\rho))}) \geq \tau$;
\end{itemize}

and for any univalent map $g:D \to \V^n$,
\begin{itemize}

\item  $g^{-1}(W^n(\rho))$ is contained in a $\rho$-neighborhood of $g^{-1}(J^n)$, and
\item the map $g: g^{-1}(W^n(\rho)) \ra W^n(\rho)$ has a bounded distortion.

\end{itemize}
\end{lemma}

\begin{pf}
  Let $W^n (\rho)= f_n^{-k} \V^n $ with sufficiently big $k$ 
(depending on the unbranched {\it a priori} bound).
\note{Maybe better to first regularize by taking the equator of
$\V^n \setminus \u^n$.} 
\end{pf} 

Let $J^n=J(f_n)$, and 
let $J^n_i= f^i (J^n)$, $i=0,1,\dots, p-1$, where $p$ is the period of $J^n$.

  If $f^k(x)$ is the first landing moment of some
$x \in U$ to $J^n$ then there exists a domain $D \subset V$ containing $x$
such that $f^k|D$ is a univalent map onto $V^n$.  We call the sets
$J^n_x=(f^k|D)^{-1}(J^n)$  a {\it copy} or {\it univalent pullbacks}  of $J^n$ of {\it depth} $k$.

Let $W^n_x(\rho)$ be the component of $f^{-k}(W^n(\rho))$ containing $J^n_x$ 
(in case of $J^n_i$, we will also use notation $W^n_i(\rho)$). 
Then
\begin{itemize}
 \item   $f^k: W^n_x(\rho) \to W^n(\rho)$ extends to a univalent map onto $\V^n$;
 \item    $f^k: W^n_x(\rho) \to W^n(\rho)$ has a bounded distortion;
 \item $W^n_x(\rho)$ is contained in the $\rho$-neighborhood of $J^n_x$.
\end{itemize}

\subsection{Primitive case}

We will first assume that none of the $f_n$ is immediately renormalizable.
In this case, the $\nu$-neighborhoods of the $J^n_i$ are disjoint for some
small, but definite, $\nu$ (depending on the unbranched {\it a priori}
bounds and possibly on the compinatorics).\note{I do not see how
to argue that there are not any combinatorial dependencies.}
Thus the neighborhoods $W^n_i(\rho)$ are pairwise disjoint
if $\rho$ is  sufficiently small (less than $\nu$).

The following estimate will play the key role in our construction:

\begin{lemma} \label {copies}

For every $\epsilon>0$, there exists $\delta>0$ such that if the
$\delta$-neighborhoods of two copies $J^n_x$ and $J^n_y$ of $J^n$ intersect
then $\diam(J^n_y)<\epsilon \diam(J^n_x) $, provided
$\depth (J^n_y) \geq \depth (J^n_x)$.

\end{lemma}

\begin{pf} 
Let $k\equiv \depth (J^n_x)\leq \depth(J^n_y)\equiv l$.
Consider the neighborhood $W^n_y(\rho)$ of $J^n_y$ defined above. 
Since the annulus $W^n_y \sm J^n_y$ has a definite modulus,
there is a $\de=\de(\rho)>0$ such that $W^n_y $ contains the  $\de(1+\eps^{-1})$-neighborhood of $J^n_y$.
 
If the $\delta$-neighborhoods of $J^n_x$ and $J^n_y$ of $J^n$ intersect but
$  \diam(J^n_y) \geq  \eps \diam(J^n_x) $,  
then $W^n_y \cap J^n(x)\not=\emptyset$. Applying  $f^l$ to $W^n_y \cap J^n(x)$ we see that 
$W^n (\rho) \cap J^n_{l-k}\not=\emptyset$, where $J^n_{l-k}\not= J^n$, 
 contradicting the choice of $\rho$.
\end{pf}

Let $f^k(x)$ be the first landing of some point $x \in U$ to
$W^n(\rho)$.  Let $W^n_x(\rho)$ be the component of
$f^{-k}(W^n(\rho))$ containing $x$. 
It follows that $f^k:W^n_x(\rho) \to W^n(\rho)$ extends to a univalent map
on $\V^n$.  Notice that $W^n_x(\rho)$ is contained
in a $\rho$-neighborhood of $J^n_x$.

\begin{lemma}\label{kappa}

For every $\kappa>0$, there exists $\rho>0$ such that if
$x_1,...,x_n \in \cup f^{-k}(J^n)$, $\cup W^n_{x_i}(\rho)$
is connected then $\cup W^n_{x_i}(\rho)$ is contained in the
$\kappa$-neighborhood of $J^n_{x_j}$, where $\depth(J^n_{x_j)})$ is
minimal.

\end{lemma}

\begin{pf}

Let $\epsilon \mapsto \delta(\epsilon)$ be the function given in Lemma~\ref
{copies}.
We may assume that $\kappa<\delta(1/10)$.  Set
$\rho=\min\{\kappa/10,\delta(\kappa)\}$.
We shall prove the result by induction on $n$.  It is obvious for $n=1$.
For the induction step, assume that $W^n_{x_1}(\rho)$ is maximal.  Let
$W$ be a connected component of $\cup_{i \geq 2} W^n_{x_i}(\rho)$.  Then $W$
intersects $W^n_{x_1}(\rho)$ and $W$ is contained in the
$\kappa$-neighborhood of some $J^n_{x_j}$ with $j \geq 2$.  In
particular, the $\kappa$-neighborhoods of $J^n_{x_j}$ and $J^n_{x_1}$
intersect, so $\diam(J^n_{x_j}) \leq \diam(J^n_{x_1})/10$.  Thus $W$ is
contained in the $\kappa$-neighborhood of $J^n_{x_1}$ as required.
\end{pf}

\begin{lemma}

For every $\kappa>0$ there exists $\rho>0$ such that if
$S$ is a connected component of $\cup f^{-k} W^n(\rho)$ containing 0
 then  $S$ is contained in a $\kappa$-neighborhood of $J^n$. 

\end{lemma}

\begin{pf}
Notice first  that there is a sequence $x_i$ such that 
$$
  S = \bigcup_{i=1}^\infty W^n_{x_i}(\rho),
$$ 
where the finite unions $S_m = \bigcup_{i=1}^m W^n_{x_i}(\rho)$ are connected. 

By Lemma \ref{kappa}, all  $S_m$ are contained in the $\kappa$-neighborhood of $J^n$,
and we are done.
\end{pf}

Given a domain $U \subset \C$, the {\it filled} $U$ is defined as
the smallest simply connected set containing $U$.  It is automatically open
and is constructed by adding to $U$ all compact components of $\C\sm U$.

Choose now $\kappa$ small, $\rho$ as in the previous lemma and let $S^n$ be
the connected component of $0$ in $\cup f^{-k} W^n(\rho)$.  We define a
sequence $Q^n$ inductively by setting $Q^0=S^0$, $Q^{n+1}=S^{n+1} \cap
f_n^{-1}(Q^n)$.  We then define $V^n$ as the filled $Q^n$.

Let us now check conditions (C1-3) and (G1-2).
\begin{enumerate}
\item [(C1)] It is enough to notice that $V^n \supset W^n(\rho)$.
\item [(C2)] Since $Q^{n+1} \subset f_n^{-1}(Q^n)$, $V^{n+1}$
is contained in the filled $f_n^{-1}(Q^n)$, which is contained in
$f_n^{-1}(V^n)$.
\item [(C3)] The boundary of $S^n$ never returns to $S^n$ (since
$S^n$ is a component of a backward invariant open set), so it follows by
induction that the boundary of $Q^n$ never returns to $Q^n$.  Since
$\partial V^n \subset \partial Q^n$, if $f^k(\partial V^n) \cap V^n \neq
\emptyset$ for some $k$ then $f^k(V^n) \subset V^n$, which contradicts $J(f)
\cap V^n \neq \emptyset$.
\item [(G1)] Notice that the annuli $\V^n \setminus
\overline \u^n$ and $f_n^{-1}(W^n(\rho)) \setminus J^n$ separate $A^n$
from the postcritical set and have definite modulus.
\item [(G2)]  Since $J^n \subset U^n \subset V^n$ which is contained in a
$1+\kappa$-neighborhood of $J^n$, we have $\diam(U^n) \asymp \diam(V^n)$.
Since $U^n \setminus J^n$ has definite modulus, $|U^n| \geq
C^{-1} \diam(J^n)^2$.  This implies that $\diam(U^n)^2 \asymp |U^n|$.  To
conclude, we have to show that $|A^n|/|V^n|$ is bounded from below.  Notice
that $f^{-k}(V^n) \subset W^n(\rho)$ for some $k$ depending on $\kappa$
and $\tau$.  Thus $V^n \setminus f^{-(k+1)}(V^n)$ has a definite modulus. 
It follows that $|V^n| \asymp |V^n \setminus f^{-(k+1)}(V^n)|=\sum_{j=0}^k
|f^{-j}(A^n)|$.  Thus  $|f^{-j}(A^n)| \asymp |V^n|$ for some $0
\leq j \leq k$, which implies that $|A^n| \asymp |V^n|$.  (An alternative
argument is to use compactness.)
\end{enumerate}

\begin{rem} \label {small domains}

Notice that this construction gives $V^n$ with additional
properties:
\begin{enumerate}
\item [(A1)] Any univalent pullback of $V^n$ contained in $V^n$ has
diameter much smaller than the diameter of $V^n$.
\item [(A2)] By bounded distortion, any univalent pullback as in (A1)
is a contraction by a factor at least $1/2$,
\item [(A3)] $U^n$ contains a round disk around $0$ of radius comparable with
$\diam(V^n)$,
\item [(A4)] $V^n$ is a neighborhood of $J^n$ ({\it this is specific to the
primitive case}).
\end{enumerate}

\end{rem}

\subsection{Modification for periodic points of renormalization}

If $f$ is a periodic
point of renormalization, $f_p(x)=\lambda f(\lambda^{-1} x)$, we may assume
that $\V^{n+p}=\lambda \V^n$.  In this case we can modify the above
construction to obtain $V^n$ satisfying $V^{n+p}=\lambda V^n$ as follows.

We may take $W^n(\rho)$ satisfying
$W^{n+p}(\rho)=\lambda W^n(\rho)$.  We then define $Q^n$ as above.
The sequence $\lambda^{-k}Q^{n+kp}$ is then
decreasing and we obtain $V^n$ as the filled
$\inter (\cap \lambda^{-k}Q^{n+kp})$.  Properties (C1-3) and (G1-2) are
easily checked as before.

\subsubsection{Immediately renormalizable case}

We now assume that some renormalizations of $f$ are immediately
renormalizable.
For each $n$ for which $f_{n-1}$ is not immediately renormalizable, we
construct domains $S^n$ as before.
For each $n$ for which $f_{n-1}$ is immediately renormalizable we will
now show how to construct domains $S^n$.  The domains $S^n$ can then
be modified to give rise to the desired $V^n$ as before.

Let $p$ and $q$ be such that $f_n=f^p=f_{n-1}^q$.
The key problem in adapting the previous argument to the immediately
renormalizable case is that Lemma~\ref {copies} fails.  Indeed, distinct
$J^n_i$ may touch, but in a controlled way.  Two copies $J^n_i$ and $J^n_j$
intersect if and only if $j-i=0 \mod p/q$.  Moreover, there exists small,
but definite, $\nu>0$ such that the $\nu$-neighborhoods of $J^n_i$, $J^n_j$
are disjoint if $j-i \neq 0 \mod p$.

\begin{lemma}

There exists a simply connected domain $\Gamma^n \subset \u^n \cap \V^{n-1}$
such that
\begin{enumerate}
\item $J^n \subset \Gamma^n$ and $\mod(\Gamma^n \setminus J^n)$ is definite,
\item If $\Gamma^n \cap J^n_i \neq \emptyset$ then there exists $j$ such
that $f_{n-1}^j(J^n)=J^n_i$.
\end{enumerate}

\end{lemma}

\begin{pf}

Take $\Gamma^n=f_n^{-l}(\u^n)$ for suitable $l$.
\end{pf}

Let us call a domain $X \subset \V^n$ admissible if
\begin{enumerate}
\item $X$ is simply connected,
\item $f_n^{-1}(X) \subset X$,
\item $J^n \setminus X$ is a finite set contained in the preorbit of
$\beta_n$..
\end{enumerate}
If $X$ is admissible, and $f^k(x) \in J^n \setminus \{\beta_n\}$ is the
first landing of $x$ to $X$ then there exists a domain $x \in D$ such that
$f^k:D \to \V^n$ is univalent.  We let $X_x=(f^k|D)^{-1}(X)$ which we will
call a copy of $X$ of depth $k$.  If $X$ is admissible then
$\cup f^{-k}(X)$ is the union of all copies of $X$.  Notice that the depth
of $X_x$ depends only on $x$ and not on $X$.

\begin{lemma} \label {boundedness}

There exists an admissible set $T^n \subset \Gamma^n$
such that
\begin{enumerate}
\item $\mod(\Gamma^n \setminus \overline T^n)$ is definite,
\item $T^n$ contains a definite neighborhood of $0$ in scale $\diam(J^n)$,
\item The hyperbolic diameter of $\partial T^n$ in $\V^n \setminus \OO(f)$
is uniformly bounded,
\item If $T^n_x \cap \partial T^n$ then $\Gamma^n_x \subset \Gamma^n
\setminus J^n$.
\end{enumerate} 

\end{lemma}

\begin{pf}

The properties described in Lemma~\ref {shapes} allow us to obtains a
straightening of $f_{n-1}:\u^{n-1} \to \V^{n-1}$ with bounded dilatation,
and thus to define equipotentials and external rays.
Let $E^{n-1}$ be the domain bounded by the equipotential at height
$2$.  Let us cut $E^{n-1}$ by the external rays landing at
$\alpha_{n-1}$, and let $T^n=T^{n,l}$ be the connected component of $0$.  It
is clear that $T^n$ is admissible for every $l$ sufficiently large
(independent of $n$), and that it satisfies the first assertion.
The second and third assertions are also clear for every large, but fixed,
$l$.

Let us deal with the fourth statement.  The only $J^n_i$
which intersect $\Gamma^n$ are of the form $f_{n-1}^j(J^n)$, so if
$\Gamma^n_x \cap J^n \neq \emptyset$ then $J^n_x$ touches $J^n$.  This
implies that $f^k:J^n_x \to J^n$, $k=\depth(J^n_x)$,
is an iterate of $f_{n-1}$.  But by construction of $T^n$ (which is nice for
$f_{n-1}$), this implies that $T^n_x \cap \partial T^n=\emptyset$.

Since $\Gamma^n_x \setminus T^n_x$ has definite moduli and $\partial T^n_x$
is very close to $J^n$, it follows that the diameter of $J^n_x$ is much
smaller than the diameter of $J^n$ and in particular $\Gamma^n_x \subset
\Gamma^n$.
\end{pf}

Let $M^n$ be the connected
component of $0$ of the interior of
$$
T^n \setminus \bigcup_{T^n_x \cap \partial T^n \neq \emptyset} T^n_x.
$$

\begin{lemma}

The domain $M^n$ is admissible.  Moreover, if a $M^n_x \cap T^n_y \neq
\emptyset$ and $\depth(T^n_x) \leq \depth(T^n_y)$ then $T^n_y \subset
T^n_x$.

\end{lemma}

\begin{pf}

It is clear that $M^n$ is simply connected.  It is easy to see that if
$T^n_x$ intersects $J^n$ then $T^n_x=T^n$, so $J^n \setminus \{\beta_n\}
\subset M^n$.  To see that $f^{-1}_n(M^n) \subset M^n$, it is enough to show
that $f_n(\partial M^n) \cap M^n=\emptyset$.  Notice that if $y \in \partial
M^n$ then either $y \in \partial T^n$ or $y \in \partial T^n_x$, where
$T^n_x \cap \partial T^n \neq \emptyset$ (this follows from the fact that
the diameters of the $T^n_x$ go to zero).  In the first case, $f_n(y) \notin
T^n \subset M^n$.  In the second case, either $f_n(T^n_x) \cap T^n=\emptyset$,
in which case $y \notin T^n \supset M^n$ again, or $f_n(T^n_x)=T^n_{f_n(x)}
\cap \partial T^n \neq \emptyset$, in which case $y \notin M^n$ by
construction.

By taking iterates, we may assume that $\depth(T^n_x)=0$, and $M^n_x=M^n$.
If $M^n \cap T^n_y \neq \emptyset$, then by construction $T^n_y \cap
\partial T^n \neq \emptyset$, so $T^n_y \subset T^n$.
\end{pf}

\begin{lemma} \label {Mn}

With respect to the hyperbolic metric in $\V^n \setminus J^n$,
$\partial M^n$ and $\partial T^n$ are at bounded distance in the Hausdorff
sense.  In particular, $M^n$ contains a definite neighborhood of $0$ in
scale $\diam(J^n)$.

\end{lemma}

\begin{pf}

Since $\partial T^n$ has uniformly bounded diameter in $\V^n \setminus
\OO(f)$, it is enough to show that if $T^n_x \cap \partial T^n \neq
\emptyset$ then $T^n_x$ has bounded hyperbolic diameter in $\V^n \setminus
\OO(f)$.  By the Schwarz Lemma, the hyperbolic diameter of $T^n_x$
inside $\Gamma^n \setminus J^n$ is smaller than the hyperbolic
diameter of $T^n_x$ in $\Gamma^n_x$ which is bounded.

The last statement follows from the behavior of the hyperbolic metric near
$0$: at every $x \in \V^n$, it is at least of order of the inverse of the
Euclidean distance to $0$.  This follows from $J^n$ connected and
$\diam(J^n) \asymp \diam(\V^n)$.
\end{pf}

Let $S^n$ be the connected component of $0$ in $\cup f^{-k}(X)$.

\begin{lemma}

$S^n \subset T^n$.

\end{lemma}

\begin{pf}

$S^n$ is an increasing union of a finite number
$M^n_{x_1}$,...,$M^n_{x_l}$ of copies of $M^n$.  One proves by
induction on $l$ that such finite union is contained in $T^n_{x_i}$ such
that $\depth(M^n_{x_i})$ is minimal.
\end{pf}

\begin{lemma}

The hyperbolic diameter of $\partial S^n$ in $\V^n \setminus \OO(f)$ is
uniformly bounded.

\end{lemma}

\begin{pf}

We have $\partial S^n \subset \overline T^n \setminus M^n$ which is a
compact subset of $\V^n \setminus \OO(f)$.  The result now
follows from Lemma~\ref {Mn}, assertion (3) of Lemma~\ref {boundedness},
and hyperbolic geometry.
\end{pf}

We define again $Q^0=S^0$, $Q^{n+1}=f_n^{-1}(Q^n) \cap S^{n+1}$ and let
$V^n$ be the filled $Q^n$.  Conditions (C1-3) and (G1-2) are easily checked
as before, as well as additional properties (A1-3).  Property (A4) no longer
holds, but can be substituted by
\begin{enumerate}
\item [(A4')] There exists a finite set $F \subset \cup f_n^{-k}(\beta_n)$
such that $V^n$ is a neighborhood of $J^n \setminus F$.
\end{enumerate}
The construction can be
also modified for periodic points of renormalization without any differences
to the previous considerations.

\begin{rem}

The only parameters determining the constants in this construction are the
unbranched {\it a priori} bounds and the combinatorics.  Actually, one can
see that a lower bound on the unbranched {\it a priori} bounds implies an
upper bound on the period of any immediate renormalizations.  As a
consequence, the combinatorics are only used to get a lower bound on $\nu$
(the spacing between non-touching $J^n_i$).  As such, they are uniform for
all real quadratic maps (or real unicritical polynomials of a fixed degree),
as well as for certain classes of quadratic maps of
``high combinatorial type'' considered in \cite {puzzle}.

\end{rem}

\comm{
\begin{lemma}

The boundary of $T^n$ has bounded hyperbolic diameter in $\V^n \setminus
\OO(f)$.  Moreover, if $T^n_x \cap \partial T^n=\emptyset$ then $T^n_x$ has
bounded hyperbolic diameter in $\V^n \setminus \{J^n\}$.

\end{lemma}

\begin{pf}

The first property is immediate from the unbranched apriori bounds.

For the second, we notice that $T^n_x \cap \partial T^n \neq \emptyset$
implies that $\Gamma^n_x \cap J^n=\emptyset$.  Indeed, the only $J^n_i$
which intersect $\Gamma^n$ are of the form $f_{n-1}^j(J^n)$, so if
$\Gamma^n_x \cap J^n \neq \emptyset$ then $J^n_x$ touches $J^n$.  This
implies that $f^k:J^n_x \to J^n$, $k=\depth(J^n_x)$,
is an iterate of $f_{n-1}$.  But by construction of $T^n$ (which is nice for
$f_{n-1}$), this implies that $T^n_x \cap \partial T^n=\emptyset$.

Since $\Gamma^n_x \setminus T^n_x$ has definite moduli and $\partial T^n_x$
is very close to $J^n$, it follows that the diameter of $J^n_x$ is much
smaller than the diameter of $J^n$ and in particular $\Gamma^n_x \subset
\Gamma^n$.  By the Schwarz Lemma, the hyperbolic diameter of $T^n_x$
inside $\V^n \setminus J^n_x$ is smaller than the hyperbolic
diameter of $T^n_x$ in $\Gamma^n_x$ which is bounded.
\end{pf}

if $T^n_x \cap T^n=\emptyset$
implies that there exists $s \leq k=\depth(T^n_x)$ such that $f^s(T^n_x)
\cap 

$J^n_x$ is not a copy of $J^n$ obtained by
taking preimages under $f_{n-1}$.   $\Gamma^n_x \cap J^n=\emptyset$.
Moreover,

\begin{enumerate}
\item $M^n$ is simply connected,
\item $f_n^{-1}(M^n) \subset M^n$,
\item $J^n \setminus \{\beta_n\} \subset M^n$,
\item $M^n$ contains a disk or radius comparable with $\diam(J^n)$ centered
around $0$,
\item If $T^n_x \cap M^n \neq \emptyset$ then $T^n_x \subset T^n$.
\end{enumerate}

\end{lemma}

Let us define copies of $M^n$ in the natural way.

\begin{lemma}

The connected component of $0$ in $\cup f^{-k}(M^n)$ is contained in $T^n$.

\end{lemma}

In the immediately renormalizable case, the copies of Julia sets are not
necessarily disjoint.  Indeed, if $f_{n-1}$ is immediately renormalizable of
period $q$ then $\beta_n$ (the non-dividing fixed point of $f_n$) belongs to
$f_{n-1}^i(J^n)$, $0 \leq i \leq q-1$, and we actually have
$f_{n-1}^i(J^n) \cap f_{n-1}^j(J^n)=\{\beta_n\}$, $0 \leq i<j \leq q-1$.
Moreover, If $J^n_x$ is a copy of $J^n$ distinct from one of the
$f^i_{n-1}(J^n)$ then $J^n_x \cap J^n=\emptyset$.

If $J^n_x$ is a copy of $J^n$ of depth $k$, we let
$\alpha(J^n(x))=(f^k|J^n)^{-1}(\alpha_n)$.  If $J^n_x \cap J^n_y \neq
\emptyset$ then we have necessarily $J^n_x \cap
J^n_y=\alpha(J^n_x)=\alpha(J^n_y)$.

Let $JB^n=\cup_{0 \leq i \leq q-1} f_{n-1}^i(J^n)$.  Then
$$
\diam(JB^n) \asymp \diam (J^n),
$$
the constant depending on $q$.  If $x \in \cup
f^{-k}(J^n)$, we will denote by $JB^n_x$ its connected component (which is
the union of $q$ copies of $J^n$), which will be called a copy of $JB^n$.

The depth of $JB^n_x$ is by definition
the minimum of the depths of the copies of $J^n$ that it intersects.  If
$\depth(JB^n_x)=k$ then $f^k:JB^n_x \to JB^n$ extends to a univalent map
onto $\cup_{0 \leq i \leq q-1} f_{n-1}^{-i} \V^n$, which is a definite
neighborhood of $JB^n$.  Thus
$$
\diam(JB^n_x) \asymp \diam(J^n_x).
$$

\begin{lemma}

For every $\epsilon>0$, there exists $\delta>0$ such that if the
$\delta$-neighborhoods of two copies $J^n_x$ and $J^n_y$ of $J^n$
intersect then $\diam(J^n_y) <
\epsilon \diam(J^n_x) $, provided $\depth (J^n_y) \geq \depth (J^n_x)$.

\end{lemma}

Let $JB^{n-1}=\cup_{0 \leq i \leq q-1} f_{n-1}^i(J^n)$.  The connected
components of $If $f^k(x) \in
JB^n$ is the first landing moment of some $x \in U$ to $JB^n$ then there
exist

\begin{lemma}

Assume that $f_{n-1}$ is immediately renormalizable (with renormalization
period $q$).  Then the connected

\comm{
In the immediately renormalizable case, Lemma~\ref {copies} is false:
indeed, if $f_n$ is immediately renormalizable then
different copies of $J^{n+1}$ touch at the preorbit of $\beta_n$.

Let $W^n(\rho)$ be defined as in Lemma~\ref {nrho}.

\begin{lemma}

For every $\epsilon>0$, there exists $\delta>0$ such that if the
$\delta$-neighborhoods of two non-touching copies
$J^n_x$ and $J^n_y$ of $J^n$ intersect
then $\diam(J^n_y) <  \epsilon \diam(J^n_x) $, provided
$\depth (J^n_y) \geq \depth (J^n_x)$.

\end{lemma}

Let $B^n$ be the union of all copies of $J^n$ which intersect $J^n$.  Thus
$B^{n+1}=J^{n+1}$ if and only if $f_n$ is not immediately renormalizable. 
Otherwise $B^{n+1}$ is the union of $p$ copies of $J^{n+1}$ where $p$ is the
period of renormalization of $f_n$.

\begin{lemma}

For every $\kappa>0$, there exists $\rho>0$ such that if
$x_1,...,x_n \in \cup f^{-k}(J^n)$, $\cup W^n_{x_i}(\rho)$
is connected then $\cup W^n_{x_i}(\rho)$ is contained in the
$\kappa$-neighborhood of $B^n_{x_j}$, where $\depth(B^n_{x_j)})$ is
minimal.

\end{lemma}

Let us cut $W^{n+1}(\rho)$ by the external rays landing at $\beta_n$
corresponding to an efficient straightening of $f_n$, and let $Z^{n+1}$ be
the connected component of $0$.  We define copies of $Z^{n+1}$ in the natural
way.  Let us call admissible the copies of $Z^{n+1}$ that do not intersect
the postcritical set.
A cluster is a connected component of admissible copies of $Z^{n+1}$.
Let us call the span of a cluster the set of $j \geq 0$ such that
$f_{n+1}^{-j}(\V^{n+1} \setminus \u^{n+1})$ intersects
the cluster.  It is an interval.

\begin{lemma}

The span of a cluster has uniformly bounded length.

\end{lemma}

\begin{pf}

Let $[j,k]$ be the span of the cluster.  Taking the image by $f_{n+1}^j$, we
may assume that $j=0$.  Let $k$ is large then the cluster intersects
$W^n(\rho)$.  This implies 

We define $S^{n+1}$ as the connected component of $0$ of
$Z^{n+1}$ minus the closure of the union of all
clusters that do not intersect the boundary of $Z^{n+1}$.  That it is
non-empty follows from the previous lemma:

\begin{lemma}

$S^{n+1}$ contains a definite neighborhood of points in $J^{n+1}$ which are
away from $\beta_n$.

\end{lemma}

Let $M^n=M^n(\rho)$ be the union of the
copies of $W^n(\rho)$ corresponding to the connected components of $B^n$.
We call copies of $B^n$ the connected
components of the union of all copies of $J^n$.  To each copy of $B^n$ we
can associate a copy of $M^n$.

Choose now $\kappa$ small, $\rho$ as in the previous lemma and let $T^n$ be
the connected component of $0$ in $\cup f^{-k} W^n(\rho)$.

If $f_n$ is not immediately renormalizable we set $S^{n+1}=T^{n+1}$.
Otherwise, let us cut $T^n$ by the external rays landing at $\beta_n$
corresponding to an efficient straightening of $f_n$ and let $Z^{n+1}$ be
the component of $0$.

In the immediately renormalizable case, Lemma~\ref {copies} is false:
indeed, if $f_n$ is immediately renormalizable then
different copies of $J^{n+1}$ touch at the preorbit of $\beta_n$.

Let $W^n(\rho)$ be defined as in Lemma~\ref {nrho}.

\begin{lemma}

For every $\epsilon>0$, there exists $\delta>0$ such that if the
$\delta$-neighborhoods of two non-touching copies
$J^n_x$ and $J^n_y$ of $J^n$ intersect
then $\diam(J^n_y) <  \epsilon \diam(J^n_x) $, provided
$\depth (J^n_y) \geq \depth (J^n_x)$.

\end{lemma}

Let $B^n$ be the union of all copies of $J^n$ which intersect $J^n$.  Thus
$B^{n+1}=J^{n+1}$ if and only if $f_n$ is not immediately renormalizable. 
Otherwise $B^{n+1}$ is the union of $p$ copies of $J^{n+1}$ where $p$ is the
period of renormalization of $f_n$.

\begin{lemma}

For every $\kappa>0$, there exists $\rho>0$ such that if
$x_1,...,x_n \in \cup f^{-k}(J^n)$, $\cup W^n_{x_i}(\rho)$
is connected then $\cup W^n_{x_i}(\rho)$ is contained in the
$\kappa$-neighborhood of $B^n_{x_j}$, where $\depth(B^n_{x_j)})$ is
minimal.

\end{lemma}

\begin{lemma}

For every $\kappa>0$ there exists $\rho>0$ such that if
$S$ is a connected component of $\cup f^{-k} W^n(\rho)$ containing 0
then  $S$ is contained in a $\kappa$-neighborhood of $B^n$. 

\end{lemma}

Let us cut $W^{n+1}(\rho)$ by the external rays landing at $\beta_n$
corresponding to an efficient straightening of $f_n$, and let $Z^{n+1}$ be
the connected component of $0$.  We define copies of $Z^n$ in the natural
way.  Let us call a cluster a connected component of the union of copies of
$Z^n$ which do not intersect the postcritical set.

Let us call the span of a cluster the number of $k$ for which
$f_{n+1}^{-k}(\V^n \setminus \u^n$ intersects the cluster.

\begin{lemma}

The span is uniformly bounded.

\end{lemma}

We define $S^{n+1}$ as the connected component of $0$ of
$Z^{n+1}$ minus the closure of the union of all
clusters that do not intersect the boundary of $Z^{n+1}$.  That it is
non-empty follows from the previous lemma:

\begin{lemma}

$S^{n+1}$ contains a definite neighborhood of points in $J^{n+1}$ which are
away from $\beta_n$.

\end{lemma}

Let $M^n=M^n(\rho)$ be the union of the
copies of $W^n(\rho)$ corresponding to the connected components of $B^n$.
We call copies of $B^n$ the connected
components of the union of all copies of $J^n$.  To each copy of $B^n$ we
can associate a copy of $M^n$.

Choose now $\kappa$ small, $\rho$ as in the previous lemma and let $T^n$ be
the connected component of $0$ in $\cup f^{-k} W^n(\rho)$.

If $f_n$ is not immediately renormalizable we set $S^{n+1}=T^{n+1}$.
Otherwise, let us cut $T^n$ by the external rays landing at $\beta_n$
corresponding to an efficient straightening of $f_n$ and let $Z^{n+1}$ be
the component of $0$.

Thus Lemma~\ref {copies} must be replaced.  Let $J^{n,\gamma}$ be the
truncation of $J^n$ obtained by
cutting a $\gamma$-neighborhood of $\beta_n$ out of $J^n$.  We define
truncation of copies in the obvious way.

\begin{lemma}

Let $\gamma>0$.  For every $\epsilon>0$, there exists $\delta>0$ such that
if the $\delta$-neighborhoods of two copies $J^{n,\gamma}_x$ and
$J^{n,\gamma}_y$ of $J^{n,\gamma}$ intersect
then either $\min \{\diam(J^{n,\gamma}_x),\diam(J^{n,\gamma}_y)\} \leq
\epsilon \max \{\diam(J^{n,\gamma}_x),\diam(J^{n,\gamma}_y)\}$.

\end{lemma}

\begin{lemma}

Let $\gamma>0$.  For every $\rho>0$, there exists $\tau>0$ and a domain
$W^n(\rho) \subset U^n$ containing $J^n \setminus \{\beta_n\}$, such that
$\overline {f_n^{-1}(W^n(\rho))}
\subset W^n(\rho)$ and $W^n(\rho)$ contains a $\tau$-neighborhood of
$J^{n,\gamma}$, and that for any univalent map $g:D \to \V^n$,
$g^{-1}(W^n(\rho))$ is contained in a $\rho$-neighborhood of $g^{-1}(J^n)$.

\end{lemma}

Choose now $\kappa$ small, $\rho$ as in the previous lemma and let $T^n$ be
the connected component of $0$ in $\cup f^{-k} W^n(\rho)$.

If $f_n$ is not immediately renormalizable we set $S^{n+1}=T^{n+1}$.
Otherwise, let us cut $T^n$ by the external rays landing at $\beta_n$
corresponding to an efficient straightening of $f_n$ and let $Z^{n+1}$ be
the component of $0$.
}
}
}
\section{Conformal measures on towers}\label{towers-sec}

\subsection{Tower limits of conformal measures}
In this section, 
let $f:\u\ra \V$ be a Feigenbaum map with stationary combinatorics,
and let $f_n: \u^n\ra \V^n$ be its pre-renormalizations
with domains satisfying Lemma \ref{shapes}.  Let $\A^n=\V^n \setminus \u^n$. 
%In the non-satellite case, the domains $U_n$ and $V_n$ 
%have been selected in such a way 
%$ that $\bar U_n \cap \bar V_{n+1}=\emptyset$, and the fundamental annuli $\bar A_n$ don't intersect the postcritical set.
Consider rescalings  
$$
   f_m^{(n)} =  T_n \circ f_{n+m}\circ T_n^{-1}  : \u^m_{(n)} \ra \V^m_{(n)},
\quad n=0,1,2,\dots; \ m=-n, -n+1, \dots, 
$$
that normalize  the maps $f_0^{(n)}  = T_n \circ f_n \circ T_n^{-1}$.   

Let $\mu$ be a $\de$-conformal measure on $J(f)$.
To study its local geometry, 
 let us push it forward by dilations $T_n$, $n\in \N$, and normalize to be 1 on
$\V^0_{(n)}$, 
$\mu_n = (T_n)_*(\mu)/ \kappa_n$. 

\begin{lem}
  The measures $\mu_n$ are uniformly bounded on compact subsets of $\C$.
\end{lem}
 
\begin{pf}
Take some $m,k\in \Z_+$, and consider domains $\V^m$ and $\u^{m+k}$. We will show that
\begin{equation}\label{comparison3}
   \mu(\V^m) \leq C_k \mu(\u^{m+k}),
\end{equation}
which after rescalings implies the assertion. 
Moreover, taking $R^mf$, we see that it is enough to prove (\ref{comparison3}) for $m=0$ with the constant $C_k$ depending
only on $\mod(f)$. 

Since the straightening of $f$ is $K$-qc where $K$ depends only on $k$ and  $\mod(f)$,  
there exists an $N$ depending only on $\mod (f)$ such that $f^N \u^k=\V$. 
Moreover, $| D f^N(z)  |\leq C$, where $C$ depends on the same data.
By $\de$-conformality of $\mu$, we conclude that 
$ \mu (\V) \leq C^\de \mu(\u^k)$.
\end{pf}  
   
Hence the family of measures $\mu_m$ is precompact in the weak$^*$
topology on compact subsets.  Let $\bmu$ be a limit of these measures. 
This measure is Radon (i.e.,  it assigns finite mass to compact sets)
and    $\de$-conformal  for the whole tower $\bar f=\{f_m\}_{m\in \Z}$. 

%Let us select quadratic-like maps $f_m: \u^m\ra V^m$ representing the tower $\bar f$
%and satisfying the  following properties:
%\begin{itemize}
%\item $U^0$ and $V^0$ are smooth Jordan curves;
%\item $A^0\equiv U^0\sm V^0$ is compactly contained in $\C\sm \om (\bar f)$;
%\item $A^m\equiv U^m\sm V^m= T^{-m} A^0$.
%\end{itemize}
%

\subsection{Ergodicity}
 A set $X\subset \C$ is called invariant
with respect to the tower $\bar f$ if $f_m(X\cap \u^m)\subset X$ for all $m\in Z$. 
 A measure $\bar \mu$ is called ergodic with respect to the tower if there is no decomposition
$\C=X\cup Y$ into two disjoint measurable sets of positive measure invariant with respect to
the tower.

\begin{thm}\label{tower ergodicity}
Let $\bar\mu$ be a conformal measure of the tower $\bar f$ which is fully supported on
either the Julia set $J(\bar f)$, or on its complement $\C\sm J(\bar f)$. Then
$\bar\mu$ is ergodic. 
\end{thm}

\begin{pf}  If $\bar\mu$ is supported on the Julia set $J(\bar f)$, then 
Theorem \ref{ergodicity} implies that it  is ergodic. 

Assume that $\bar \mu$ is supported on $\C\sm J(\bar f)$. Assume that there is a 
decomposition $\C\sm J(\bar f)=X\cup Y$ into two invariant measurable sets with positive
measure. Without loss of generality we can assume that the inequality
\begin{equation}\label{choice}
\bar\mu(Y\cap \A^m)   \geq    \bar\mu (X\cap \A^m)
\end{equation}
holds for infinitely many levels $m\in -\N$.

Let us cut each $\A^m$ by a  smooth arc $\gamma$ of zero $\bar\mu$-measure to obtain a 
(topological) rectangle $\Delta^m$.   Then there exists a compact subset
$S\subset X\cap \u^0$ of positive measure whose forward orbit does not intersect the cuts
$\gamma_m$. It follows that for each $m\in -\N$, this set can be covered by  disjoint preimages
$\Pi^m_i$ of $\De^m$ under iterates  of $f_m$.

By  \cite[Prop. 6.9]{McM-towers}, $\diam \Pi^m_i\to 0$ as $m\to -\infty$.
Hence the unions $\Pi^m\equiv\cup_i\Pi^m_i$ shrink to $S$ as $m\to -\infty$. It follows that
$\bar\mu(\Pi^m\sm S)\to 0$. Hence there exists a rectangle  $\Pi^m_{i(m)}$
such that $\dens(S|\Pi^m_{i(m)})\to 1$ as $m\to -\infty$.

On the other hand, since the rectangles $\Pi^m_i$ are mapped onto the $\De^m$ with bounded
distortion, $\dens(S|\Pi^m_i)\leq 1-\eps<1$ for all $m$ satisfying (\ref{choice}). 
This contradiction proves the result.
\end{pf} 

\begin{cor}\label{uniqueness}
 Given an exponent $\de$, a tower can have at most one normalized $\de$-conformal measure
supported on the Julia set $J(\bar f)$, and at most one such a measure supported on the
complement $\C\sm J(\bar f)$.
\end{cor}

Thus,  any tower conformal measure $\bar \mu$ may have at most two ergodic components,
$\bmu| J(\bar f)$ and $\bmu| \C\sm J(\bar f)$.

\subsection{Scaling}

In what follows, $\rho\in (0,1)$ will stand for the scaling factor of the renormalization fixed point. 

\begin{lem}[Scaling covariance]\label{scaling covariance}
Let $\bar\mu$ be a conformal measure for the tower $\bar f$
 fully supported either on $J(\bar f)$,
or on $\C\sm J(\bar f)$.
Then there exists a $\si=\si(\bmu)>0$ such that $T_*\bmu=\rho^\si\bmu$,  and for all $r>0$,
$\bar\mu(\D_r) \asymp r^\si$.
\end{lem}

\begin{pf} 
Since the Feigenbaum dilation $T z\mapsto \rho z$ conjugates $f_m$ to $f_{m-1}$, the measure
$T_*\bar\mu$ is also conformal with the same exponent.
 By Corollary \ref{uniqueness}, $T_*\bar\mu=\kappa\bar\mu$ with some $\kappa>0$.
  Hence $\bar\mu(\D_r)=\kappa\bar\mu(\D_{\rho r})$ for all $r>0$. Thus
$$
   \bar\mu(\D_{\rho^m})=  \kappa^{-m} \bar\mu(\D),\quad m=\pm 1, \pm 2, \dots,
$$
and the statement follows with $\si=-(\log\kappa/\log\rho)$. 
\end{pf}

If $\bmu$ is a tower conformal measure with two ergodic components,
then let 
$$\sigma_-\equiv \sigma_-(\bmu)=\max\{ \sigma(\bmu| J(\bar f)), \si(\bmu | \C\sm J(\bar f)) \};$$
$$\sigma_+\equiv \sigma_+(\bmu)=
\min \{ \sigma(\bmu| J(\bar f)), \si(\bmu | \C\sm J(\bar f)) \}.$$
In the case of one ergodic component we set $\si_+=\si_-=\si(\bmu)$.

\begin{cor}
For any tower $\de$-conformal measure $\bmu$,
  $\bmu(\D_r)\asymp r^{\si_+}$ for $r\leq 1$ and
   $\bmu(\D_r)\asymp r^{\si_-}$ for $r\geq 1$, where the  constants are independent 
  of $\bmu$.   
\end{cor}

Let $\bmu_+$ be the tower conformal measure on the Julia set (if it exists), 
and $\bmu_-$ be such a measure on the complement (if exists). 
Both measures are normalized to be 1 on $\V$. 

\begin{cor}\label{behave at 0}
Let $\mu$ be a conformal measure of $f$. Then 
$$
     \si_+\leq \lim_{r\to 0}\inf  {\frac {\log \mu(\D_r)} {\log r}}\leq 
     \lim_{r\to 0}\sup {\frac {\log \mu(\D_r)} {\log r}} \leq  \si_-.
$$  
\end{cor}

\begin{pf} 
  The limit set of the measures $\mu_n$ consists of convex combinations 
$t\bmu_+ + (1-t) \bmu_-$.  For any such combination  $\bmu$, we have:
$$ 
     \rho^{\si_-} \bmu (\D_r)  \leq   \bmu(\D_{\rho r}) \leq \rho^{\si_+} \bmu (\D_r)
$$ 
for any $r>0$. Hence 
$$
          \rho^{\si_- + \eps}   \mu (\D_r)  \leq   \mu(\D_{\rho r})   \leq  \rho^{\si_+ - \eps } \bmu (\D_r),
$$
provided $r$ is sufficiently small (depending on $\eps>0$),  which implies the assertion.
\end{pf}

\begin{lem}\label{aux}
 If a tower conformal measure $\bmu$ is not ergodic, then
$$\si_+=\si(\bmu | J(\bar f))< \si (\bmu | \C\sm J(\bar f)) = \si_- .$$
\end{lem}

\begin{pf}
Let $X=\C\sm J(\bar f)$, $Y=J(\bar f)$.  Assume that
$$\si_+ = \si(\bmu | X ) \leq \si(\bmu | Y  )= \si_- .$$
 Then by Lemma \ref{scaling covariance},
$$
   \bmu (X\cap \A^m) = \bmu(X\cap \A^0) \rho^{m\si_+}\quad  {\mathrm and} \quad
   \bmu (Y\cap \A^m) = \bmu(Y\cap \A^0)  \rho^{m\si_-}. 
$$
Hence there exists a $c>0$ such that $\bmu(Y\cap \A^m)\geq c \bmu(X\cap \A^m)$ for all $m\in -\N$. 
This leads to a contradiction in the same way as (\ref{choice}) led to a contradiction in the
proof of Theorem \ref{tower ergodicity}. 
\end{pf}

\begin{lem}\label{de and si} 
Assume that $\area(J(f))=0$. Let $\bmu$ be a tower $\de$-conformal measure.
\begin{itemize}
\item [(i)] If  $\de<  2$ then $\si_-(\bmu) < \de$.
\item [(ii)] If $\de =2$ then $\si_-(\bmu)  \leq 2$. 
  Moreover, $\si_+(\bmu) < 2$  unless $\bmu$ is the Lebesgue measure.%
 \footnote{Thus, if $\bmu$ is supported on $J(\bar f)$ then $\si(\bmu) < 2$. }   
\end{itemize}   
\end{lem}

\begin{pf} 
Since $\area(J(f))=0$, there exists a compact set $S\subset \D\sm J(\bar f)$ of positive
Lebesgue measure.  

Let $\De^m\subset \A^m$ be the topological rectangles obtained by cutting the annuli
$\A^m$
by arcs of zero $\bar\mu$-measure (as in the proof
of Theorem \ref{tower ergodicity}).  
Take a point  $z\in S$, and consider a sequence of moments
$n_m=n_m(z)\to \infty$ 
such that $\zeta_m\equiv f^{n_m} z\in \De^m$, $m\leq 0$. Pulling this rectangle
back by $f^{n_m}$, we obtain a rectangle $\Pi^m=\Pi^m(z)$ containing $z$. 
 
\msk
 (i) Let $\de < 2$. 
 Assume $\si_- \geq  \de$. Then 
 \begin{equation}\label{eq1}
    \bar\mu(\De^m) \geq c (\diam \De^m)^{\si_-} \geq  c (\diam \De^m)^\de, \quad m\in -\N.        
\end{equation}
Let $g_m=f^{n_m}: \Pi^m\ra \De^m$.
By the dynamical $\de$-covariance and bounded distortion, 
\begin{equation}\label{eq2}
    \bar\mu(\Pi^m)\geq a{\frac {\bmu(\De^m)} {|D g_m(z)|^\de}}\geq 
     ac {\frac {(\diam\De^m)^\de} {|Dg_m(z)|^\de}}\asymp  (\diam \Pi^m)^\de.
\end{equation}
% \end{equation}
Hence 
\begin{equation}\label{eq3}
\mbox{ $\bar\mu(\Pi^m)/\area(\Pi^m)\to \infty$ as $m\to -\infty$.}
\end{equation} 
But for a given $m$, different  rectangles $\Pi^m(z)$, $z\in S$, are disjoint
and cover the whole set $S$.
% Moreover, by \cite{McM-towers},  $\diam \Pi^m\to 0$ as $m\to \infty$.
It follows that $\bar\mu(S)=\infty$ - contradiction. 

\msk
 (ii) Let  $\de  = 2$. 
 Note that the Lebesgue measure is a 2-conformal measure supported on $\C\ssm J(\bar f)$.
By ergodicity  (Theorem \ref{tower ergodicity}),
 $\bmu\, |\,  \C\sm J(\bar f)$ is proportional to it.
 Let  $\bnu= \bmu+\area$. By Lemma \ref{aux}, 
$\si_-(\bnu) = \si (\area) =2$ and $\si_+(\bnu)<2$ unless $\bnu$
 (and hence $\bmu$) is Lebesgue.    
\end{pf}

%\begin{cor}\label{r to si} Let $\mu$ be a $\de$-conformal measure on $J$, $\de\in (1,2]$.
%Then  $ \mu(\D_r) =  r^\si \phi(r),$
%where $\phi(r) \to \infty$ and $\log \phi(r) / \log r \to 0$  as $r\to 0$. 
%Moreover,  $\si\leq \de$, and $\si<\de$  for $\de<2$.
%\end{cor}

\begin{lem}\label{scaling of mu}
  Assume $\area(J)=0$ but $\HD(J)=2$. Let $\mu$ be a $2$-conformal measure on $J$. 
Then
$$
    \frac{\log \mu(\D_r)}{\log r} \to 2\quad \text{as}\quad r\to 0.
$$
\end{lem}

%\begin{cor}\label{area exp decay}
%Assume  $\area J=0$. 
%Let $f_n: U^n\ra V^n$ be an $n$-fold pre-renormalization of $f$.
%% with nice $V^n$;
%$X_n =\cup_{k=0}^\infty f^{-k} V^n$. Then
% $$
%     \eta_n\equiv \area(X_n) \asymp \rho^{n(2-\si)} / \phi (\rho^n) ,
%$$
% where $\rho\in (0,1)$ is the scaling factor and $\si$ is the exponent associated with the $2$-conformal measure $\mu$
%on $J$\footnote{which exists by  Proposition \ref{spec}}.    
%\end{cor} 

\begin{pf}
By Corollary \ref{behave at 0} and Lemma \ref{de and si},
$$
    \lim_{r\to 0} \sup \frac{\log \mu(\D_r)}{\log r} \leq 2.
$$
Assume
$$
\lim_{r\to 0} \inf \frac{\log \mu(\D_r)}{\log r} <  2,
$$    
so that $\mu(\D_{r_k}) \geq r_k^\si$ for some $\si<2$ and a sequence $r_k\to 0$.
Then
$$
   \mu(\V^{n_k}) \geq (\diam \V^{n_k})^\si  \asymp \rho^{\si n_k }
$$
for some sequence of $n_k\to \infty$.

Let $X_n =\cup_{k=0}^\infty f^{-k} \V^n$.
Let us consider the landing map $\La_n: X_n \ra \V^n$ to $\V^n$
and the associated Poincar\'e series  $ \Xi_2 (\La_n, z)$, $z\in \V^n$.
% By Proposition \ref{spec}, there exists a  2-conformal measure $\mu\equiv \mu_2$ on $J$.
For $z\in \V^n$, we have: 
$$
       \mu(\V^n)\,  \Xi_2(\La_n, z)\asymp   \mu (X_n) =1,   
$$
and  hence 
$$
    \Xi_2(\La_{n_k}, z)\leq  \rho^{-\si {n_k}} .
$$
   But then
$$
  \eta_{n_k}\equiv \area(X_{n_k}) \asymp \area (\V^{n_k}) \,  \Xi_2(\La_{n_k}, z)\leq \rho^{2 n_k}\rho^{-\si {n_k}},  
$$ 
contradicting  Theorem \ref{balanced case}.
\end{pf}

Let $\mu$ be a $2$-conformal measure on the Julia set $J$. 
We say that the Julia set is {\it measure-theoretic hairy} at the critical point if 
$$
  {\frac {\pi} {\mu(\D_r)}}\, \mu\circ r  \to \area \quad \text{as}\quad r\to 0.
$$
Here convergence is understood in the weak$^*$-topology on compact sets.
%In other words, for any continuous function $\phi$ with compact

\begin{prop}\label{hairy}
  If $\HD(J) = 2$ then  $J$ is measure-theoretic hairy.
\end{prop} 
 
\begin{pf}
  If $\area(J)>0$ then the assertion is obvious since the critical point is the density point of $J$.
So, assume $\area(J)=0$.
Then any limit measure $\bar\mu$ on the tower $\bar f$ is 2-covariant not only by the
tower dynamics but by the $\rho$-scaling as well. 
Hence $\bar\mu(\D_r)\asymp r^2$.
Spreading it around using distortion estimates, we see that $\bar\mu(\D_r(z))\asymp r^2$ 
near any point $z\in \C$.
Hence $\bar\mu$  is equivalent to the Lebesgue measure. 
But the only absolutely continuous  $2$-conformal measure is Lebesgue.
\end{pf}

Let us consider a function $\phi: \R_+\ra \R_+$ such that $\phi(r)\to \infty$ as $r\to 0$.
  Let us say that it  is {\it slowly varying} (near the origin)
if for any $\kappa<1$, $\phi(\kappa r)/\phi(r)\to 1$ as $r\to 0$. Clearly, such a function has 
sub-polynomial growth:
$$
    { \frac {\log \phi(r)} {\log r}}  \to 0 \quad \text{as} \quad  r\to 0.
$$  

\begin{cor}
   Assume $\area(J) = 0$ but $\HD(J) =2$. Let $\mu$ be a $2$-conformal measure on $J$.
Then $\mu(\D_r) = r^2 \phi(r)$, where $\phi$ is a slowly varying (growing to $\infty$) function.
\end{cor} 

\begin{pf}
  Let us show that $\phi(r)\to \infty$. Otherwise there would be a sequence of levels $n_k\to \infty$
such that $\mu(\V^{n_k} ) \leq C \area (\V^{n_k})$. Hence for $z\in \V^{n_k}$, 
$$
 \mu(X_{n_k}) \asymp \mu(\V^{n_k})\, \Xi_2 (z) \leq C \area (\V^{n_k} )\, \Xi_2 (z) \asymp \area(X_{n_k}).
$$
(Here $X_n$ is as in Lemma \ref{scaling of mu}.)
But this is impossible since $\mu(X_{n_k})=1$ while   $\area(X_{n_k})\to 0$ as $k\to \infty$.

\ssk
   Furthermore,   by Proposition \ref{hairy},
$$
    {\frac {\mu (\D_{\kappa r} )} {\mu (\D_r)}} \to \kappa^2 \quad \text{as}\quad r\to 0,
$$ 
which implies that $\phi$ is slowly varying.
\end{pf}

In fact, we have:

\begin{prop}
    Under the above circumstances, $\phi(r)\asymp \log \frac{1}{r} $ for sufficiently small $r$. 
\end{prop}

\begin{pf}
  As in the proof of Lemma \ref{scaling of mu}, we have:
$$
       \Xi_2(\La_n, z)\asymp  \frac{1}{\mu(\V^n)} \asymp \frac{1} {\rho^{2n} \phi(\rho^n})   
$$
and
$$
  \eta_n \equiv \area(X_n) \asymp \area (\V^n) \,  \Xi_2(\La_{n_k}, z) \asymp \frac{\rho^{2 n } } { \rho^{2n } \phi(\rho^n)}= 
               \frac{1}{\phi(\rho^n)}.  
$$ 
Since by  Theorem \ref{balanced case}, $ \eta_n  \asymp 1/n$,
we obtain the right asymptotics for $r=\rho^n$.
Since $\phi$ is slowly varying, we can interpolate it to the intermediate scales.
\end{pf}

\end{document}